\documentclass{amsart}
\usepackage{amssymb}
\usepackage{amsmath}
\usepackage{amsfonts}

\newtheorem{theorem}{Theorem}
\newtheorem{lemma}{Lemma}[section]
\theoremstyle{definition}
\newtheorem{definition}[lemma]{Definition}
\newtheorem{example}[lemma]{Example}

\theoremstyle{remark}
\newtheorem{remark}[lemma]{Remark}
\numberwithin{equation}{section}
\theoremstyle{plain}

\newtheorem{conjecture}[lemma]{Conjecture}
\newtheorem{corollary}[lemma]{Corollary}

\newtheorem{problem}[lemma]{Problem}
\newtheorem{proposition}[lemma]{Proposition}

\input{tcilatex}

\begin{document}
\title[Uniqueness for the Signature]{Uniqueness for the Signature of a Path
of Bounded Variation and the Reduced Path Group}
\author{B.M. Hambly}
\address[B.M. Hambly and Terry J. Lyons]{Mathematical Institute\\
Oxford University\\
24-29 St. Giles\\
Oxford OX1 3LB\\
England}
\email[B.M. Hambly]{hambly@maths.ox.ac.uk}
\author{Terry J. Lyons}
\email[Terry J. Lyons]{tlyons@maths.ox.ac.uk}
\thanks{The authors gratefully acknowledge EPSRC support: GR/R29628/01,
GR/S18526/01}
\date{July, 2006}

\begin{abstract}
We introduce the notions of
tree-like path and tree-like equivalence between paths and prove that the
latter is an equivalence relation for paths of finite length. We show
that the equivalence classes form a group with some
similarity to a free group, and that in each class there is one special tree
reduced path. The set of these paths is the Reduced Path Group. It is a
continuous analogue to the group of reduced words.
The signature of the path is a power series whose coefficients
are definite iterated integrals of the path. 
We identify the paths with trivial signature as the tree-like paths,
and prove that two paths are in tree-like equivalence if and only if
they have the same signature. In this way, we extend Chen's theorems
on the uniqueness of the sequence of iterated integrals associated
with a piecewise regular path to finite length paths and identify the
appropriate extended meaning for reparameterisation in the general
setting.
It is suggestive to think of this result as a non-commutative analogue of
the result that integrable functions on the circle are determined, up to
Lebesgue null sets, by their Fourier coefficients.
As a second theme we give quantitative versions of Chen's theorem in the
case of lattice paths and paths with continuous derivative, and as a
corollary derive results on the triviality of exponential products in the
tensor algebra.
\end{abstract}

\maketitle

\section{Introduction}

\subsection{Paths with finite length}

Paths, that is to say (right) continuous functions $\gamma $ mapping a
non-empty interval $J$ $\subset \mathbb{R}$ into a topological space $V$,
are fundamental objects in many areas of mathematics, and capture the
concept of an ordered evolution of events.

If $\left( V,d_{V}\right) $ is a metric space, then one of $\gamma $'s most
basic properties is its length $\left\vert \gamma \right\vert _{J}$. This
can be defined as 
\begin{equation*}
\left\vert \gamma \right\vert _{J}:=\sup_{\mathcal{D}\subset J}\sum 
_{\substack{ t_{i}\in \mathcal{D}  \\ i\neq 0}}d_{V}\left( \gamma
_{t_{i-1}},\gamma _{t_{i}}\right)
\end{equation*}%
where the supremum is taken over all finite partitions $\mathcal{D=}\left\{
t_{0}<t_{1}<\cdots <t_{r}\right\} $ of the interval $J$. It is clear that $%
\left\vert \gamma \right\vert $ is positive (although possibly infinite) and
independent of the parameterisation for $\gamma $. Letting $\tau \left(
t\right) =\left\vert \gamma \right\vert _{\left[ 0,t\right] }$ and setting $%
\eta \left( \tau \left( t\right) \right) =\gamma \left( t\right) $ one sees
that any continuous path of finite length can always be parameterised to
have unit speed. Paths of finite length are often said to be those of
bounded or finite variation.

\begin{definition}
We denote the set of paths of bounded variation by $BV,$ $BV$-paths with
values in $V$ by $BV(V),$ and those defined on $J$ by $BV\left( J,V\right) $.
\end{definition}

If $V$ is a vector space, then for any $\gamma \in BV\left( [0,t],V\right) $ and $\tau \in
BV\left( [0,s],V\right) $ we can form the concatenation $\gamma \ast \tau
\in BV\left( \left[ 0,s+t\right] ,V\right) $ 
\begin{eqnarray*}
\gamma \ast \tau \left( u\right) &=&\gamma \left( u\right) ,\ u\in \left[ 0,s%
\right] \\
\gamma \ast \tau \left( u\right) &=&\tau \left( u-s\right) +\gamma \left(
s\right) -\tau \left( 0\right) ,\ u\in \left[ s,s+t\right] .
\end{eqnarray*}%
The operation * is associative, and if $V$is a normed space, then $\left\vert
\gamma \right\vert +\left\vert \tau \right\vert =\left\vert \gamma \ast \tau
\right\vert$.


\subsection{Differential Equations}

One reason for looking at $BV(V)$ is that one can do calculus with these
paths, while at the same time the set of paths with $\left\vert \gamma
\right\vert _{J}\leq l$ is closed under the topology of pointwise
convergence (uniform convergence, \ldots ). Differential equations allow one
to express relationship between paths in $BV$. If $f^{i}$ are Lipschitz
vector fields on a space $W$ and $\gamma _{t}=\left( \gamma _{1}\left(
t\right) ,\ldots ,\gamma _{d}\left( t\right) \right) $ $\in BV\left( \mathbb{%
R}^{d}\right) $ then the differential equation 
\begin{eqnarray}
\frac{dy}{dt} &=&\sum_{i}f^{i}\frac{d\gamma _{i}}{dt}=f\left( y\right) \cdot 
\frac{d\gamma }{dt}  \label{eq:bvde} \\
y_{0} &=&a.  \notag
\end{eqnarray}%
has a unique solution for each $\gamma $. $BV$ is a natural class here for,
unless the vector fields commute, there is no meaningful way to make sense
of this equation if the path $\gamma $ is only assumed to be continuous.

If $\left( y,\gamma \right) $ solves the differential equation and $\left( 
\tilde{y},\tilde{\gamma}\right) $ are simultaneous reparameterisations, then
they also solve the equation and so it is customary to drop the $dt$ and
write 
\begin{equation*}
dy=\sum_{i}f^{i}\dot{\gamma}_{i}dt=\sum_{i}f^{i}d\gamma _{i}=f\left(
y\right) \cdot d\gamma .
\end{equation*}%
We can regard the location of $y_{s}$ as a variable and consider the
diffeomorphism $\pi _{st}$ defined by $\pi _{st}\left( y_{s}\right) :=y_{t}$%
. Then $\pi _{st}$ is a function of $\gamma |_{\left[ s,t\right] }$. One
observes that the map $\gamma |_{\left[ s,t\right] }\rightarrow \pi _{st}$
is a homomorphism from $\left( BV\left( V\right) ,\ast \right) $ to the
group of diffeomorphisms of the space $W$.

\subsection{Iterated integrals and the signature of a path}

One could ask which are the key features of $\gamma |_{\left[ s,t\right] }$,
which, with $y_{s}$, accurately predict the value $y_{t}$ in equation (\ref%
{eq:bvde}). The answer to this question can be found in a map from $BV$ into
the free tensor algebra!

\begin{definition}
Let $\gamma $ be a path of bounded variation on $\left[ S,T\right] $ with
values in a vector space $V$. Then its signature is the sequence of definite
iterated integrals 
\begin{align*}
\mathbf{X}_{S,T}& =\left( 1+X_{S,T}^{1}+\ldots +X_{S,T}^{k}+\ldots \right) \\
& =\left( 1+\int_{S<u<T}d\gamma _{u}+\ldots +\int_{S<u_{1}<\ldots
<u_{k}<T}d\gamma _{u_{1}}\otimes \ldots \otimes d\gamma _{u_{k}}+\ldots
\right)
\end{align*}
regarded as an element of an appropriate closure of the tensor algebra $%
T(V)=\bigoplus_{n=0}^{\infty }V^{\otimes n}$. 
\end{definition}

The signature is the \emph{definite} integral over the fixed interval where $%
\gamma $ is defined; re-parameterising $\gamma $ does not change its
signature. The first term $X_{\left[ S,t\right] }^{1}$ produces the path $%
\gamma $ (up to an additive constant). For convenience of notation, when we
have many paths, we will sometimes use a symbol such as $Y_{t}$ (instead of $%
\gamma _{t})$ for our path, $Y_{S,T}^{i}$ for the $i$-th coordinate of the
signature of $Y_{t}$, and $\mathbf{Y}_{S,T}$ for the signature of the path.
In some circumstances we will drop the time interval and just write $Y$ for
the path and $\mathbf{Y}$ for its signature. We call this map the signature
map and sometimes denote it by $S:X\rightarrow S\left( X\right) $ when this
helps our presentation.

The signature of $X$ is a natural object to study. The map $X\rightarrow 
\mathbf{X}$ is a homomorphism (c.f. Chen's identity \cite{lyons}) from the
monoid of paths with concatenation to (a group embedded in) the algebra $%
T\left( V\right) $. 
The signature $\mathbf{X}(=\mathbf{X}_{0,T})$ can be computed by solving the
differential equation 
\begin{eqnarray*}
d\mathbf{X}_{0,u} &=&\mathbf{X}_{0,u}\otimes dX_{u} \\
\mathbf{X}_{0,0} &=&\left( 1,0,0,\ldots \right) ,
\end{eqnarray*}%
and, in particular, paths with different signatures will have different
effects for some choice of differential equation.

There is a converse, although this is a consequence of our main theorem. If $%
X$ controls a system through a differential equation 
\begin{eqnarray*}
dY_{u} &=&f(Y_{u})dX_{u}, \\
Y_{0} &=&a,
\end{eqnarray*}%
and $f$ is Lipschitz, then the state $Y_{T}$ of the system after the
application of $X|_{\left[ 0,T\right] }$ is completely determined by the
signature $\mathbf{X}_{0,T}$ and $Y_{0}$. In other words the signature $%
\mathbf{X}_{S,T}$ is a truly fundamental representation for the bounded
variation path defined on $\left[ S,T\right] $ that captures its effect on
any non-linear system.

This paper explores the relationship between a path and its signature. We
determine a precise \emph{geometric} relation $\sim $ on bounded variation
paths, we prove that two paths of finite length are $\sim $-equivalent if
and only if they have the same signature: 
\begin{equation*}
X|_{J}\,\sim \,Y|_{K}\iff \mathbf{X}_{J}=\mathbf{Y}_{K}
\end{equation*}%
and hence prove that $\sim$ is an equivalence relation and identify
the sense in which the signature of a path determines the path.

The first detailed studies of the iterated integrals of paths are due to K.
T. Chen. In fact Chen \cite{KTCHEN01} proves the following theorems which
are clear precursors to our own results:

\textbf{Chen Theorem 1:} Let $d\gamma _{1},\cdots ,d\gamma _{d}$ be the
canonical 1-forms on $\mathbb{R}^{d}$. If $\alpha ,\beta \in \lbrack
a,b]\rightarrow R^{d}$ are irreducible piecewise regular continuous paths,
then the iterated integrals of the vector valued paths $\int_{\alpha \left(
0\right) }^{\alpha \left( t\right) }d\gamma $ and $\int_{\beta \left(
0\right) }^{\beta \left( t\right) }d\gamma $ agree if and only if there
exists a translation T of $\mathbb{R}^{d}$, and a continuous increasing
change of parameter $\lambda :[a,b]\rightarrow \lbrack a,b]$ such that $%
\alpha =T\beta \lambda $.

\textbf{Chen Theorem 2:} Let $G$ be a Lie group of dimension $d$, and let $%
\omega _{1}\cdots \omega _{d}$ be a basis for the left invariant $1$-forms
on G. If $\alpha ,\beta \in \lbrack a,b]\rightarrow G$ are irreducible
piecewise regular continuous paths, then the iterated integrals of the
vector valued paths $\int_{\alpha \left( 0\right) }^{\alpha \left( t\right)
}d\omega $ and $\int_{\beta \left( 0\right) }^{\beta \left( t\right)
}d\omega $ agree if and only if there exists a translation T of $G$, and a
continuous increasing change of parameter $\lambda :[a,b]\rightarrow \lbrack
a,b]$ such that $\alpha =T\beta \lambda $.\footnote{%
We borrow these formulations from the Math Review of the paper but include
the precise smoothness assumptions.}

In particular, Chen characterised piecewise regular paths in terms of their
signatures.

\subsection{The main results}

There are two essentially independent goals in this paper.

\begin{enumerate}
\item To provide quantitative versions of some of Chen's results. If $\gamma 
$ is continuous, of bounded variation and parameterised at unit speed, then
we will obtain lower bounds on the coefficients in the signature in terms of
the modulus of continuity of $\dot{\gamma}$ and the length of the path. For
example Theorem \ref{Thm:recoverlength} shows how one can recover the length
of a path $\gamma $ using the asymptotic magnitudes of these coefficients
(c.f. Tauberian theorems in Fourier Analysis). A detailed discussion is to
be found in Sections~2 and~3.

\item To prove a uniqueness theorem characterising paths of bounded
variation in terms of their signatures (c.f. the characterisation of
integrable functions in terms of their Fourier series) extending Chen's
theorem to the bounded variation setting.
\end{enumerate}

For this second goal we need a notion of tree-like path, our definition
codes $R$-trees by positive continuous functions on the line, as developed,
for instance, in \cite{snake}.

\begin{definition}
\label{def:tree-like} $X_{t},\;t\in \left[ 0,T\right] $ is a tree-like path
in $V$ if there exists a positive real valued continuous function $h$
defined on $\left[ 0,T\right] $ such that $h\left( 0\right) =h\left(
T\right) =0$ and such that 
\begin{equation*}
\left\Vert X_{t}-X_{s}\right\Vert _{V}\leq h\left( s\right) +h\left(
t\right) -2\inf_{u\in \left[ s,t\right] }h\left( u\right) .
\end{equation*}%
The function $h$ will be called a height function for $X$. We say $X$ is a
Lipschitz tree-like path if $h$ can be chosen to be of bounded variation. 
\end{definition}

\begin{definition}
Let $X,Y \in BV(V)$. We say $X\sim Y$ if the
concatenation of $X$ and $Y$ `run backwards' is a Lipschitz tree-like path.
\end{definition}

We now focus on $\mathbb{R}^d$ and state our main results.

\begin{theorem}
\label{thm:main} Let $X\in BV(\mathbb{R}^{d})$.
The path $X$ is tree-like if and only if the signature of $X$ is $\mathbf{0}%
=(1,0,0,\dots )$.
\end{theorem}

As the map $X\rightarrow \mathbf{X}$ is a homomorphism, and running a path
backwards gives the inverse for the signature in $T\left( V\right) $, an
immediate consequence of Theorem $\ref{thm:main}$ is

\begin{corollary}
If $X,Y \in BV({\mathbb R}^d)$, then $\mathbf{X}=\mathbf{Y}$ if and
only if the concatenation of $X$ and `$Y$ run backwards' is a Lipschitz
tree-like path.
\end{corollary}

\begin{corollary}
\label{cor:equiv} For $X,Y \in BV({\mathbb R}^d)$ the relation $X\sim
Y$ is an equivalence relation.
Concatenation respects $\sim $ and the equivalence classes $\Sigma $ form a
group under this operation.
\end{corollary}

There is an analogy between the space of paths of finite length in $\mathbb{R%
}^{d}$ and the space of words $a^{\pm 1}b^{\pm 1}\ldots c^{\pm 1}$ where the
letters $a,b,\ldots ,c$ are drawn from a $d$-letter alphabet $A$. Every such
word has a unique reduced form. This reduction respects the concatenation
operation and projects the space of words onto the free group. We extend
this result from paths on the integer lattice (words) to the bounded
variation case.

\begin{corollary}
\label{cor:minimisers} For any $X\in BV(\mathbb{R}^d)$ there exists a
unique path of minimal length, $\bar{X}$, called the reduced path, with the same
signature $\mathbf{X}=\mathbf{\bar{X}}$.
\end{corollary}

Taking these results together we see that the reduced paths form a group.
The multiplication operation is to concatenate the paths and then reduce the
result. One should note that this reduction process is not unique (although
we have proved that the reduced word one ultimately gets is). This group is
at the same time rather natural and concrete (a collection of paths of
finite length), but also very different to the usual finite dimensional Lie
groups. It admits more than one natural topology, and multiplication is not
continuous for the topology of bounded variation. 

We can restate these results in different language. The space $BV$
with $\ast$, the operation of concatenation, is a monoid.
Let $\mathcal{T}$ be
the set of tree-like paths in $BV$. Then $\mathcal{T}$ is also closed under
concatenation. If $\gamma \in BV$ and we use the notation $\gamma ^{-1}$ for $%
\gamma $ run backwards. It is clear from the definition that $\gamma
^{-1}\mathcal{T}\gamma \subset \mathcal{T}$ for all $\gamma \in BV.$ As we have proved that
tree-like equivalence is an equivalence relation $BV/\mathcal{T}$ is well defined,
closed under multiplication, and has inverses; it is a group.

We have the following picture%
\begin{equation*}
0\rightarrow \mathcal{T}\rightarrow BV\overset{\dashleftarrow }{\rightarrow }\Sigma
\rightarrow 0
\end{equation*}%
where one can regard $\Sigma $ as the $\sim$-equivalence classes of
paths or as the subgroup of the tensor algebra. The map $\dashleftarrow $
takes the class to the reduce path which is an element of $BV$. As
$\mathcal{T}$ has no natural $BV$-normal sub-monoids, one should expect that any continuous
homomorphism of $BV$ into a group will factor through $\Sigma $ if it is
trivial on the tree-like elements. 
It is clear that the set 
\begin{equation*}
\hat{\mathcal{T}}=\left\{ \left( \gamma ,h\right) ,\gamma \in \mathcal{T},\ h\text{ a height
function for }\gamma \right\}
\end{equation*}%
is contractable. An interesting question is whether $\mathcal{T}$ itself is contractable.

We prove in Lemma $\ref{lemma:plapprox}$ that any $\gamma \in \mathcal{T}$ is the
limit of weakly piecewise linear tree-like paths and hence
$\mathcal{T}$ is the \emph{%
smallest} multiplicatively closed and topologically closed set containing
the trivial path. This universality suggests that $\Sigma $ has similarities
to the Free group. One characterising property of the free group is that
every function from the alphabet $A$ into a group can be extended to a map
from words made from $A$ into paths in the group. The equivalent map for
bounded variation paths is Cartan development. Let $\theta $ be a linear map
of $\mathbb{R}^{d}\ $to the Lie algebra $\mathfrak{g}$ of a Lie group $G$
and let $X_{t}|_{t\leq T}$ be a bounded variation path. Then Cartan
development provides a canonical projection of $\theta \left( X\right) $ to
a path $Y$ starting at the origin in $G$ and we can define $\tilde{\theta}%
:X\rightarrow Y_{T}$. This map $\tilde{\theta}$ is a homomorphism from $%
\Sigma $ to $G$.

It is an exercise to prove that this map $\tilde{\theta}$ takes all
tree-like paths to the identity element in the group $G.$ As a consequence, $%
\tilde{\theta}$ is a map from paths of finite variation to $G$ which is
constant on each $\sim$ equivalence class and so defines a map from $\Sigma $
to $G$.

Let $X_{t}|_{t\leq T}$ be a path of bounded variation in $\mathbb{R}^{d}$
and suppose that for every linear map $\theta $ into a Lie algebra $%
\mathfrak{g}$, that $\tilde{\theta}\left( X\right) $ is trivial. As the
computation of the first $n$ terms in the signature is itself a development
(into the free $n$-step nilpotent group) we conclude that $\mathbf{X}_{0,T}=
\left( 1,0,0,\ldots \right) $ and so X is tree-like. In this way we have a

\begin{corollary}
A path of bounded variation is tree-like if and only if its development into
every finite dimensional Lie group is trivial.
\end{corollary}

The observation that any linear map of $\mathbb{R}^{d}\ $to the Lie algebra $%
\mathfrak{g}$ defines a map from $\Sigma $ to the Lie group is a universal
property of a kind giving further evidence that $\Sigma $ is some sort of continuous
analogue of the free group. However, $\Sigma $ is not a Lie group although
it has a Lie algebra and it is not characterised by this property. (Chen's
piecewise regular paths provide another example since they are paths of
bounded variation and are dense in the unit speed paths of finite length).

\subsection{Questions and Remarks}

How important to these results is the condition that the paths have finite
length? Does anything survive if one only insists that the paths are
continuous?

The space of continuous paths with the uniform topology is another
natural generalisation of words - certainly concatenation makes them a
monoid. However, despite their popularity in homotopy theory, there seems
little hope that a natural closed equivalence relation could be found on
this space that transforms it into a continuous `free group' in the sense we
mapped out above. The notion of tree-like makes good sense (one simply drops
the assumption that the height function $h$ is Lipschitz). With this
relaxation,

\begin{problem}
Does $\sim $ define an equivalence relation on continuous paths?
\end{problem}

Homotopy is the correct deformation of paths if one wants to preserve the
line integral of a path against a closed one-form. On the other hand
tree-like equivalence is the correct deformation of paths if one wants to
preserve the line integral of a path against any one form. As we mention
elsewhere in this paper, integration of continuous functions against general
one forms makes little sense. This is perhaps evidence to suggest the answer
to the problem is in the negative. The problem is in the transitivity of the
relation.

\begin{problem}
Is there a unique tree reduced path associated to any continuous path?
\end{problem}

For smooth paths $\gamma =\left( \gamma _{1},\gamma _{2}\right) $ in $%
\mathbb{R}^{2}$ Cartan development into the Heisenberg group is the map $%
\left( \gamma _{1},\gamma _{2}\right) \rightarrow \left( \gamma _{1},\gamma
_{2},\int \gamma _{1}d\gamma _{2}\right) .$ One knows \cite[Proposition 1.29]%
{LyStFlour} that there is no continuous bilinear map extending this
definition to any Banach space of paths which carries the Wiener measure. We
also know from Levy, that there are many \textquotedblleft almost sure"
constructions for this integral made in similar ways to \textquotedblleft
Levy area\textquotedblright . All are highly discontinuous and can give
different answers for the same Brownian path in $\mathbb{R}^{2}$. This wide
choice for the case of Brownian paths (which have finite $p$-variation for
every $p>2$) makes it clear there cannot be a canonical development for all
continuous paths.

The paper \cite{lyons} sets out a close relationship between differential
equations, the signature, and the notion of a geometric rough path. These
\textquotedblleft paths\textquotedblright\ also form a monoid under
concatenation and any linear map from $\mathbb{R}^{d}$ into the $\left(
p+\varepsilon \right) $-Lipschitz vector fields on a manifold $M$ induces a
canonical homomorphism of the $p$-rough paths with concatenation into the
group of diffeomorphisms of $M$ so they certainly have the analogy to the
Cartan development property. Similarly, every rough path has a signature,
and the map is a homomorphism.

\begin{problem}
Given a path $\gamma $ of finite $p$-variation for some $p>1$, is the
triviality of the signature of $\gamma $ equivalent to the path being
tree-like?
\end{problem}

Our theorem establishes this in the context of $p=1$ or
bounded variation paths but our proof uses the one dimensionality of the image of
the path in an essential way. An extension to $p$-rough paths with $p>1$
would 
require new ideas to account for the fact that these rougher paths are of
higher \textquotedblleft dimension\textquotedblright .

There seem to be many other natural questions.

By Corollary~\ref{cor:minimisers} among paths of finite length with the same
signature there is a unique shortest one - the reduced path. Successful
resolution of the following question could have wide ramifications in
numerical analysis and beyond. The question is interesting even for lattice
paths.

\begin{problem}
How does one effectively reconstruct the reduced path from its signature?
\end{problem}

A related question is to:

\begin{problem}
Identify those elements of the tensor algebra that are signatures of paths
and relate properties of the paths (for example their smoothness) to the
behaviour of the coefficients in the signature.
\end{problem}

Some interesting progress in this direction can be found in \cite{fawcett}.

We conclude with some wider comments.

\begin{enumerate}
\item There is an obvious link between these reduced paths and geometry
since each connection defines a closed subgroup of the group of reduced
paths (the paths whose developments are loops).

\item It also seems reasonable to ask about the extent to which the
intrinsic structure of the space of reduced paths (with finite length) in $%
d\geq 2$ dimensions changes as$\ d$ varies.
\end{enumerate}

\subsection{Outline}

We begin in Section 2 by discussing the lattice case. In this setting we can
obtain our first quantitative result on the signature. We do not have best
possible estimates, but we can prove that a word in the free group of length 
$L$ in two generators is completely reducible if the first $\left\lfloor
e\log \left( 1+\sqrt{2}\right) L\right\rfloor $ terms in the signature of
the path in the lattice corresponding to the word are zero. The case of
words in $d$ generators is also treated and if the first $c(d)L$ terms
in the signature are zero, the word is reducible, where the constant
$c(d)$ grows logarithmically in $d$.

In Section~3 we extend these quantitative estimates to finite length paths.
In order to do this we need to discuss the development of a path into a
suitable version of hyperbolic space - a technique that has more recently
proved useful in \cite{LySiDoob}. Using this idea we obtain a quantitative
estimate on the difference between the length of the developed path and its
chord in terms of the modulus of continuity of the derivative of the path.
This allows us to obtain, in the case where the derivative is
continuous, some estimates on the coefficients in the signature
and also shows how to recover the length of the path from the signature.

We can also prove for example that any path with bounded local curvature and
the first $N$ terms in the signature zero must be rather long or trivial - a
sort of rigidity theorem. We can obtain explicit bounds depending only on
the curvature bounds and $N$. However they are far from sharp as we
can see from the figure of 8, a path with curvature at most $%
4\pi$ and length one. It is clear that the first two terms in its
signature zero, but our results indicate that it cannot have
all of the first 115 terms in the signature zero! 

After this we return to the proof of our uniqueness result, the extension of
Chen's theorem. Our proof relies on various analytic tools (the Lebesgue
differentiation theorem, the area theorem), and particularly we introduce a
mollification of paths that retain certain deeply non-linear properties of
these paths to reduce the problem to the case where $\gamma $ is piecewise
linear. Piecewise linear paths are irreducible piecewise regular paths in
the sense of Chen and thus the result follows from Chen's Theorem. The
quantitative estimates we obtained give an independent proof for this
piecewise linear result.

In Section~\ref{SectiontreeLike} we establish the key properties for
tree-like paths that we need. In Section~5 we prove that any path $%
X_{t}|_{t\in \left[ 0,T\right] }\in BV$ and with trivial
signature can, after re-parameterisation, be uniformly approximated by
(weakly) piecewise linear paths with trivial signature. This is an
essentially non-linear result as the constraint of trivial signature
corresponds to an infinite sequence of polynomial constraints of increasing
complexity. In Section~6 we show that, by our quantitative version of 
Chen's theorem, such piecewise
linear paths must be reducible and so tree-like in our language.

This certainly gives us enough to show, in Section~7, that any weakly
piecewise linear path with trivial signature is tree-like. It is clear from
the definitions that uniform limits of tree-like paths with uniformly
bounded length are themselves tree-like. Applying the results of section 5
the argument is complete. We draw together all the parts to give the proofs
of our main Theorem and Corollaries in Section~8.

\section{Paths on the integer lattice}

\label{Section2}

\subsection{A discrete case of Chen's theorem}

Consider an alphabet $A$ and new letters $A^{-1}=\left\{ a^{-1},a\in
A\right\} $. Let $\Omega $ be the set of words in $A\cup A^{-1}$. Then $%
\Omega $ has a natural multiplication (concatenation) and an equivalence
relation that respects this multiplication.

\begin{definition}
\label{defn:reducedword}A word $w\in \Omega $ is said to cancel to the empty
word if, by applying successive applications of the rule 
\begin{equation*}
a\ldots bcc^{-1}d\ldots e\rightarrow a\ldots bd\ldots e,\quad \quad
a,b,c,d,e,\ldots \in A\cup A^{-1}
\end{equation*}
one can reduce $w$ to the empty word. We will say that $\left( a\ldots
b\right) $ is equivalent to $\left( e\ldots f\right) $ 
\begin{equation*}
\left( a\ldots b\right) \,\symbol{126}\,\left( e\ldots f\right)
\end{equation*}
if $\left( a\ldots bf^{-1}\ldots e^{-1}\right) $ cancels to the empty word.
\end{definition}

An easy induction argument shows that $\symbol{126}$ is an equivalence
relation. It is well known that the free group $F_{A}$ can be identified as $%
\Omega /\symbol{126}$. There is an obvious bijection between words in $%
\Omega $, and lattice paths, that is to say the piecewise linear paths $%
x_{u} $ which satisfy $x_{0}=0$ and $\left\Vert x_{k}-x_{k+1}\right\Vert =1$%
, are linear on each interval $u\in \left[ k,k+1\right] $, and have $%
x_{k}\in \mathbb{Z}^{\left\vert A\right\vert }$ for each $k$. The length of
the path is an integer equal to the number of letters in the word. The
equivalence relation between words can be re-articulated in the language of
lattice paths: Consider two lattice paths $x$ and $y$, and let $z$ be the
concatenation of $x$ with $y$ traversed backwards. Clearly, if $x$ and $y$
are equivalent then, keeping its endpoints fixed, $z$ can be
\textquotedblleft retracted\textquotedblright\ step by step to a point while
keeping the deformations inside what remains of the graph of $z$. The
converse is also true: if $U$ is the universal cover of the lattice, and we
identify based path segments in the lattice with points in $U$ then the
words equivalent to the empty word correspond with paths $x_{t}$ in the
lattice that lift to loops in $U$ . They are the paths that can be factored
into the composition of a loop in a tree with a projection of that tree into
the lattice. A loop in a tree is a tree-like path, as one can use the
distance from the basepoint of the loop as a height function. 


Chen's theorem tells us that any path that is not retractable to a point in
the sense of the previous paragraph has a non-trivial signature. Our
quantitative approach allows us to prove an algebraic version of this
result. Let $\gamma _{w}$ be the lattice path associated to the word $%
w=a_{1}^{\sigma _{1}}\ldots a_{L}^{\sigma _{L}}$ (where $\sigma
=(\sigma_1,\dots,\sigma_L) \in \left\{
\pm 1\right\} ^{L}$ gives the signs associated to each letter). As the
signature is a homomorphism, we have $S\left( \gamma _{w}\right) =S\left(
\gamma _{a_{1}^{\sigma _{1}}}\right) \ldots S\left( \gamma _{a_{L}^{\sigma
_{L}}}\right) $. Since $\gamma _{a_{i}^{\sigma _{i}}}$ is a path that moves $%
\sigma _{i}$ units in a straight line in the $a_{i}$ direction, its
signature is the exponential and $S\left( \gamma _{w}\right)
=e^{\sigma _{1}a_{1}}\otimes \ldots \otimes e^{\sigma _{L}a_{L}}.$

Our quantitative approach will show in Theorem~\ref{prop:intlat} 
that for a word of length $L$ in a two letter alphabet, if 
\begin{equation*}
e^{\sigma _{1}a_{1}}\otimes \ldots \otimes e^{\sigma _{L}a_{L}}=\left(
1,0,0,\ldots ,0,X^{N\left( L\right) +1},X^{N\left( L\right) +2}\ldots
\right) ,\ \ \sigma \in \left\{ \pm 1\right\} ^{L}
\end{equation*}%
where $N\left( L\right) =\left\lfloor e\log \left( 1+\sqrt{2}\right)
L\right\rfloor $, then there is an $i$ for which $a_{i}=-a_{i+1}$ and by
induction the reduced word is trivial.

The proof is based on regarding $\mathbb{R}^{d}$ as the tangent space to a
point in $d$-dimensional hyperbolic space $\mathbb{H}$, scaling the path $%
\gamma $ and developing it into hyperbolic space. There are two ways to view
this development of the path, one of which yields analytic information out
of the iterated integrals, the other geometric information. Together they
quickly give the result. We work in two dimensional hyperbolic space and, at
the end, show that the general case can be reduced to this one.

\subsection{The universal cover as a subset of $\mathbb{H}$}

Let $X$ be a lattice path in $\mathbb{R}^{2}$, $\theta \geq 0$, and $%
X^{\theta }=\theta X$ be the re-scaled lattice path. The development $%
Y^{\theta }$ of $X^{\theta }$ into $\mathbb{H}$ moves along successive
geodesic segments of length $\theta $ in $\mathbb{H}$, each time $X^{\theta
} $ turns a corner, so does $Y^{\theta }$ and angles are preserved.

For a fixed choice of $\theta $ we can trace out in $\mathbb{H}$ the four
geodesic segments from the origin, the three segments out from each of
these, and the three from each of these, and so on. It is clear that if the
scale $\theta $ is large enough, the negative curvature forces the image to
be tree. This will happen exactly when the path that starts by going along
the real axis and then always turns anti-clockwise never hits its reflection
in the line $x=y$.

The successive moves can be expressed as iterations of a Mobius transform, 
\begin{align*}
m\left( x\right) & :=\frac{-ir+x}{-i-rx} \\
x_{n}& =m^{n}\left( 0\right),
\end{align*}%
and if $r=1/\sqrt{2}$, then the trajectory eventually ends at $\left(
1+i\right) /\sqrt{2}$. Hyperbolic convexity ensures that all these
trajectories are (after the first linear step) always in the region
contained by the geodesic from $\left( 1+i\right) /\sqrt{2}$ to $\left(
1-i\right) /\sqrt{2}.$ In particular they never intersect the trajectories
whose first move is from zero to $i,$ to $-i$, or to $-1$. Now, there is
nothing special about zero in this discussion, and using conformal
invariance it is easy to see that

\begin{lemma}
If $\theta$ is at least equal to the hyperbolic distance from $0$ to $1/%
\sqrt{2}$ in $\mathbb{H}$, then the path $Y^{\theta}_{t}$ takes its values in
a tree. This value $1/\sqrt{2}$ is sharp.
\end{lemma}

We have developed $X^{\theta }$ into a tree in $\mathbb{H}$; we
have already observed that a loop in a tree is tree-like. If we can prove
that $Y_{T}^{\theta }=Y_{0}^{\theta }$, the $Y^{\theta }$ will be tree-like
and hence so will $X^{\theta }$and $X$. To achieve this we must use the
assumption that the path has finite length and that all its iterated
integrals are zero from a different perspective.

\subsection{Cartan development as a linear differential equation}

If $G$ is a closed subgroup of the matrices, and $X_{t}|_{t\leq T}$ is a
path in its Lie algebra $\mathfrak{g}$, then the equation for the Cartan
development $M_{T}\in G$ of $X_{t}|_{t\leq T}\in \mathfrak{g}$ is given by
the differential equation%
\begin{equation*}
M_{t+\delta t}\approx M_{t}\exp \left( \delta X_{t}\right) 
\mbox{
  or equivalently }dM_{t}=M_{t}dX_{t}.
\end{equation*}%
The development of a smooth path in the tangent space to $0$ in $\mathbb{H}$ 
\begin{equation*}
\mathbb{H=}\left\{ z\in \mathbb{C}|\left\Vert z\right\Vert <1\right\} ,
\end{equation*}%
is also expressible as a differential equation. However, it is easier to
express this development in terms of Cartan development in the group of
isometries regarded as matrices in $GL(2,\mathbb{C})$ rather than on the points of $%
\mathbb{H}$. We identify $\mathbb{R}^{2}$ with the Lie subspace 
\begin{equation*}
\left( 
\begin{array}{cc}
0 & x+iy \\ 
x-iy & 0%
\end{array}%
\right) .
\end{equation*}%
In this representation, the equation for $M_{t}$ is linear and so we have an
expansion for $M$: 
\begin{align*}
M_{T}& =M_{0}\left(
I+\int_{0<u<T}dX_{u}+\int_{0<u_{1}<u_{2}<T}dX_{u_{1}}dX_{u_{2}}+\ldots
\right)  \\
& =M_{0}\times \left( 
\begin{array}{cc}
a & b \\ 
\bar{b} & \bar{a}%
\end{array}%
\right) ,
\end{align*}%
where 
\begin{eqnarray*}
a &=&1+\sum_{k}\int_{0<u_{1}<u_{2}<\ldots <u_{2k}<T}dX_{u_{1}}d\bar{X}%
_{u_{2}}\ldots dX_{u_{2k-1}}d\bar{X}_{u_{2k}} \\
b &=&\sum_{k}\int_{0<u_{1}<u_{2}<\ldots <u_{2k-1}<T}dX_{u_{1}}d\bar{X}%
_{u_{2}}\ldots dX_{u_{2k-1}}
\end{eqnarray*}%
and $\int_{0<u_{1}<\dots <u_{2k}<T}dX_{u_{1}}d\bar{X}_{u_{2}}\ldots
dX_{u_{2k-1}}d\bar{X}_{u_{2k}}$ is now, with an abuse of notation, a complex
number. We have an a priori bound:

\begin{lemma}
\label{lem:intpathapriori} If $X$ is a path of length exactly $\theta L$,
then 
\begin{equation*}
\left| \int_{0<u_{1}<u_{2}<\ldots<u_{2k}<T}dX_{u_{1}}d\bar{X}_{u_{2}}\ldots
dX_{u_{2k-1}}d\bar{X}_{u_{2k}}\right| <\frac{\left( \theta L\right) ^{2k}}{%
\left( 2k\right) !}.
\end{equation*}
\end{lemma}

To use this lemma we need to be able to estimate the tail of an exponential
series. The following lemma (based on Stirling's formula) articulates a
convenient inequality.

\begin{lemma}
Let $x\geq 1/e$. \label{lem:stirling} 
(1) $\frac{x^{m}}{m!}<\frac{\xi_{0}}{m^{1/2}}$ holds for all $m\geq ex$. 
\newline
(2) If for any $k$ one has $m\geq ex+k$, then 
\begin{align*}
\sum_{r\geq m}\frac{x^{r}}{r!} & \leq\frac{e^{\frac{1}{2}}}{\sqrt{2\pi}(e-1)}
e^{-k}x^{-1/2} \simeq 0.38 \;e^{-k}x^{-1/2}.
\end{align*}
\end{lemma}

\begin{proof}
By Stirling's formula $\lim_{y\rightarrow \infty }\frac{e^{-y}y^{\frac{1}{2}%
+y}}{y!} = \frac{1}{\sqrt{2\pi}}$ and is approached monotonely from below.
It is an upper bound and also a good global approximation to $\frac{e^{-y}y^{%
\frac{1}{2}+y}}{y!}$ valid for all $y\geq 1$. Putting $y=ex$ gives 
\begin{align*}
\frac{e^{-ex}\left( ex\right) ^{\frac{1}{2}+ex}}{\left( ex\right) !}& < 
\frac{1}{\sqrt{2\pi}} \\
\frac{x^{ex}}{\left( ex\right) !}& <e^{-\frac{1}{2}}x^{-\frac{1}{2}} \frac{1%
}{\sqrt{2\pi}}.
\end{align*}
Moreover the recurrence relation for the $!$ function implies, for every $%
k\in \mathbb{Z}$ with $ex+k>0$, that 
\begin{equation*}
\frac{x^{ex}}{\left( ex\right) !}\geq e^{k}\frac{x^{ex+k}}{\left(
ex+k\right) !}
\end{equation*}
and so 
\begin{equation*}
\frac{x^{ex+k}}{\left( ex+k\right) !}<e^{-k-\frac{1}{2}}x^{-\frac{1}{2}} 
\frac{1}{\sqrt{2\pi}},
\end{equation*}
establishing the first claim. Now summing this bound we have 
\begin{align*}
\sum_{k\geq 0}\frac{x^{ex+k}}{\left( ex+k\right) !}& \leq e^{-\frac{1}{2}%
}x^{-\frac{1}{2}}\frac{1}{\sqrt{2\pi}}\sum_{k\geq 0}e^{-k} \\
& =\frac{e^{\frac{1}{2}}}{e-1}\frac{1}{\sqrt{2\pi}} x^{-\frac{1}{2}}.
\end{align*}
Since for $ex>0$ the function $k\rightarrow \frac{x^{ex+k}}{\left(
ex+k\right) !}$ is monotone decreasing on $\mathbb{R}^{+}$ we see that 
\begin{equation*}
\sum_{m\geq ex}\frac{x^{m}}{m!}\leq \frac{e^{\frac{1}{2}}}{e-1}\frac{1}{%
\sqrt{2\pi}} x^{-\frac{1}{2}},
\end{equation*}
completing the proof of the lemma. Finally we note the approximate value of
the constant: 
\begin{equation*}
\frac{e^{\frac{1}{2}}}{\sqrt{2\pi}(e-1)} \simeq 0.38.
\end{equation*}
\end{proof}

\subsection{The signature of a word of length $L$}

We deduce the following totally algebraic corollary for paths $X$ that have
traversed at most $L$ vertices.

\begin{theorem}
\label{prop:intlat} If a path of length $L$ in the two dimensional integer
lattice (corresponding to a word with $L$ letters drawn from a two letter
alphabet and its inverse), has the first $\lfloor e\log (1+\sqrt{2})L\rfloor 
$ $GL(2,\mathbb{C})$-iterated integrals \footnote{$GL(2,\mathbb{C})$%
-iterated integrals: since our path is in a vector subspace of the algebra $%
GL(2,\mathbb{C})$ we may compute the iterated integrals in the algebra $GL(2,%
\mathbb{C})$ or in the tensor algebra over the vector subspace. There is a
natural algebra homomorphism of the tensor algebra onto $GL(2,\mathbb{C})$.
The $GL(2,\mathbb{C})$-iterated integrals are the images of those in the
tensor algebra under this projection and \`{a} priori contain less
information.} zero, then all iterated integrals (in the tensor algebra) are
zero, the path is tree-like, and the corresponding reduced word is trivial.
\end{theorem}

\begin{proof}
Any Mobius transformation preserving the disk can be expressed as 
\begin{equation*}
M=\left( 
\begin{array}{cc}
z & 0 \\ 
0 & \bar{z}%
\end{array}%
\right) \left( 
\begin{array}{cc}
\frac{1}{\sqrt{1-r^{2}}} & \frac{r}{\sqrt{1-r^{2}}} \\ 
\frac{r}{\sqrt{1-r^{2}}} & \frac{1}{\sqrt{1-r^{2}}}%
\end{array}%
\right) \left( 
\begin{array}{cc}
\omega  & 0 \\ 
0 & \bar{\omega}%
\end{array}%
\right) ,
\end{equation*}%
where $\left\vert z\right\vert =\left\vert \omega \right\vert =1$ and $r$ is
the Euclidean distance from $0$ to $M0.$ Now 
\begin{eqnarray*}
Tr\left[ A\overline{B}\right]  &=&\sum_{i}\sum_{j}a_{ij}\overline{b}_{ji} \\
&=&\overline{Tr\left[ B\overline{A}^{t}\right] }
\end{eqnarray*}%
and 
\begin{equation*}
\left( 
\begin{array}{cc}
\omega  & 0 \\ 
0 & \bar{\omega}%
\end{array}%
\right) \overline{\left( 
\begin{array}{cc}
\omega  & 0 \\ 
0 & \bar{\omega}%
\end{array}%
\right) ^{T}}=\left( 
\begin{array}{cc}
1 & 0 \\ 
0 & 1%
\end{array}%
\right) 
\end{equation*}%
hence 
\begin{align*}
Tr\left[ M\overline{M}^{t}\right] & =Tr\left[ \left( 
\begin{array}{cc}
\frac{1}{\sqrt{1-r^{2}}} & \frac{r}{\sqrt{1-r^{2}}} \\ 
\frac{r}{\sqrt{1-r^{2}}} & \frac{1}{\sqrt{1-r^{2}}}%
\end{array}%
\right) ^{2}\right]  \\
& =\frac{2\left( 1+r^{2}\right) }{\left( 1-r^{2}\right) }.
\end{align*}%
Letting $r=1/\sqrt{2}$ we see that if 
\begin{equation*}
Tr\left[ M\overline{M}^{t}\right] <6,
\end{equation*}%
then the image of $0$ under the Mobius transformation must lie in the circle
of radius $1/\sqrt{2}$. On the other hand 
\begin{equation*}
Tr\left[ \left( 
\begin{array}{cc}
a & b \\ 
\bar{b} & \bar{a}%
\end{array}%
\right) \left( 
\begin{array}{cc}
\bar{a} & b \\ 
\bar{b} & a%
\end{array}%
\right) \right] =2\left( \left\vert a\right\vert ^{2}+\left\vert
b\right\vert ^{2}\right) 
\end{equation*}%
and in our context, where the first $N$ iterated integrals are
zero, this gives the inequality 
\begin{eqnarray}
&&\left\vert 1+\sum_{k>N}\int_{0<u_{1}<\dots <u_{2k}<T}dX_{u_{1}}d%
\bar{X}_{u_{2}}\ldots dX_{u_{2k-1}}d\bar{X}_{u_{2k}}\right\vert ^{2}+  \notag
\\
&&\qquad \qquad \left\vert \sum_{k>N}\int_{0<u_{1}<\dots
<u_{2k-1}<T}dX_{u_{1}}d\bar{X}_{u_{2}}\ldots dX_{u_{2k-1}}\right\vert ^{2}<6.
\label{eq:rineq}
\end{eqnarray}%
Using our a priori estimate from Lemma~\ref{lem:intpathapriori} we have that
the inequality will hold if 
\begin{equation*}
\left( 1+\sum_{k>N}\frac{(\theta L)^{2k-1}}{(2k-1)!}\right)
^{2}+\left( \sum_{k>N}\frac{(\theta L)^{2k}}{(2k)!}\right) ^{2}<6.
\end{equation*}%
Observing that, as we will choose $N>\theta L$, the terms in the
sums are decreasing, we have 
\begin{equation*}
\sum_{k>N}\frac{(\theta L)^{2k-1}}{(2k-1)!}<s
\end{equation*}%
and then 
\begin{equation*}
\sum_{k>N}\frac{(\theta L)^{2k}}{(2k)!}<s.
\end{equation*}%
We see that (\ref{eq:rineq}) will always be satisfied if we choose $s$ such
that 
\begin{equation*}
2\left( s+s^{2}\right) <5\mbox{ and }\sum_{k>N}\frac{(\theta
L)^{2k-1}}{(2k-1)!}<s
\end{equation*}%
Hence, if 
\begin{equation*}
\sum_{k>N}\frac{(\theta L)^{2k-1}}{(2k-1)!}<\frac{\sqrt{11}-1}{2},
\end{equation*}%
then $Y_{T}=0$. By Lemma~\ref{lem:stirling} (2) with $x=\log (1+\sqrt{2})L$,
we have if $N\geq e\log (1+\sqrt{2})L$, 
\begin{eqnarray*}
\sum_{m\geq N}\frac{\left( \log (1+\sqrt{2})L\right) ^{m}}{m!}
&\leq &0.38..(\log (1+\sqrt{2})L)^{-1/2} \\
&<&\frac{\sqrt{11}-1}{2}
\end{eqnarray*}%
for all $L\geq 1$.

Observe that if $\theta \geq \log \left[ 1+\sqrt{2}\right] $ then $%
Y_{t}=M_{t}0$ lies in a tree, and the development of $Y$ is such that every
vertex of the tree is at least a distance $\theta $ from the origin except
the origin itself. By our hypotheses and the above argument $d\left[ Y_{T},0%
\right] <\theta $ and hence $Y_{T}=0$. Therefore $Y$ is tree-like and the
reduced word is trivial.
\end{proof}




Finally we note that the case of the free group with two generators is
enough to obtain a general result as the free group on $d$ generators can be
embedded in it.

\begin{lemma}
Suppose that $\Gamma _{d}$ is the free group on $d$ letters $e_{i}$ and that 
$\Gamma $ is the free group on the letters $a,\;b$. Then we can identify $%
f_{i}\in \Gamma $ so that the homomorphism induced by $e_{i}\rightarrow
f_{i} $ from $\Gamma _{d}$ to $\Gamma $ is an isomorphism and so that the
length of the reduced words $f_{i}$ are at most $\left| f_{i}\right| \leq
2\left\lceil \log _{3}\frac{d}{2}\right\rceil +3.$
\end{lemma}

\begin{proof}
It is enough to show that we can embed $\Gamma _{23^{l-1}}$ into $\Gamma $
so that each $f_{i}$ has length $l$. Consider the collection of all reduced
words of length $l$ in $\Gamma $. There are $43^{l-1}$ of them if $l>0$.
Partition them into pairs, so that the left most letter of each of the words
in a pair is the same up to inverses. Order them lexicographically. Now
consider the space which is the ball in the Cayley graph of $\Gamma $
comprising reduced words with length at most $l$. It is obvious that this is
a contractable space. Now adjoin new edges connecting the ends of our pairs.
Associate with each of the new edges the alternate letter and orient the
edge to point from the lower to the higher word in the lexicographic order.
Then this new space $\Delta $ is contractable to $23^{l-1}$ loops and so has
the free group $\Gamma _{23^{l-1}}$ as its fundamental group. On the other
hand, we can obviously lift any path in $\Delta $ to the Cayley graph of $%
\Gamma$; the map from loops in $\Delta $ to $\Gamma $ is a homomorphism. The
monodromy theorem tells us that this homomorphism induces a homomorphism of
the homotopy group of $\Delta $ to $\Gamma $. As $\Gamma $ is a tree,
any two lifts of paths with the same endpoint in $\Gamma $ are homotopic
relative to those endpoints in graph $\Gamma $. Therefore we can associate
every point in the homotopy group of $\Delta $ with a unique element of $%
\Gamma $ and see that the homomorphism is injective. So we see that the
image is a copy of the free group $\Gamma _{23^{l-1}}$. The generators of
the classes in $\Delta $ clearly lift to paths of length $2l+1$ in $\Gamma $
and we take the end points of these paths to be the $f_{i}$.
\end{proof}

\begin{theorem}
If $X$ is a path of length $L$ in the $d$-dimensional integer lattice and
the projections into $GL(2,\mathbb{C})$ of the first $\lfloor (2\left\lceil
\log _{3}\frac{d}{2}\right\rceil +3)e\log (1+\sqrt{2})L\rfloor $ iterated
integrals are zero, then the path is tree-like.
\end{theorem}

In this section, our arguments depend on the tree-like nature of the
development of the path in the lattice and little else - this is a property
of the development into any rank one symmetric space but is still plausible,
if less obvious for general homogeneous spaces. Each space will give rise to
a different class of iterated integrals that are sufficient to determine the
tree-like nature of a path in a `jungle gym'. One should note that computing
the iterated integrals is not the most efficient way to determine if a word
is reducible if the word, as opposed to its signature, is presented.

\section{Quantitative versions of Chen's Theorem}

\label{Section3}We work in the hyperboloid model for $\mathbb{H}$\ (which
embeds the space $\mathbb{H}$ into a $d+1$-dimensional Lorentz space)
because the isometries of $\mathbb{H}$ extend to linear maps. 

Consider the quadratic form on $\mathbb{R}^{d+1}$defined by 
\begin{equation*}
I_{d}\left( x,y\right) =\sum_{1}^{d}x_{i}y_{i}-x_{d+1}y_{d+1}
\end{equation*}%
and the surface 
\begin{equation*}
\mathbb{H=}\left\{ x,\;I_{d}\left( x,x\right) =-1\right\} .
\end{equation*}%
Then $\mathbb{H}$ is hyperbolic space with the metric obtained by
restricting $I_{d}$ to the tangent spaces to $\mathbb{H}$. (If $x\in \mathbb{%
H}$ then $\left\{ y|I_{d}\left( y,x\right) =0\right\} $ is the tangent space
to $\mathbb{H}$ in $\mathbb{R}^{d+1}$ and moreover $I_{d}\left( z,z\right) $
is positive definite on $z\in \left\{ y|I_{d}\left( y,x\right) =0\right\} $
and so this inner product is a Riemannian structure on $\mathbb{H)}$. In
fact, (see \cite{cfkp}, p83) distances in $\mathbb{H}$ can be calculated
using $I_{d}$ 
\begin{equation}
-\cosh d\left( x,y\right) =I_{d}\left( x,y\right)  \label{eqn:metric}
\end{equation}%
If $SO\left( I_{d}\right) $ denotes the group of matrices with positive
determinant preserving\footnote{%
Precisely, $M\in SO\left( I_{n}\right) $ if $I_{d}\left( \left(
My^{t}\right) ^{t},\left( Mx^{t}\right) ^{t}\right) \equiv I_{d}\left(
y,x\right) $} the quadratic form $I_{d}$ then one can prove this is exactly
the group of orientation preserving isometries of $\mathbb{H}$. The Lie
algebra of $SO\left( I_{d}\right) $ is easily recognised as the $d+1$
dimensional matrices that are antisymmetric in the top left $d\times d$
block and symmetric in the last column and bottom row and zero in the bottom
right corner. Then the development of a path $\gamma \in \mathbb{R}^{d}$ to $%
SO\left( I_{d}\right) $ and $\mathbb{H}$ (chosen to commute with the action
of multiplication on the right in $SO\left( I_{d}\right) )$ is given by
solving the following differential equation 
\begin{equation}
d\Gamma _{t}=\left( 
\begin{array}{cccc}
0 & \cdots & 0 & d\gamma _{t}^{1} \\ 
\vdots & \ddots & \vdots & \vdots \\ 
0 & \cdots & 0 & d\gamma _{t}^{d} \\ 
d\gamma _{t}^{1} & \cdots & d\gamma _{t}^{d} & 0%
\end{array}%
\right) \Gamma _{t}.  \label{eq:development}
\end{equation}%
We define $X$ to be the development of the path $\gamma $ to the path in $%
\mathbb{H}$ starting at $o=\left( 0,\cdots ,0,1\right) ^{t}\ $and given by 
\begin{equation*}
X_{t}=\Gamma _{t}o.
\end{equation*}%
Now we can write $d\Gamma _{t}=F(d\gamma _{t})\Gamma _{t}$ where 
\begin{equation*}
F:x\rightarrow \left( 
\begin{array}{cccc}
0 & \cdots & 0 & x_{1} \\ 
\vdots & \ddots & \vdots & \vdots \\ 
0 & \cdots & 0 & x_{d} \\ 
x_{1} & \cdots & x_{d} & 0%
\end{array}%
\right)
\end{equation*}%
is a map from $\mathbb{R}^{d}$ to $\mathrm{Hom}\left( \mathbb{R}^{d+1},%
\mathbb{R}^{d+1}\right) $, where for precision we choose the Euclidean norm
on $\mathbb{R}^{d}$ and $\mathbb{R}^{d+1}$ and the operator norm on $\mathrm{%
Hom}\left( \mathbb{R}^{d+1},\mathbb{R}^{d+1}\right) $.

\begin{lemma}
\label{lem:Fnorm} In fact $\left| \left| F\right| \right|_{Hom\left( \mathbb{%
R}^{d+1}, \mathbb{R}^{d+1}\right)}=1. $
\end{lemma}

\begin{proof}
Let $e\in \mathbb{R}^{d}$ and $f\in \mathbb{R}$. Then for $x\in\mathbb{R}^d$ 
\begin{equation*}
F\left( x\right) \left( 
\begin{array}{c}
e \\ 
f%
\end{array}
\right) =\left( 
\begin{array}{c}
fx \\ 
e.x%
\end{array}
\right)
\end{equation*}
and computing norms 
\begin{equation*}
\left| \left| \left( 
\begin{array}{c}
fx \\ 
e.x%
\end{array}
\right) \right| \right| ^{2} \leq f^{2}\left| \left| x\right| \right|
^{2}+\left| \left| e\right| \right| ^{2}\left| \left| x\right| \right| ^{2}
= \left| \left| \left( 
\begin{array}{c}
e \\ 
f%
\end{array}
\right) \right| \right| ^{2}\left| \left| x\right| \right| ^{2}
\end{equation*}
and hence $\|F\|_{Hom\left(\mathbb{R}^{d+1}, \mathbb{R}^{d+1}\right)}=1$. 
\end{proof}

\subsection{Paths close to a geodesic}

We are interested in developing paths $\gamma $ of fixed length $l$ into
paths $\Gamma $ in $SO\left( I_{d}\right) $ and in the function 
\begin{equation*}
\varrho \left( \gamma \right) :=d\left( o,\Gamma o\right)
\end{equation*}%
giving the length of the chord connecting the beginning and end of the
development of $\gamma $ into Hyperbolic space. Amongst these paths $\gamma $
of fixed length, straight lines maximise $\varrho $ as the developments are
geodesics. The function $\varrho $ is a smooth function on path space \cite{LixLy2006}.
Therefore one would expect that for some constant $K$ 
\begin{equation*}
\varrho \left( \gamma \right) \geq l-K\varepsilon ^{2}
\end{equation*}%
whenever $\gamma $ is in the $\varepsilon $-neighbourhood (for the
appropriate norm) of a straight line. We will make this precise using
Taylor's theorem.

Suppose our straight line is in the direction of a unit vector $v$. If our
path $\gamma $ is parameterised at unit speed we can represent it by 
\begin{equation*}
d\gamma _{t}=\Theta _{t}vdt,
\end{equation*}%
where $\Theta _{t}$ is a path in the isometries of $\mathbb{R}^{d}$. In this
discussion we assume that $\Theta _{t}$ is continuous and has modulus of
continuity $\delta $. Of course, $\gamma $ is close to $t\rightarrow tv$ if $%
\Theta $ is uniformly close to the identity. Consider the development $%
\Gamma _{t}o=\left( \hat{x}_{t},x_{t}\right) $ of $\gamma $ into the
hyperboloid model of $\mathbb{H}$ defined by 
\begin{eqnarray*}
x_{t} &\in &\mathbb{R} \\
\hat{x}_{t} &\in &\mathbb{R}^{d} \\
dx_{t} &=&\hat{x}_{t}.\Theta _{t}vdt\ \ \ \ \ \ x_{0}=1 \\
d\hat{x}_{t} &=&x_{t}\Theta _{t}vdt\ \ \ \ \ \ \ \text{~}\hat{x}_{0}=0
\end{eqnarray*}%
We know that $\left\Vert \hat{x}_{t}\right\Vert ^{2}+1=\left\vert
x_{t}\right\vert ^{2}$ and that 
\begin{eqnarray*}
\cosh d\left( \Gamma _{t}o,o\right) &=&-I_{d}\left( \left( 
\begin{array}{c}
\hat{x}_{t} \\ 
x_{t}%
\end{array}%
\right) ,\left( 
\begin{array}{c}
0 \\ 
1%
\end{array}%
\right) \right) \\
&=&x_{t}
\end{eqnarray*}%
in other words $\cosh \varrho \left( \gamma |_{\left[ 0,t\right] }\right)
=x_{t}$.

\begin{proposition}
\label{prop:geod} Suppose the one can express $\Theta _{t}$ in the form e$%
^{A_{t}}$ where $A_{t}$ is a continuously varying anti-symmetric matrix and
that $\left\Vert A\right\Vert _{\infty }\leq \eta < 1$. Then 
\begin{equation*}
\left\vert \cosh T-x_{T}\right\vert \leq 4^{T}\frac{\left\Vert A\right\Vert
_{\infty }^{2}}{2}
\end{equation*}
\end{proposition}

\begin{proof}
Suppose $\varepsilon \in \left[ -1,1\right] $. We can introduce a family of
paths $\gamma _{t}^{\varepsilon }$ with $\gamma _{t}^{1}\equiv $ $\gamma
_{t} $ and with $\gamma _{t}^{0}$ the straight line $tv$ by setting 
\begin{eqnarray*}
d\gamma _{t}^{\varepsilon } &=&e^{\varepsilon A_{t}}vdt \\
\gamma _{0}^{\varepsilon } &=&0.
\end{eqnarray*}%
We can then consider the real valued function $f$ on $\left[ -1,1\right] $
comparing the length of the development of $\gamma ^{\varepsilon }$ and the
straight line 
\begin{equation*}
f\left( \varepsilon \right) :=\cosh \varrho \left( \gamma _{t}^{\varepsilon
}|_{t\in \left[ 0,T\right] }\right) -\cosh T.
\end{equation*}%
Of course $f\left( 0\right) =0$ and $f\leq 0$. Now \cite[Theorem 2.2]%
{LixLy2006} proves that development of a path $\gamma $ is Frechet
differentiable as a map from paths to paths in all $p$- variation norms with 
$p\in \left[ 1,2\right) .$

It is elementary that%
\begin{eqnarray*}
d\left( \gamma _{t}^{\varepsilon }-\gamma _{t}^{\varepsilon +h}\right)
&=&e^{\varepsilon A_{t}}\left( 1-e^{hA_{t}}\right) vdt \\
&=&hA_{t}e^{\varepsilon A_{t}}vdt+\frac{1}{2}\tilde{h}_{t}^{2}A_{t}^{2}e^{%
\varepsilon A_{t}}vdt
\end{eqnarray*}%
where $\tilde{h}_{t}\in \left[ 0,h\right] $. Working towards the $1$%
-variation derivative 
\begin{eqnarray*}
\int_{t\in \left[ 0,T\right] }\left\vert d\left( \gamma _{t}^{\varepsilon
}-\gamma _{t}^{\varepsilon +h}\right) -hA_{t}e^{\varepsilon
A_{t}}vdt\right\vert &\leq &\int_{t\in \left[ 0,T\right] }\left\vert \frac{1%
}{2}\tilde{h}_{t}^{2}A_{t}^{2}e^{\varepsilon A_{t}}vdt\right\vert \\
&\leq &\frac{h^{2}}{2}\int_{t\in \left[ 0,T\right] }\left\vert
A_{t}^{2}\right\vert dt
\end{eqnarray*}%
and $\varepsilon \rightarrow \gamma ^{\varepsilon }$ is differentiable with
derivative 
\begin{equation*}
d\gamma ^{\left( 1\right) ,\varepsilon }:=A_{t}e^{\varepsilon A_{t}}vdt,
\end{equation*}%
providing $\int_{t\in \left[ 0,T\right] }\left\vert A_{t}^{2}\right\vert
dt<\infty $. A similar estimate shows that the derivative of $\gamma
^{\left( 1\right) ,\varepsilon }$ exists and is 
\begin{equation*}
d\gamma ^{\left( 2\right) ,\varepsilon }:=A_{t}^{2}e^{\varepsilon A_{t}}vdt,
\end{equation*}%
providing $\int_{t\in \left[ 0,T\right] }\left\vert A_{t}^{3}\right\vert
dt<\infty $. From \cite[Theorem 2.2]{LixLy2006} we know that the development
map is certainly twice differentiable in the $1$-variation norm and applying
the chain rule it follows that $f$ is a twice differentiable function on $%
\left[ -1,1\right] $. On the other hand $f\left( 0\right) =0$ and $f\left(
\varepsilon \right) \leq 0$ for $\varepsilon \in \left[ -1,1\right] $ so
that $f^{\prime }\left( 0\right) =0$ and applying Taylor's theorem 
\begin{equation*}
0\geq f\left( 1\right) \geq \inf_{\varepsilon \in \left[ 0,1\right] }\frac{%
\varepsilon ^{2}}{2}f^{\prime \prime }\left( \varepsilon \right) .
\end{equation*}%
In fact the derivatives in $\varepsilon $ form a simple system of
differential equations. If 
\begin{equation*}
\left( 
\begin{array}{c}
\hat{x}_{t}^{\varepsilon +h} \\ 
x_{t}^{\varepsilon +h}%
\end{array}%
\right) =\left( 
\begin{array}{c}
\hat{x}_{t}^{\varepsilon } \\ 
x_{t}^{\varepsilon }%
\end{array}%
\right) +h\left( 
\begin{array}{c}
\hat{y}_{t}^{\varepsilon } \\ 
y_{t}^{\varepsilon }%
\end{array}%
\right) +\frac{h^{2}}{2}\left( 
\begin{array}{c}
\hat{z}_{t}^{\varepsilon } \\ 
z_{t}^{\varepsilon }%
\end{array}%
\right) +o\left( h^{2}\right) ,
\end{equation*}%
then 
\begin{eqnarray*}
\left( 
\begin{array}{c}
d\hat{x}_{t}^{\varepsilon } \\ 
d\hat{y}_{t}^{\varepsilon } \\ 
d\hat{z}_{t}^{\varepsilon }%
\end{array}%
\right) &=&\left( 
\begin{array}{ccc}
e^{\varepsilon A_{t}}vdt & 0 & 0 \\ 
A_{t}e^{\varepsilon A_{t}}vdt & e^{\varepsilon A_{t}}vdt & 0 \\ 
A_{t}^{2}e^{\varepsilon A_{t}}vdt & 2A_{t}e^{\varepsilon A_{t}}vdt & 
e^{\varepsilon A_{t}}vdt%
\end{array}%
\right) \left( 
\begin{array}{c}
x_{t}^{\varepsilon } \\ 
y_{t}^{\varepsilon } \\ 
z_{t}^{\varepsilon }%
\end{array}%
\right) \\
\left( 
\begin{array}{c}
dx_{t}^{\varepsilon } \\ 
dy_{t}^{\varepsilon } \\ 
dz_{t}^{\varepsilon }%
\end{array}%
\right) &=&\left( 
\begin{array}{ccc}
e^{\varepsilon A_{t}}vdt & 0 & 0 \\ 
A_{t}e^{\varepsilon A_{t}}vdt & e^{\varepsilon A_{t}}vdt & 0 \\ 
A_{t}^{2}e^{\varepsilon A_{t}}vdt & 2A_{t}e^{\varepsilon A_{t}}vdt & 
e^{\varepsilon A_{t}}vdt%
\end{array}%
\right) \bullet \left( 
\begin{array}{c}
\hat{x}_{t}^{\varepsilon } \\ 
\hat{y}_{t}^{\varepsilon } \\ 
\hat{z}_{t}^{\varepsilon }%
\end{array}%
\right)
\end{eqnarray*}%
with the initial conditions 
\begin{equation*}
\begin{array}{cc}
\hat{x}_{0}^{\varepsilon }=0 & x_{0}^{\varepsilon }=1 \\ 
\hat{y}_{0}^{\varepsilon }=0 & y_{0}^{\varepsilon }=0 \\ 
\hat{z}_{0}^{\varepsilon }=0 & z_{0}^{\varepsilon }=0.%
\end{array}%
\end{equation*}%
The simple exponential bound on the solution of a linear equation shows that 
\begin{equation*}
\left\vert z_{t}^{\varepsilon }\right\vert \leq 4^{\max \left\{ T,\int_{t\in 
\left[ 0,T\right] }\left\vert A_{t}\right\vert dt,\int_{t\in \left[ 0,T%
\right] }\left\vert A_{t}^{2}\right\vert dt\right\} .}
\end{equation*}%
Applying Taylor's theorem we have that 
\begin{equation*}
f\left( \varepsilon \right) \geq -\frac{\varepsilon ^{2}}{2}4^{\max \left\{
T,\int_{t\in \left[ 0,T\right] }\left\vert A_{t}\right\vert dt,\int_{t\in %
\left[ 0,T\right] }\left\vert A_{t}^{2}\right\vert dt\right\} }.
\end{equation*}%
If $\left\Vert A\right\Vert _{\infty }\leq 1$ then $f\left( \varepsilon
\right) >-\frac{\varepsilon ^{2}}{2}4^{T}$, and as $\left\Vert A\right\Vert
_{\infty }\leq \eta <1$, then we can replace $A$ by $\eta ^{-1}A$ and
evaluate $f_{\eta ^{-1}A}$ at $\eta $ to deduce that $f_{A}\left( 1\right) >-%
\frac{\eta ^{2}}{2}4^{T}$ giving us the uniform estimate we seek.
\end{proof}

\subsection{Some estimates from hyperbolic geometry}

We require some simple hyperbolic geometry. Fix $A$ (in hyperbolic space),
and consider two other points $B$ and $C.$ Let $\theta _{A}$, $\theta _{B}$,
and $\theta _{C}$ be the angles at $A$, $B$, and $C$ respectively. Let $a$, $%
b$, and $c$ be the hyperbolic lengths of the opposite sides. Recall the
hyperbolic cosine rule
\begin{equation*}
\sinh (b)\sinh (c)\cos (\theta _{A})=\cosh (b)\cosh (c)-\cosh (a)
\end{equation*}
and note the following simple lemmas:

\begin{lemma}
If the distance $c$ from $A$ to $B$ is at least $\log \left( \frac{\cos
\left| \theta _{A}\right| +1}{1-\cos \left| \theta _{A}\right| }\right) $,
then 
\begin{equation*}
\left| \theta _{B}\right| \leq \left| \theta _{A}\right| .
\end{equation*}
\end{lemma}

\begin{proof}
Fix $c$ and the angle $\theta _{A}$, the angle $\theta _{B}$ is zero if $b=0$
and monotone increasing as $b\rightarrow \infty $. Suppose that $\left\vert
\theta _{B}\right\vert >\left\vert \theta _{A}\right\vert $. We may reduce $%
b $ so that $\left\vert \theta _{B}\right\vert =\left\vert \theta
_{A}\right\vert $, now the triangle has two equal edges and applying the
cosine rule to compute the base length: 
\begin{eqnarray*}
\sinh (a)\sinh (c)\cos (\theta _{A}) &=&\cosh (a)\cosh (c)-\cosh (a) \\
c &=&\log \left( -\frac{\left( \cos \left\vert \theta _{A}\right\vert
\right) e^{2a}+e^{2a}-\cos \left\vert \theta _{A}\right\vert +1}{%
-e^{2a}+\left( \cos \left\vert \theta _{A}\right\vert \right) e^{2a}-\cos
\left\vert \theta _{A}\right\vert -1}\right) \\
&<&\lim_{a\rightarrow \infty }\log \left( -\frac{\left( \cos \left\vert
\theta _{A}\right\vert \right) e^{2a}+e^{2a}-\cos \left\vert \theta
_{A}\right\vert +1}{-e^{2a}+\left( \cos \left\vert \theta _{A}\right\vert
\right) e^{2a}-\cos \left\vert \theta _{A}\right\vert -1}\right) \\
&=&\log \left( \frac{\cos \left\vert \theta _{A}\right\vert +1}{1-\cos
\left\vert \theta _{A}\right\vert }\right) .
\end{eqnarray*}
\end{proof}

\begin{lemma}
\label{lem:mindist} We have $a\geq b+c - \log \frac{2}{1-\cos \theta _{A}}$.
and thus if $\max \left( b,c\right) \geq $ $\log \frac{2}{1-\cos \theta _{A}}
$, then $a>\min \left( b,c\right) $.
\end{lemma}

\begin{proof}
Suppose consider triangles with fixed angle $\theta _{A}$ and with side
lengths $\lambda b$, $\lambda c$ and resulting length $a\left( \lambda
\right) $ for the opposite. Then 
\begin{equation*}
\lambda b+\lambda c-a\left( \lambda \right)
\end{equation*}%
is monotone increasing in $\lambda $ with a finite limit. Now 
\begin{eqnarray*}
\sinh (\lambda b)\sinh (\lambda c)\cos (\theta _{A}) &=&\cosh (\lambda
b)\cosh (\lambda c)-\cosh (a\left( \lambda \right) ) \\
\frac{\cosh (\lambda b)\cosh (\lambda c)}{\sinh (\lambda b)\sinh (\lambda c)}%
-\cos (\theta _{A}) &=&\frac{\cosh (a\left( \lambda \right) )}{\sinh
(\lambda b)\sinh (\lambda c)} \\
\lim_{\lambda \rightarrow \infty }\log \frac{\cosh (a\left( \lambda \right) )%
}{\sinh (\lambda b)\sinh (\lambda c)} &=&\lim_{\lambda \rightarrow \infty
}\left( a\left( \lambda \right) -\lambda b-\lambda b\right) +\log 2 \\
\lambda b+\lambda c-a\left( \lambda \right) &\leq &\lim_{\lambda \rightarrow
\infty }\left( \lambda b+\lambda c-a\left( \lambda \right) \right) \\
&=&\log \frac{2}{1-\cos \theta _{A}}.
\end{eqnarray*}%
Thus 
\begin{equation*}
a\geq b+c-\log \frac{2}{1-\cos \theta _{A}}.
\end{equation*}%
Also, providing $\max \left( b,c\right) \geq $ $\log \frac{2}{1-\cos \theta
_{A}}$, one has $a\geq \min \left( b,c\right) $.
\end{proof}

\begin{corollary}
\label{cor:mindist} If the distance $c$ from $A$ to $B$ is at least $\log
\left( \frac{2}{1-\cos \left\vert \theta _{A}\right\vert }\right) $, then 
\begin{equation*}
\left\vert \theta _{B}\right\vert \leq \left\vert \theta _{A}\right\vert ,
\end{equation*}%
and $a\geq b$.
\end{corollary}

The above lemma is useful in the case where the angles of interest are
acute. But in some contexts we are interested in one angle is very obtuse in
which case the following lemma gives much better information.

\begin{lemma}
\label{cor:mindist2}Suppose that $\theta _{A}>\pi /2$ and that the distance $%
c$ from $A$ to $B$ is at least $\log \left( \sqrt{2}+1\right) $then $\theta
_{B}<\left( \pi -\theta _{A}\right) /2.$
\end{lemma}

\begin{proof}
The second hyperbolic cosine rule states that 
\begin{eqnarray*}
\sin (\theta _{B})\sin (\theta _{A})\cosh (c) &=&\cos (\theta _{C})+\cos
(\theta _{B})\cos (\theta _{A}) \\
\cosh (c) &=&\frac{\cos (\theta _{C})+\cos (\theta _{B})\cos (\theta _{A})}{%
\sin (\theta _{B})\sin (\theta _{A})}
\end{eqnarray*}%
Fix $\theta _{A}>\pi /2$. By our assumptions $\cosh \left( c\right) \geq 
\sqrt{2}$. and so 
\begin{equation*}
\frac{\cos (\theta _{C})+\cos (\theta _{B})\cos (\theta _{A})}{\sin (\theta
_{B})\sin (\theta _{A})}\geq \sqrt{2}.
\end{equation*}%
Since the sum of interior angles in a Hyperbolic triangle is less than $\pi $
one can conclude that $\theta _{B}=\alpha \left( \pi -\theta _{A}\right) $
where $0<\alpha <1$ and that $\theta _{B}$ and $\theta _{C}$ are in $\left[
0,\pi /2\right) $. To prove this lemma we need to show further, that $\alpha
\leq \frac{1}{2}$. It is enough to demonstrate that, in the case $\theta
_{A}>\pi /2,$ and $\frac{1}{2}<\alpha <1,$we have 
\begin{equation*}
\frac{\cos (\theta _{C})+\cos (\theta _{B})\cos (\theta _{A})}{\sin (\theta
_{B})\sin (\theta _{A})}<\sqrt{2}.
\end{equation*}

It is enough to prove that%
\begin{equation*}
\frac{1+\cos (\theta _{B})\cos (\theta _{A})}{\sin (\theta _{B})\sin (\theta
_{A})}<\sqrt{2}.
\end{equation*}%
Replacing $\left( \pi -\theta _{A}\right) $ by $\tau $ and rewriting 
\begin{equation*}
f\left( \alpha ,\tau \right) :=\frac{1-\cos (\alpha \tau )\cos (\tau )}{\sin
(\alpha \tau )\sin (\tau )}
\end{equation*}%
it is enough to prove that $f\left( \alpha ,\tau \right) <\sqrt{2}$ if $\tau
<\pi /2$ and $\frac{1}{2}<\alpha <1$. The derivative in $\alpha $ of $%
f\left( \alpha ,\tau \right) $ is%
\begin{equation*}
\mathbb{\tau }\frac{\left( \cos \left( \mathbb{\tau }\right) -\cos \left(
\alpha \mathbb{\tau }\right) \right) }{\sin \left( \mathbb{\tau }\right)
\sin \left( \alpha \mathbb{\tau }\right) ^{2}}
\end{equation*}%
and so $f$ is strictly decreasing in $\alpha $ in our domain. Hence, if $%
\alpha >1/2$ then 
\begin{equation*}
f\left( \alpha ,\tau \right) <f\left( \frac{1}{2},\tau \right)
\end{equation*}%
The derivative of $f\left( \frac{1}{2},\tau \right) $ is readily computed as 
\begin{equation*}
\frac{1}{8}\frac{\left( 1+2\cos \left( \mathbb{\tau }/2\right) \right) \tan
\left( \mathbb{\tau }/4\right) }{\cos \left( \mathbb{\tau }/4\right)
^{2}\cos \left( \mathbb{\tau }/2\right) ^{2}}
\end{equation*}%
and this is seen to be positive so that%
\begin{equation*}
f\left( \alpha ,\tau \right) <f\left( \frac{1}{2},\tau \right) <f\left( 
\frac{1}{2},\pi /2\right) =\sqrt{2}
\end{equation*}%
which completes the argument.
\end{proof}




\begin{lemma}
\label{spaced_out} Let $0=T_{0}<\ldots <T_{i}<\ldots T_{n}=T$ be a partition
of $\left[ 0,T\right] $. Let $\left( X_{t}\right) _{t\in \left[ 0,T\right] }$
be a continuous path, geodesic on the intervals $\left[ T_{i},T_{i+1}\right]
|_{i=0,\ldots ,n-1}$ in hyperbolic space with $n\geq 1$ where, at each $%
T_{i} $, the angle between the two geodesic segments: $\angle
X_{i-1}X_{T_{i}}X_{T_{i+1}}$ is in $\left[ 2\theta,\pi \right] $. Suppose
that each geodesic segment has length at least $K\left( \theta \right) =\log
\left( \frac{2}{1-\cos \left\vert \theta\right\vert }\right) $.

\begin{enumerate}
\item $d\left( X_{0},X_{T_{i}}\right) $ is increasing in $i$ and for each $%
i\leq n$ 
\begin{eqnarray}
d\left( X_{0},X_{T_{i}}\right) &\geq &d\left( X_{0},X_{T_{i-1}}\right)
+d\left( X_{T_{i-1}},X_{T_{i}}\right) -K\left( \theta \right)  \label{eq:dit}
\\
&\geq &K\left( \theta\right)  \notag
\end{eqnarray}%
and the angle between $\overrightarrow{X_{T_{i-1}},X_{T_{i}}}$ and $%
\overrightarrow{X_{0}X_{T_{i}}}$ is at most $\theta$.

\item We also have 
\begin{equation*}
0\leq \sum_{i=1}^{n}d\left( X_{T_{i-1}},X_{T_{i}}\right) -d\left(
X_{0},X_{T_{n}}\right) \leq (n-1) K(\theta).
\end{equation*}
\end{enumerate}
\end{lemma}

\begin{proof}
We proceed by induction. Suppose $d\left( X_{0},X_{T_{i}}\right) \geq
K\left( \theta \right) $ and the angle $\angle X_{T_{i-1}}X_{T_{i}}X_{T_{0}}$
is at most $\theta $. Now the angle $\angle X_{T_{i-1}}X_{T_{i}}X_{T_{i+1}}$
is at least $2\theta $ so that the angle $\angle
X_{T_{0}}X_{T_{i}}X_{T_{i+1}}$ is at least $\theta $. As $d\left(
X_{0},X_{T_{i}}\right) \geq K\left( \theta \right) $ and our supposition $%
d\left( X_{T_{i}},X_{T_{i+1}}\right) \geq K\left( \theta \right) $, Lemma~%
\ref{lem:mindist} and Corollary~\ref{cor:mindist} imply%
\begin{equation*}
d\left( X_{0},X_{T_{i+1}}\right) \geq d\left( X_{0},X_{T_{i}}\right)
+d\left( X_{T_{i}},X_{T_{i+1}}\right) -K\left( \theta \right) ,
\end{equation*}%
and that $\angle X_{0}X_{T_{i+1}}X_{T_{i}}\leq \theta $ proving the main
inequality. Using the induction one also has the second part of the
inequality%
\begin{equation*}
d(X_{0},X_{T_{i+1}})\geq K(\theta ).
\end{equation*}

The second claim is obtained by iterating (\ref{eq:dit}), 
\begin{equation*}
d\left( X_{0},X_{T_{i+1}}\right) \geq d\left( X_{0},X_{T_1}\right)
+\sum_{j=1}^i d\left( X_{T_{j}},X_{T_{j+1}}\right) -iK\left( \theta\right).
\end{equation*}
Now rearrange to get the result.
\end{proof}

\subsection{The main quantitative estimate}

Let $\gamma $ in $\mathbb{R}^{d}$ be a continuous path of finite length $l$,
and parameterised at unit speed. With this parameterisation $\dot{\gamma}$
can be regarded as a path on the unit sphere in $\mathbb{R}^{d}$. We
consider the case where $u\rightarrow \dot{\gamma}\left( u\right) $ is
continuous with modulus of continuity $\delta _{\gamma }$. If $\alpha \in 
\mathbb{R}$, then the path $\gamma _{\alpha }:=t\rightarrow \alpha \gamma
\left( t/\alpha \right) $ is also parameterised at unit speed, its length is 
$\alpha l$ and its derivative has modulus of continuity $\delta _{\gamma
_{\alpha }}\left( \alpha h\right) =\delta _{\gamma }\left( h\right) $. Its
development from the identity matrix (defined in (\ref{eq:development}))
into $SO\left( I_{d}\right) $ is denoted by $\Gamma _{\alpha }$.

The goal of this section is to provide a quantitative understanding for $%
\Gamma _{a}$ as we let $\alpha \rightarrow \infty $. Our estimates will only
depend on $\delta $ and the length of the path. We let $R_0=\log(1+\sqrt{2})$

\begin{proposition}
\label{lem:quantest} Let $\gamma $ in $\mathbb{R}^{d}$ be a continuous path
of length $l$. For each $C<1$ and $1\leq M\in {\mathbb{N}}$ then for any $%
\alpha $ chosen large enough that $\alpha l\geq MR_{0}$ and $\delta \left( 
\frac{M+1}{M}\frac{R_{0}}{\alpha }\right) <\sqrt{2\left( \sqrt{2}-\sqrt{%
1+C^{2}}\right) 4^{-\frac{M+1}{M}R_{0}}}$ one has 
\begin{equation*}
\left\vert d\left( o,\Gamma _{\alpha }o\right) -\alpha l\right\vert \leq
\left( \frac{4^{\frac{M+1}{M}R_{0}}}{2C}+\frac{16\log 2}{\pi ^{2}}\right) 
\frac{\alpha l}{R_{0}}\delta _{\gamma }\left( \frac{M+1}{M}\frac{R_{0}}{%
\alpha }\right) ^{2}
\end{equation*}%
and if $\alpha =1$, $l\geq MR_{0}$, providing $\delta \left( \frac{M+1}{M}%
R_{0}\right) <\sqrt{2\left( \sqrt{2}-\sqrt{1+C^{2}}\right) 4^{-\frac{M+1}{M}%
R_{0}}}$ we have for paths of any length%
\begin{equation*}
\left\vert d\left( o,\Gamma o\right) -l\right\vert \leq \left( \frac{4^{%
\frac{M+1}{M}R_{0}}}{2C}+\frac{16\log 2}{\pi ^{2}}\right) \frac{l}{R_{0}}%
\delta _{\gamma }\left( \frac{M+1}{M}R_{0}\right) ^{2}.
\end{equation*}
\end{proposition}

We set $D_{1}\left( C,M\right) =\left( \frac{4^{\frac{M+1}{M}R_{0}}}{2C}+%
\frac{16\log 2}{\pi ^{2}}\right) \frac{l}{R_{0}}$, and $D_{2}\left( M\right)
=\frac{M+1}{M}R_{0}$ so that the inequality becomes 
\begin{equation}
\left\vert d\left( o,\Gamma _{\alpha }o\right) -l\alpha \right\vert \leq
D_{1}\delta _{\gamma }\left( D_{2}/\alpha \right) ^{2}\alpha l.
\end{equation}


We note that $R_{0}\approx .881374$, $4^{R_0} \approx 3.34393 \leq
4^{(M+1)R_{0}/M} \leq 4^{2R_{0}}\approx 11.5154$, $\frac{16\log 2}{\pi ^{2}}%
\approx 1.12369.$ Fixing $M=1,$ one immediately sees that the distance $%
d\left( o,\Gamma _{\alpha }o\right) $ grows linearly with the scaling and
the chordal distance $d\left( o,\Gamma _{\alpha }o\right) $ behaves like the
length of the path $\gamma $ as $\alpha \rightarrow \infty $. We also note
that the shape of this result is reminiscent of the elegant result of Fawcett%
\cite[Lemma 68]{fawcett}: that among $C^{2}-$curves $\gamma $ with modulus
of continuity $\delta _{\gamma }\left( h\right) \leq \kappa h$ one has sharp
estimates on the minimal value of $d\left( o,\Gamma o\right) $ given by 
\begin{equation*}
\inf_{\gamma }\cosh \left( d\left( o,\Gamma o\right) \right) =\frac{\cosh
\left( \alpha l\sqrt{1-\kappa ^{2}}\right) -\kappa ^{2}}{1-\kappa ^{2}}.
\end{equation*}%
A natural question to ask is whether our estimate (which is
non-infinitesimal and only needs information about $\delta _{\gamma }\left(
2R_{0}\right) $) can be improved to this shape and even to this sharp form.

\begin{proof}
The path $\gamma _{\alpha }:=t\rightarrow \alpha \gamma \left( t/\alpha
\right) $ is of length $\alpha l$ and parameterised at unit speed; its
derivative has modulus of continuity $\delta _{\alpha }:t\rightarrow \delta
_{\gamma }\left( t/\alpha \right) $. Because $\alpha l\geq MR_{0}$ we can
fix $R=\alpha l/N$, where $R\in \left[ R_{0},\frac{M+1}{M}R_{0}\right] $ and 
$N$ is a positive integer depending on $\alpha $. Let $t_{i}=iR$ where $i\in %
\left[ 0,N\right] $. Let $G_{i}\in SO\left( I_{d}\right) $ be the
development of the path segment $\gamma _{\alpha }|_{\left[ t_{i-1},t_{i}%
\right] }$ into $SO\left( I_{d}\right) $ and $\Gamma _{\alpha ,t}$ be the
development of the path segment $\gamma _{\alpha }|_{\left[ 0,t\right] }$.
We define $X_{0}:=o\in \mathbb{H}$ and $X_{j}:=G_{t_{j}}X_{j-1}\in \mathbb{H}
$. Then $X_{j}$ are the points $\Gamma _{\alpha ,t_{j}}o$ on the path $%
\Gamma _{a,t}o.$


As the length of the path is greater than any chord%
\begin{eqnarray}
\left\vert \alpha l-d\left( o,\Gamma _{\alpha }o\right) \right\vert
&=&\alpha l-\sum_{i=1}^{N}d\left( X_{i-1},X_{i}\right)  \notag \\
&&+\sum_{i=1}^{N}d\left( X_{i-1},X_{i}\right) -d\left( X_{0},X_{N}\right) .
\label{eq:dgam}
\end{eqnarray}%
and
\begin{equation*}
\alpha l-\sum_{i=1}^{N}d\left( X_{i-1},X_{i}\right) \geq 0 \;\;\;\;
\sum_{i=1}^{N}d\left( X_{i-1},X_{i}\right) -d\left(
X_{0},X_{N}\right) \geq 0.
\end{equation*}
We now estimate each of these terms from above. %
%

For the first term we use our result on paths close to a geodesic. By
Proposition~\ref{prop:geod} we have 
\begin{eqnarray*}
\cosh d\left( X_{i-1},X_{i}\right) &\geq &\cosh R-\frac{\delta _{\gamma
_{\alpha }}(R)^{2}}{2}4^{R} \\
&=&\cosh R-\frac{\delta _{\gamma }(\frac{R}{\alpha })^{2}}{2}4^{R}
\end{eqnarray*}%
Thus, using the convexity of cosh and hyperbolic trig identities, 
\begin{eqnarray*}
\frac{\delta _{\gamma }(\frac{R}{\alpha })^{2}}{2}4^{R} &\geq &\cosh R-\cosh
d\left( X_{i-1},X_{i}\right) \\
&\geq &\left( R-d\left( X_{i-1},X_{i}\right) \right) \sinh d\left(
X_{i-1},X_{i}\right) \\
&=&\left( R-d\left( X_{i-1},X_{i}\right) \right) \sqrt{\cosh d\left(
X_{i-1},X_{i}\right) ^{2}-1} \\
&\geq &\left( R-d\left( X_{i-1},X_{i}\right) \right) \sqrt{\left( \cosh R-%
\frac{\delta _{\gamma }(\frac{R}{\alpha })^{2}}{2}4^{R}\right) ^{2}-1} \\
&\geq &\left( R-d\left( X_{i-1},X_{i}\right) \right) \sqrt{\left( \sqrt{2}-%
\frac{\delta _{\gamma }(\frac{R}{\alpha })^{2}}{2}4^{R}\right) ^{2}-1}.
\end{eqnarray*}%
Now, for $C<1$, providing 
\begin{equation*}
\left( \sqrt{2}-\frac{\delta _{\gamma }(\frac{R}{\alpha })^{2}}{2}%
4^{R}\right) ^{2}\geq 1+C^{2},
\end{equation*}%
we have 
\begin{equation*}
\frac{\delta _{\gamma }(\frac{R}{\alpha })^{2}}{2C}4^{R}\geq \left(
R-d\left( X_{i-1},X_{i}\right) \right) .
\end{equation*}%
This will follow if we choose $\alpha $ large enough such that our
condition 
\begin{equation*}
\delta _{\gamma }(\frac{M+1}{M}\frac{R_{0}}{\alpha })^{2}\leq 2\left( \sqrt{2%
}-\sqrt{1+C^{2}}\right) 4^{-\frac{M+1}{M}R_{0}},
\end{equation*}%
holds.

Hence, summing over all the pieces, we have 
\begin{equation*}
\left( \frac{\alpha l}{R}\right) \frac{\delta _{\gamma }(\frac{R}{\alpha }%
)^{2}}{2C}4^{R}\geq \left( \alpha l-\sum_{i=1}^{N}d\left(
X_{i-1},X_{i}\right) \right) .
\end{equation*}%
We now use our bounds on $R$ to obtain 
\begin{equation}
\left( \frac{\alpha l}{R_{0}}\right) \frac{\delta _{\gamma }(\frac{M+1}{M}%
\frac{R_{0}}{\alpha })^{2}}{2C}4^{\frac{M+1}{M}R_{0}}\geq \left( \alpha
l-\sum_{i=1}^{N}d\left( X_{i-1},X_{i}\right) \right) .  \label{eq:partia}
\end{equation}

Applying Lemma~\ref{cor:mindist2} and Lemma~\ref{spaced_out} we see that as $%
R\geq R_{0}=\log (1+\sqrt{2})$, then 
the angle $X_{0}X_{T_{n}}X_{T_{n+1}}$ is at least $\pi -2\delta _{\gamma
}\left( \frac{R}{\alpha }\right) $ for each $0<n<N-1$. 

Thus 
\begin{eqnarray*}
\sum_{i=1}^{N}d\left( X_{i-1},X_{i}\right) -d\left( X_{0},X_{T_{N}}\right)
&\leq &(N-1)K(\pi -2\delta _{\gamma }(\frac{R}{\alpha })) \\
&=&\left( \frac{\alpha l}{R}-1\right) \log \left( \frac{2}{1-\cos \left( \pi
-2\delta _{\gamma }\left( \frac{R}{\alpha }\right) \right) }\right) .
\end{eqnarray*}%
Since $\left( \frac{1}{u^{2}}\log \left( \frac{2}{1-\cos \left\vert \pi
-u\right\vert }\right) \right) $ is increasing in $u$ for $u\leq \pi $, we
have 
\begin{eqnarray*}
\log \left( \frac{2}{1-\cos \left\vert \pi -2\delta _{\gamma }\left( \frac{R%
}{\alpha }\right) \right\vert }\right) &\leq &4\delta _{\gamma }\left( \frac{%
R}{\alpha }\right) ^{2}\left( \frac{1}{\left( \pi /2\right) ^{2}}\log \left( 
\frac{2}{1-\cos \left\vert \pi /2\right\vert }\right) \right) \\
&\leq &\delta _{\gamma }\left( \frac{R}{\alpha }\right) ^{2}\frac{16\log 2}{%
\pi ^{2}}\ \ \ \text{since} \\
\ \delta _{\gamma }\left( \frac{R}{\alpha }\right) &\leq &\sqrt{2\left( 
\sqrt{2}-\sqrt{1+C^{2}}\right) 4^{-\frac{M+1}{M}R_{0}}}
\end{eqnarray*}%
which is less than $\pi /4$ for all $C$ and $M.$

Observing that $\delta _{\gamma }(R)$ is increasing and that $\frac{M+1}{M}%
R_{0}\geq R$ gives the second part of our estimate 
\begin{eqnarray}
\sum_{i=1}^{n}d\left( X_{i-1},X_{i}\right) -d\left( X_{0},X_{T_{n}}\right)
&\leq &\frac{16\log 2}{\pi ^{2}}\left( \frac{\alpha l}{R}-1\right) \delta
_{\gamma }\left( \frac{R}{\alpha }\right) ^{2}. \\
&\leq &\frac{16\log 2}{\pi ^{2}}\left( \frac{\alpha l}{R_{0}}-1\right)
\delta _{\gamma }\left( \frac{M+1}{M}\frac{R_{0}}{\alpha }\right) ^{2}.
\label{eq:partib}
\end{eqnarray}

Combining the estimates (\ref{eq:partia}) and (\ref{eq:partib}) completes
proof.
\end{proof}

\subsection{Recovering the length of the path from its signature}

From the last section we know that if $\alpha $ is large enough then 
\begin{equation}
\left\vert d\left( o,\Gamma _{\alpha }o\right) -l\alpha \right\vert \leq
D_{1}\delta _{\gamma }\left( D_{2}/\alpha \right) ^{2}\alpha l
\label{eq:near_linear}
\end{equation}%
and in particular will go to zero as $\alpha \rightarrow \infty $ if $\delta
_{\gamma }\left( \varepsilon \right) =o\left( \varepsilon ^{1/2}\right) $.

The lower bound on $d\left( o,\Gamma _{\alpha }o\right) $ implicit in (\ref%
{eq:near_linear}) leads to a lower bound on the norm of $\Gamma _{\alpha }$
as a matrix. We will compare it with the upper bound that comes from
expressing the matrix $\Gamma _{\alpha }$ as a series whose coefficients are
iterated integrals. We have an upper bound for each coefficient in the
series, and taken together these provide a bound for the sum. This bound is
so close to the lower bound that it allows us to conclude a lower bound for
each coefficient and relate the decay rate for the norms of the iterated
integrals directly to the length of $\gamma $.

It is an open question as to whether signatures with given decay rate
correspond to paths of finite length. 

\begin{proposition}
\label{prop:lbg}Let $G\in SO\left( I_{d}\right) .$ Then $\left\Vert
G\right\Vert \geq e^{d\left( o,Go\right) }$ where $\left\Vert G\right\Vert $%
is the operator norm for $G\in Hom\left( \mathbb{R}^{d+1},\mathbb{R}%
^{d+1}\right) $ where $\mathbb{R}^{d+1}$ has the Euclidean norm.
\end{proposition}

\begin{proof}
If%
\begin{equation*}
F_{\rho }:=\left( 
\begin{array}{ccccccc}
1 & 0 & \cdots & \cdots & 0 & 0 & 0 \\ 
0 & 1 & \ddots & \ddots & \vdots & \vdots & \vdots \\ 
\vdots & \ddots & \ddots & \ddots & \vdots & \vdots & \vdots \\ 
\vdots & \ddots & \ddots & 1 & 0 & \vdots & \vdots \\ 
0 & \cdots & \cdots & 0 & 1 & 0 & 0 \\ 
0 & \cdots & \cdots & \cdots & 0 & \cosh \rho & \sinh \rho \\ 
0 & \cdots & \cdots & \cdots & 0 & \sinh \rho & \cosh \rho%
\end{array}%
\right),
\end{equation*}%
then $F_{\rho }F_{\tau }=F_{\rho +\tau }$ and the set of such elements forms
a (maximal) abelian subgroup of $SO\left( I_{d}\right) $. Any element $G$ of 
$SO\left( I_{d}\right) $ can be factored into a Cartan Decomposition $%
KF_{\rho }\tilde{K}$ where $K$ and $\tilde{K}$ are built out of rotations $%
\Theta $ of $\mathbb{R}^{d}$%
\begin{equation*}
\left( 
\begin{array}{cc}
\Theta & \mathbf{0} \\ 
\mathbf{0}^{t} & 1%
\end{array}%
\right)
\end{equation*}%
and $\rho \in \mathbb{R}_{+}.$ As an operator on Euclidean space, $G$ has
norm $\Vert G\Vert =\left\Vert KF_{\rho }\tilde{K}\right\Vert =$ $\left\Vert
F_{\rho }\right\Vert $ since $K$, $\tilde{K}$ are isometries. In addition,
the matrix $F_{\rho }$ is symmetric and hence has a basis comprising
eigenfunctions; its norm is at least as large as its largest eigenvalue.
Computation shows that the eigenvalues of $F_{\rho }$ are $\left\{ e^{\rho
},e^{-\rho },1,\cdots ,1\right\} $ so that, given $\rho >0$, one has 
\begin{equation*}
\left\Vert G\right\Vert \geq e^{\rho }.
\end{equation*}%
On the other hand 
\begin{eqnarray*}
-\cosh d\left( o,Go\right) &=&I_{d}\left( \left( 
\begin{array}{c}
x_{1} \\ 
\vdots \\ 
x_{d} \\ 
\cosh \rho%
\end{array}
\right) ,\left( 
\begin{array}{c}
0 \\ 
\vdots \\ 
0 \\ 
1%
\end{array}
\right) \right) \\
&=&-\cosh \rho
\end{eqnarray*}%
and so $\left\Vert G\right\Vert \geq e^{d\left( o,Go\right) }.$
\end{proof}

If $\gamma $ is a path of finite length then the development (\ref%
{eq:development}) into hyperbolic space $\mathbb{H}$ is defined by 
\begin{equation*}
d\Gamma _{t}=F\left( d\gamma _{t}\right) \Gamma _{t},
\end{equation*}%
where by Lemma~\ref{lem:Fnorm} $F:\mathbb{R}^{d}\rightarrow \hom \left( 
\mathbb{R}^{d+1},\mathbb{R}^{d+1}\right) $ has norm one as a map from
Euclidean space to the operators on Euclidean space. As a result the
development of a path $X$ is given by 
\begin{equation*}
G=I+\int_{0<u<T}F\left( dX_{u}\right) +\ldots +\int_{0<u_{1}<\ldots
<u_{k}<T}F\left( dX_{u_{1}}\right) \otimes \ldots \otimes F\left(
dX_{u_{k}}\right) +\ldots
\end{equation*}%
and, as in Lemma~\ref{lem:intpathapriori}, if $X$ is a path of length $\theta $,
then we have an a priori bound 
\begin{equation*}
\left\Vert \int_{0<u_{1}<\ldots <u_{k}<T}dX_{u_{1}}\otimes \ldots \otimes
dX_{u_{k}}\right\Vert \leq \frac{\theta ^{n}}{n!}.
\end{equation*}%
Applying this to $\alpha \gamma $, we conclude that 
\begin{eqnarray*}
e^{d\left( o,\Gamma _{\alpha }o\right) } &\leq &\left\Vert \Gamma _{\alpha
}\right\Vert \\
&\leq &\sum \alpha ^{n}\left\Vert \int_{0<u_{1}<\ldots <u_{k}<1}F\left(
d\gamma _{u_{1}}\right) \otimes \ldots \otimes F\left( d\gamma
_{u_{k}}\right) \right\Vert \\
&\leq &\sum \alpha ^{n}\left\Vert \int_{0<u_{1}<\ldots <u_{k}<1}d\gamma
_{u_{1}}\otimes \ldots \otimes d\gamma _{u_{k}}\right\Vert \\
&\leq &e^{\alpha l},
\end{eqnarray*}%
where $l$ is the length of $\gamma $. Letting 
\begin{equation*}
b_{n}=n!\left\Vert \int_{0<u_{1}<\ldots <u_{k}<1}d\gamma _{u_{1}}\otimes
\ldots \otimes d\gamma _{u_{k}}\right\Vert
\end{equation*}%
one has for all $\alpha $ that 
\begin{eqnarray*}
e^{d\left( o,\Gamma _{\alpha }o\right) -\alpha l} &\leq &e^{-\alpha
l}\sum_{k=0}^{\infty }\frac{\alpha ^{k}}{k!}b_{k}\leq 1 \\
0 &\leq &b_{k}\leq l^{k}.
\end{eqnarray*}%
Thus the expectation of $b_{n}$ with respect to a Poisson measure with mean $%
\alpha l$ is close to one while at the same time the $b_{n}$ are all bounded
above by one and positive. In particular%
\begin{eqnarray*}
\sum_{k=0}^{\infty }\frac{\alpha ^{k}}{k!}\left\vert l^{k}-b_{k}\right\vert
&\leq &e^{\alpha l}-e^{d\left( o,\Gamma _{\alpha }o\right) } \\
&\leq &e^{\alpha l}\left( 1-e^{-D_{1}\delta \left( D_{2}/\alpha \right)
^{2}\alpha l}\right)
\end{eqnarray*}%
and so 
\begin{equation*}
\left\vert l^{k}-b_{k}\right\vert \leq \inf_{\alpha >1}k!\alpha
^{-k}e^{\alpha l}\left( 1-e^{-D_{1}\delta \left( D_{2}/\alpha \right)
^{2}\alpha l}\right)
\end{equation*}%
applying Stirling's formulae that $k!=e^{k\log k-k+\frac{1}{2}\log k+C_{k}}$
where $C_{k}=o\left( 1\right) $ and setting $\alpha =k/l$ gives 
\begin{eqnarray*}
\left\vert l^{k}-b_{k}\right\vert &\leq &e^{C_{k}}l^{k}\sqrt{k}\left(
1-e^{-D_{1}\delta \left( D_{2}l/k\right) ^{2}k}\right) \\
&\leq &l^{k}\tilde{C}\delta \left( lD_{2}/k\right) ^{2}k\sqrt{k}
\end{eqnarray*}%
where $\tilde{C}=D_{1}e^{C_{k}}$ and so we see that, if $\delta _{\gamma
}\left( lD_{2}/k\right) ^{2}k^{3/2}\rightarrow 0$ as $k\rightarrow \infty $,
then $b_{k}/l^{k}\rightarrow 1$. Thus we have shown the following

\begin{theorem}\label{thm:strongrecover}
For any path of finite length with $\delta _{\gamma }\left( \varepsilon
\right) =o\left( \varepsilon ^{3/4}\right) $, 
\begin{equation*}
l^{-k}k!\left\Vert \int_{0<u_{1}<\ldots <u_{k}<1}d\gamma _{u_{1}}\otimes
\ldots \otimes d\gamma _{u_{k}}\right\Vert \rightarrow 1,
\end{equation*}
as $k\to\infty$,
\end{theorem}

This is of course quite a strong result obtained by making strong
assumptions. One could ask less and so we give a weaker but more widely
applicable result.

\begin{theorem}
\label{Thm:recoverlength}Let $\gamma $ be a path of finite length $l$, and
suppose its derivative, when parameterised at unit speed, is continuous.
Then the Poisson averages $C_{\alpha }$ of the $b_{k}$ defined by 
\begin{equation*}
C_{\alpha }=e^{-\alpha }\sum_{k=0}^{\infty }\frac{\alpha ^{k}}{k!}b_{k}
\end{equation*}%
satisfy 
\begin{equation*}
\lim_{\alpha \rightarrow \infty }\frac{1}{\alpha }\log C_{\alpha }=l-1
\end{equation*}
\end{theorem}

Note that the $C_{\alpha }$ are averages of the $b_{k}$ against Poisson
measures; it is standard that these are close to Gaussian with mean $\alpha $
and variance $\alpha$.

\begin{proof}
Note that 
\begin{equation*}
e^{d\left( o,\Gamma _{\alpha }o\right) -\alpha l}\leq e^{-\alpha
l}\sum_{k=0}^{\infty }\frac{\alpha ^{k}}{k!}b_{k}\leq 1
\end{equation*}%
and so 
\begin{equation}
\frac{d\left( o,\Gamma _{\alpha }o\right) }{\alpha }-l\leq \frac{1}{\alpha }%
\log \left( e^{-\alpha }\sum_{k=0}^{\infty }\frac{\alpha ^{k}}{k!}%
b_{k}\right) +1-l\leq 0  \label{eq:sandwich}
\end{equation}%
and using (\ref{eq:near_linear}) we have 
\begin{equation*}
\left\vert \frac{d\left( o,\Gamma _{\alpha }o\right) }{\alpha }-l\right\vert
\leq D_{1}\delta _{\gamma }\left( D_{2}/\alpha \right) ^{2}l
\end{equation*}%
and so the left hand side in (\ref{eq:sandwich}) goes to zero.
\end{proof}

In particular we see that the high order coefficients of the signature
already determine the length of the path and in fact one can obtain
quantitative estimates in terms of the modulus of continuity for the
derivative of $\gamma $.


\begin{conjecture}
The length of $\tilde{\gamma}$ can be recovered from the asymptotic
behaviour of averages of the $b_{k}$.
\end{conjecture}

It might be that $\lim_{\alpha \rightarrow \infty }1+\frac{1}{\alpha }\log
C_{\alpha }$ gives the length of $\tilde{\gamma}$ directly although
the Poisson averages may have to be replaced in some way.

We conclude with an analogous result to that proved for the lattice case in
Proposition~\ref{prop:intlat}.

\begin{theorem}
Let $\gamma $ be a path of length $l$ parameterised at unit speed, and let $%
\delta _{\gamma }$ be the modulus of continuity for $\dot{\gamma}$. Fix $C<1$
and $1\leq M\in \mathbb{N}$. 
Suppose that $\delta _{\gamma }\left( 0\right) <\frac{1}{\sqrt{D_{1}\left(
C,M\right) }}$, then there is an integer $N\left( l,\delta \right) $ such
that at least one of the first $N(l,\delta)$ terms in the signature must be non-zero.
\end{theorem}

\begin{proof}
In the case where the first $e\alpha l$ coefficients in the signature of the
path $\gamma $ are zero, by Lemma~\ref{lem:stirling}, we have some explicit constant
$C_{1}$ such that 
\begin{eqnarray*}
\left\Vert \Gamma _{\alpha }\right\Vert &\leq &1+\sum_{m>e\alpha l}\frac{%
\left( \alpha l\right) ^{m}}{m!} \\
&\leq &1+C_{1}(\alpha l)^{-1/2}.
\end{eqnarray*}%
By Proposition \ref{prop:lbg}, and letting $\alpha$ be sufficiently large
so that we can apply (\ref{eq:near_linear}), 
\begin{eqnarray*}
\left\Vert \Gamma _{\alpha }\right\Vert &\geq &e^{d\left( o,\Gamma _{\alpha
}o\right) } \\
&\geq &e^{l\alpha -D_{1}\delta _{\gamma }\left( D_{2}/\alpha \right)
^{2}\alpha l} \\
&\geq &1+l\alpha -D_{1}\delta _{\gamma }\left( D_{2}/\alpha \right)
^{2}\alpha l.
\end{eqnarray*}%
These two statements lead to a contradiction if for large $\alpha$, we have $l\alpha
-D_{1}\delta _{\gamma }\left( D_{2}/\alpha \right) ^{2}\alpha l>C_{1}(\alpha
l)^{-1/2}$ or 
\begin{equation*}
\alpha ^{3/2}\left( 1-D_{1}\delta _{\gamma }\left( D_{2}/\alpha \right)
^{2}\right) >C_{1}l^{-3/2}.
\end{equation*}%
Thus providing $1>D_{1}\delta _{\gamma }\left( D_{2}/\alpha \right) ^{2}$ for
some large $\alpha $ (continuity of the derivative is enough) then the left
hand side goes to infinity as $\alpha \rightarrow \infty $. This
always gives a contradiction and shows that existence of $N(l,\delta)$.
\end{proof}

An explicit estimate is $N(l,\delta) = \lceil e\alpha l\rceil$, where
one chooses the smallest $\alpha \geq \frac{MR_{0}}{l}$, large enough so
that 
\begin{equation*}
\delta _{\gamma }\left( \frac{M+1}{M}\frac{R_{0}}{\alpha }\right) <\sqrt{%
2\left( \sqrt{2}-\sqrt{1+C^{2}}\right) 4^{-\frac{M+1}{M}R_{0}}}
\end{equation*}%
and so that $\alpha ^{3/2}\left( 1-D_{1}\delta _{\gamma }\left( D_{2}/\alpha
\right) ^{2}\right) >C_{1}l^{-3/2}$.
To give an idea of the numerical size of $N$, the number of iterated integrals
required for this result, we note that if the path has $\delta (h)\leq h$,
then the optimal value of $\alpha $ is around 15.2 at $C=0.8875$ (with
$M=1$) and the
number $N$ is the integer greater than $41.38l$ (for large $l$). With a more
careful optimization of the constants (varying $M)$ our estimate can be
reduced to the integer greater than $13.28l$ (for large $l$).

\begin{remark}
\textrm{We note that this proof did not require that $\delta _{\gamma
}\left( 0\right) =0$ or that $\dot{\gamma}$ is continuous.}
\end{remark}

\begin{remark}
\textrm{An easy way to produce a path with each of the first $N$ iterated
integrals zero is take two paths with the same signature up to the level of
the $N$'th iterated integral and to take the first path concatenated with
the second with time run backwards. Since these paths will, except at the
point of joining, have the same smoothness as they did before, all focus
goes to the point where they join. One could hope that a development of
these ideas would prove that the two paths must be nearly tangential. If
this were exactly true, then it would give a reconstruction theorem.}
\end{remark}

We have obtained quantitative lower bounds on the signature of $\gamma $\
when $\gamma $ is parameterised at unit speed and $\dot{\gamma}$ is close to
continuously differentiable. In fact one could obtain estimates whenever the 
$\dot{\gamma}$ is piecewise continuous and the jumps are less than $\pi $.
However the main extra idea is already visible in the case where $\dot{\gamma}$
is piecewise constant. We give an explicit estimate in
Theorem~\ref{cor:nondeg} in Section~\ref{section6}.

\section{Tree-Like paths}

\label{SectiontreeLike}\label{Section4}

We now turn to our proof of the extension of Chen's theorem to the case of
finite length paths. In this section we suppose that $X_{t\in \left[ 0,T%
\right] }$ is a path in a Banach or metric space $E$ and we recall our
definition~\ref{def:tree-like} of tree-like paths in this more general
setting.

\begin{theorem}
\label{prop:h}If $X$ is a tree-like path with height function $h$ and, if $X$
is of bounded variation, then there exists a new height function $\tilde{h}$
having bounded variation and hence $X$ is a \emph{Lipschitz} tree like path;
moreover, the variation of $\tilde{h}$ is bounded by the variation of $X$.
\end{theorem}

\begin{proof}
\ The function $h$ allows one to introduce a partial order and tree
structure on $\left[ 0,T\right] $. Let $t\in \left[ 0,T\right] $. Define the
continuous and monotone function $g_{t}\left( .\right) $ by 
\begin{equation*}
g_{t}\left( v\right) =\inf_{v\leq u\leq t}h\left( u\right) ,\ v\in \left[ 0,t%
\right] .
\end{equation*}%
The intermediate value theorem ensures that $g_{t}$ maps $\left[ 0,t\right] $
onto $\left[ 0,h\left( t\right) \right] $. 
Let $\tau _{t}$ be a maximal inverse of $h$ in that 
\begin{equation}
\tau _{t}\left( x\right) =\sup \left\{ u\in \left[ 0,t\right] |g_{t}\left(
u\right) =x\right\} ,\;\;x\in \lbrack 0,h(t)].  \label{eq:sup}
\end{equation}%
As $\ g_{t}$ is monotone and continuous 
\begin{equation}
\tau _{t}\left( x\right) =\inf \left\{ u\in \left[ 0,t\right] |g_{t}\left(
u\right) >x\right\}   \label{eq:inf}
\end{equation}%
for $x<h\left( t\right) $.

Now say $s\preceq t$ if and only if $\ s$ is in the range of $\tau _{t}$;
that is to say if there is an $x\in \left[ 0,h\left( t\right) \right] $ so
that $s=\tau _{t}\left( x\right) $. Since $\tau _{t}\left( h\left( t\right)
\right) =\sup \left\{ u\in \left[ 0,t\right] |g_{t}\left( u\right) =h\left(
t\right) \right\} $, it follows that $\tau _{t}\left( h\left( t\right)
\right) =t$ and so $t\preceq t$. Since $h\left( \tau _{t}\left( x\right)
\right) =x$ \ for \ $x\in \left[ 0,h\left( t\right) \right] $ we see there
is an inequality-preserving bijection between the $\left\{ s|s\preceq
t\right\} $ and $\left[ 0,h\left( t\right) \right] $.

Suppose $t_{1}\preceq t_{0}$ and that they are distinct; then $h\left(
t_{1}\right) <h\left( t_{0}\right) $. We may choose $x_{1}\in \left[
0,h\left( t_{0}\right) \right) $ so that $t_{1}=\tau _{t_{0}}\left(
x_{1}\right) $, it follows that 
\begin{eqnarray*}
t_{1} &=&\tau _{t_{0}}\left( x_{1}\right)  \\
&=&\inf \left\{ u\in \left[ 0,t_{0}\right] |g_{t_{0}}\left( u\right)
>x_{1}\right\} ,
\end{eqnarray*}%
and that 
\begin{equation*}
h\left( t_{1}\right) =x_{1}<h\left( u\right) ,~u\in \left( t_{1},t_{0}\right]
.
\end{equation*}%
Of course 
\begin{eqnarray*}
g_{t_{0}}\left( t_{1}\right)  &=&\inf_{t_{1}\leq u\leq t_{0}}h\left(
u\right)  \\
&=&h\left( t_{1}\right)  \\
&=&g_{t_{1}}\left( t_{1}\right) ,
\end{eqnarray*}%
and hence $g_{t_{0}}\left( u\right) =g_{t_{1}}\left( u\right) $ for all $%
u\in \left[ 0,t_{1}\right] $. Hence, $\tau _{t_{0}}\left( x\right) =\tau
_{t_{1}}\left( x\right) $ for any $x<g_{t_{1}}\left( t_{1}\right) =h\left(
t_{1}\right) =x_{1}$; we have already seen that $\tau _{t_{1}}\left( h\left(
t_{1}\right) \right) =t_{1}=\tau _{t_{0}}\left( x_{1}\right) $. It follows
that the range $\tau _{t_{1}}\left( \left[ 0,h\left( t_{1}\right) \right]
\right) $ is contained in the range of $\tau _{t_{0}}$. In particular, we
deduce that if $t_{2}\preceq t_{1}$ and $t_{1}\preceq t_{0}$ then $%
t_{2}\preceq t_{0}$.

We have shown that $\preceq $ is a partial order, and that $\left\{
t|t\preceq t_{0}\right\} $ is totally ordered under $\preceq $, and in one
to one correspondence with $\left[ 0,h\left( t_{0}\right) \right] $.

Now, consider two generic times $s<t$. Let $x_{0}=\inf_{s\leq u\leq
t}h\left( u\right) $ and $I=\left\{ v\in \left[ s,t\right] |h\left( v\right)
=x_{0}\right\} $. Since $h$ is continuous and $\left[ s,t\right] $ is
compact the set $I$ is non-empty and compact. 
By the construction of the function $g_{t}$ it is obvious that $g_{t}\leq
g_{s}$ on $\left[ 0,s\right] $ and that if $g_{t}\left( u\right)
=g_{s}\left( u\right) $, then $g_{t}\left( v\right) =g_{s}\left( v\right) $
for $v\in \left[ 0,u\right] $. Thus, there will be a unique $r\in \left[ 0,s%
\right] $ so that $g_{s}=g_{t}$ on $\left[ 0,r\right] $ and $g_{t}<g_{s}$ on 
$\left( r,s\right] $. Observe that $g_{t}\left( r\right) =x_{0}$ and that $%
\tau _{t}\left( x_{0}\right) =\sup I$ and, essentially as above $\tau
_{s}=\tau _{t}$ on $\left[ 0,h\left( r\right) \right) $. Observe also that
if $\tilde{t}\in \left[ s,t\right] $ then $g_{s}=g_{\tilde{t}}$ on $\left[
0,r\right] $ so that $\tau _{s}=\tau _{\tilde{t}}$ on $\left[ 0,h\left(
r\right) \right) $.

Having understood $h$ and $\tau $ to the necessary level of detail, we
return to the path $X$. For $x$, $y\in \left[ 0,h\left( t\right) \right] $
one has, for $x<y$, 
\begin{eqnarray*}
\left\Vert X_{\tau _{t}\left( x\right) }-X_{\tau _{t}\left( y\right)
}\right\Vert &\leq &h\left( \tau _{t}\left( x\right) \right) +h\left( \tau
_{t}\left( y\right) \right) -2\inf_{u\in \left[ \tau \left( x\right) ,\tau
\left( y\right) \right] }h\left( u\right) \\
&\leq &x+y-2\inf_{z\in \left[ x,y\right] }h\left( \tau _{t}\left( z\right)
\right) \\
&=&y-x
\end{eqnarray*}
so we see that $X_{\tau _{t}\left( .\right) }$ is continuous and of bounded
variation.

The intuition is that $X_{\tau _{t}\left( .\right) }$ is the branch of a
tree corresponding to the time $t$. Consider two generic times $s<t$, then $%
X_{\tau _{s}\left( .\right) }$ and $X_{\tau _{t}\left( .\right) }\ $agree on
the initial segment \ $\left[ 0,h\left( r\right) \right) $ but thereafter $%
\tau _{s}\left( .\right) \in \left[ r,s\right] $ while $\tau _{t}\left(
.\right) \in \left[ \sup I,t\right] $. The restriction of $X_{\tau
_{t}\left( .\right) }\ $to the initial segment \ $\left[ 0,h\left( r\right)
\right) $ is the path $X_{\tau _{\sup I}\left( .\right) }$. As $h\left(
r\right) =\inf \left[ h\left( u|u\in \left[ s,t\right] \right) \right]$ they
have independent trajectories after $h\left( r\right)$.

Let $\tilde{h}\left( t\right) $ be the total $1$-variation of the path $%
X_{\tau _{t}\left( .\right) }$. The claim is that $\tilde{h}$ has total $1$%
-variation bounded by that of $X$ and is also a height function for $X$.

As the paths $X_{\tau _{s}\left( .\right) }$ and $X_{\tau _{t}\left(
.\right) }$ share the common segment $X_{\tau _{r}\left( .\right) }$ we have 
\begin{equation*}
\left\Vert X_{s}-X_{t}\right\Vert =\left\Vert X_{\tau _{s}\left( s\right)
}-X_{\tau _{t}\left( t\right) }\right\Vert \leq \tilde{h}\left( t\right) -%
\tilde{h}\left( r\right) +\tilde{h}\left( s\right) -\tilde{h}\left( r\right)
,
\end{equation*}
and in particular 
\begin{equation*}
\left\Vert X_{s}-X_{t}\right\Vert \leq \tilde{h}\left( s\right) +\tilde{h}%
\left( t\right) -2\tilde{h}\left( r\right) .
\end{equation*}
On the other hand $\tilde{h}\left( r\right) =\tilde{h}\left( \sup I\right)
=\inf_{s\leq u\leq t}\left( \tilde{h}\left( u\right) \right) $ and so 
\begin{equation*}
\left\Vert X\left( s\right) -X\left( t\right) \right\Vert \leq \tilde{h}%
\left( s\right) +\tilde{h}\left( t\right) -2\inf_{s\leq u\leq t}\left( 
\tilde{h}\left( u\right) \right) .
\end{equation*}
and $\tilde{h}$ is a height function for $X$.

Finally we control the total variation of $\tilde{h}$ by $\omega_X$, the
total variation of the path. In fact, 
\begin{eqnarray*}
\left\vert \tilde{h}\left( s\right) -\tilde{h}\left( t\right) \right\vert
&\leq &\tilde{h}\left( s\right) +\tilde{h}\left( t\right) -2\inf_{s\leq
u\leq t}\left( \tilde{h}\left( u\right) \right) \\
&\leq &\omega _{X}\left( s,t\right),
\end{eqnarray*}
where $\omega_X(s,t) = \sup_{D\in \mathcal{D}} \sum_{D} \left\Vert
X_{t_{i+1}}-X_{t_i}\right\Vert|$, with $\mathcal{D}$ denoting the set of all
partitions of $[s,t]$ and for $D\in \mathcal{D}$, then $D= \{s\leq
\dots<t_i<t_{i+1}<\dots \leq t\}$. The first of these inequalities is
trivial, but the second needs explanation. As before, notice that the paths $%
X_{\tau _{s}\left( .\right) }$, and $X_{\tau _{t}\left( .\right) }$ share
the common segment $X_{\tau _{r}\left( .\right) }$ and that $\inf_{s\leq
u\leq t}\left( \tilde{h}\left( u\right) \right) =\tilde{h}\left( r\right) $.
So $\tilde{h}\left( s\right) +\tilde{h}\left( t\right) -2\inf_{s\leq u\leq
t}\left( \tilde{h}\left( u\right) \right) $ is the total length of the two
segments $X_{\tau _{s}\left( .\right) }|_{\left[ h\left( r\right) ,h\left(
s\right) \right] }$ and $X_{\tau _{t}\left( .\right) }|_{\left[ h\left(
r\right) ,h\left( t\right) \right] }$. Now the total variation of $X_{\tau
_{t}\left( .\right) }|_{\left[ h\left( r\right) ,h\left( t\right) \right] }$
is obviously bounded by $\omega _{X}\left( \sup I,t\right) $, as the path $%
X_{\tau _{t}\left( .\right) }|_{\left[ h\left( r\right) ,h\left( t\right) %
\right] }$ is a time change of $X|_{\left[ \sup I,t\right] }$.

It is enough to show that the total length of $X_{\tau _{s}\left( .\right)
}|_{\left[ h\left( r\right) ,h\left( s\right) \right] }$ is controlled by $%
\omega _{X}\left( s,\inf I\right) $ to conclude that 
\begin{equation*}
\tilde{h}\left( s\right) +\tilde{h}\left( t\right) -2\inf_{s\leq u\leq
t}\left( \tilde{h}\left( u\right) \right) \leq \omega _{X}\left( s,t\right) .
\end{equation*}
In order to do this we work backwards in time. Let 
\begin{eqnarray*}
f_{s}\left( u\right) &=&\inf_{s\leq v\leq u}h\left( v\right) \\
\rho _{t}\left( x\right) &=&\inf \left\{ u\in \left[ s,T\right] |f_{s}\left(
u\right) =x\right\}
\end{eqnarray*}
then, because $X$ is tree-like 
\begin{equation*}
X_{\tau _{s}\left( .\right) }|_{\left[ 0,h\left( s\right) \right] }=X_{\rho
_{s}\left( .\right) }|_{\left[ 0,h\left( s\right) \right] },
\end{equation*}
and in particular, the path segment $X_{\rho _{s}\left( .\right) }|_{\left[
h\left( r\right) ,h\left( s\right) \right] }$ is a time change (but
backwards) of $X|_{\left[ s,\inf I\right] }$.
\end{proof}

The property of being tree-like is re-parameterisation invariant. We see
informally that a tree-like path $X$ is the composition of a contraction on
the $R$-tree defined by $h$ and the based loop in this tree obtained by
taking $t\in \left[ 0,T\right] $ to its equivalence class under the metric
induced by $h$ (for definitions and a proof see the Appendix).

Any path that can be factored through a based loop of finite length in an $R$
-tree and a contraction of that tree to the space $E$ is a Lipschitz
tree-like path. If $0$ is the root of the tree and $\phi $ is the based loop
defined on $\left[ 0,T\right] $, then define $h\left( t\right) =d\left(
0,\phi \left( t\right) \right) $. This makes $\phi $ a tree-like path. Any
Lipschitz image of a tree-like path is obviously a Lipschitz tree-like path.

We have the following trivial lemma.

\begin{lemma}
\label{lem:h}A Lipschitz tree-like path $X$ always has bounded variation
less than that of any height function $h$ for $X$.
\end{lemma}

\begin{proof}
Let $\mathcal{D}=\left\{ t_{0}<\ldots<t_{n}\right\} $ be a partition of $%
\left[ 0,T\right] $. Choose $u_{i}\in\left[ t_{i-1},t_{i}\right] $
maximising $h\left( t_{i}\right) +h\left( t_{i-1}\right) -2h\left(
u_{i}\right) $ and let $\mathcal{\tilde{D}}=\left\{ t_{0}\leq
u_{1}\ldots\leq t_{n-1}\leq u_{n}\leq t_{n}\right\} .$ Relabel the points of 
$\mathcal{\tilde{D}}=\left\{ v_{0}\leq v_{1}\ldots\leq v_{m}\right\} $ Then 
\begin{equation}
\sum_{\mathcal{D}}\left\| X_{t_{i}}-X_{t_{i-1}}\right\| \leq\sum _{\mathcal{%
\tilde{D}}}\left| h\left( v_{i}\right) -h\left( v_{i-1}\right) \right|.
\label{eq:001}
\end{equation}
\end{proof}

We now prove a compactness result.

\begin{lemma}
\label{lem:cmpctness} Suppose that $\{h_{n}\}$ are a sequence of height
functions on $\left[ 0,T\right] $ for a sequence of tree-like paths $%
\{X_{n}\}$. Suppose further that the $h_{n}$ are parameterised at speeds of
at most one and that the $X_{n}$ take their values in a common compact set
within $E$. Then we may find a subsequence $\left( X_{n\left( k\right)
},h_{n\left( k\right) }\right) $ converging uniformly to a Lipschitz
tree-like path $\left( Y,h\right) $. The speed of traversing $h$ is at most
one.
\end{lemma}

\begin{proof}
The $h_{n}$ are equi-continuous, and in view of (\ref{eq:001}) the $X_{n}$
are as well. Our hypotheses are sufficient for us to apply the Arzela-Ascoli
theorem to obtain a subsequence $\left( X_{n\left( k\right) },h_{n\left(
k\right) }\right) $ converging uniformly to some $\left( Y,h\right) $. In
view of the fact that the Lip norm is lower semi-continuous in the uniform
topology, we see that $h$ is a bounded variation function parameterised at
speed at most one and that $Y$ is of bounded variation; of course $h$ takes
the value 0 at both ends of the interval $\left[ 0,T\right] $.

Now $h_{n\left( k\right) }$ converge uniformly to $h$ and hence $\inf _{u\in%
\left[ s,t\right] }h_{n\left( k\right) }\left( u\right) \rightarrow\inf_{u\in%
\left[ s,t\right] }h\left( u\right) ;$ meanwhile the $h_{n}$ are height
functions for the tree-like paths $X_n$ and hence we can take limits through
the definition to show that $h$ is a height function for $Y$.


\end{proof}

\begin{corollary}
Every Lipschitz tree-like path $X$ has a height function $h$ of minimal
total variation and its total variation measure is boundedly absolutely
continuous with respect to the total variation measure of any other height
function.
\end{corollary}

\begin{proof}
We see that this is an immediate corollary of Proposition \ref{prop:h} and
Lemma \ref{lem:h} and \ref{lem:cmpctness}.
\end{proof}

There can be more than one minimiser $h$ for a given $X$.

\section{Approximation of the path}

\label{Section5}

\subsection{Representing the path as a line integral against a rank one form}

Let $\gamma $ be a path of finite variation in a finite dimensional
Euclidean space $V$ with total length $T$ and parameterised at unit speed.
Its parameter set is $[0,T]$. We note that the signature of $\gamma $ is
unaffected by this choice of parameterisation.

\begin{definition}
Let $\gamma \left( \left[ 0,T\right] \right) $ denote the range of $\gamma $
in $V$ and let the occupation measure $\mu$ on $(V,\mathcal{\ B}(V)) $ be
denoted 
\begin{equation*}
\mu\left( A\right) =\left| \left\{ s<T|\gamma \left( s\right) \in A\right\}
\right|,\;\; A\subset V.
\end{equation*}
\end{definition}

Let $n\left( x\right) $ be the number of points on $\left[ 0,T\right] $
corresponding under $\gamma $ to $x\in E$. By the area formulae \cite%
{ohtsuka} p125-126, one has the total variation, or length, of the path $%
\gamma $ is given by 
\begin{equation}
Var\left( \gamma \right) =\int n\left( x\right) \Lambda _{1}\left( dx\right),
\label{eq:areaformula}
\end{equation}
where $\Lambda _{1}$ is one dimensional Hausdorff measure. Moreover, for any
continuous function $f$ 
\begin{equation*}
\int f\left( \gamma \left( t\right) \right) dt=\int f\left( x\right) n\left(
x\right) \Lambda _{1}\left( dx\right).
\end{equation*}
Note that $\mu=n\left( x\right) \Lambda _{1}$ and that $n$ is integrable.

\begin{lemma}
\label{pullback_null} The image under $\gamma $ of a Lebesgue null set is
null for $\mu$. That is to say $\mu(\gamma(N)) = \left| \gamma ^{-1}\gamma
\left( N\right) \right| =0$ if $\left| N\right| =0$.
\end{lemma}

\begin{definition}
We will say that $N\subset \left[ 0,T\right] $ is $\gamma $-stable if $%
\gamma ^{-1}\gamma \left( N\right) =N$.
\end{definition}

As a result of Lemma~\ref{pullback_null} we see that any null set can always
be enlarged to a $\gamma $-stable null set.

The Lebesgue differentiation theorem tells us that $\gamma $ is
differentiable at almost every $u$ in the classical sense, and with this
parameterisation the derivative will be absolutely continuous and of modulus
one.

\begin{corollary}
There is a set $\ G$ of full $\mu$ measure in $V$ so that $\gamma $ is
differentiable with $\left| \gamma ^{\prime }\left( t\right) \right| =1$
whenever $\gamma \left( t\right) \in G$. We set \ $M=\gamma ^{-1}G$. $M$ is $%
\gamma $-stable.
\end{corollary}

Now it may well happen that the path visits the same point $x\in G$ more
than once. A priori, there is no reason why the direct ions of the
derivative on $\left\{ t\in M|\gamma \left( t\right) =m\right\} $ should not
vary. However this can only occur at a countable number of points.

\begin{lemma}
The set of pairs $\left( s,t\right) $ of distinct times in $\ M\times M$ for
which 
\begin{eqnarray*}
\gamma \left( s\right) &=&\gamma \left( t\right) \\
\gamma ^{\prime }\left( s\right) &\neq & \pm \gamma ^{\prime }\left( t\right)
\end{eqnarray*}
is countable.
\end{lemma}

\begin{proof}
If $\gamma(s_*)=\gamma(t_*)$ but $\gamma^{\prime}(s_*) \neq \pm
\gamma^{\prime}(t_*)$ then, by a routine transversality argument, there is
an open neighbourhood of $(s_*,t_*)$ in which there are no solutions of $%
\gamma(s)=\gamma(t)$ except $s=s_*, t=t_*$.
\end{proof}

Up to sign and with countably many exceptions, the derivative of $\gamma $
does not depend on the occasion of the visit to a point, only the location.
Sometimes we will only be concerned with the unsigned or projective
direction of $\gamma $ and identify $v\in S$ with $-v$.

\begin{definition}
For clarity we introduce $\symbol{126}_{\pm }$ as the equivalence relation
that identifies $v$ and $-v$ and let $[\gamma^{\prime}]_{\symbol{126}_{\pm}
}\in S/\symbol{126}_{\pm }$ denote the unsigned direction of $\gamma$.
\end{definition}

$\gamma ^{\prime \pm }$ is defined on the full measure subset of $\left[ 0,T %
\right] $ where $\gamma ^{\prime }$ is defined and in $S$.

\begin{corollary}
There is a function $\phi $ defined on $G$ with values in the projective
sphere $S/\symbol{126}_{\pm }$ so that $\phi \left( \gamma \left( t\right)
\right) =\left[ \gamma ^{\prime }\left( t\right) \right] _{\symbol{126}_{\pm
}}$.
\end{corollary}

As a result we may define a useful vector valued 1-form $\mu $-almost
everywhere on $G$. If $\xi $ is a vector in $S$, then $\left\langle \xi
,u\right\rangle \xi $ is the linear projection of $u$ onto the subspace
spanned by $\xi .$ As $\left\langle \xi ,u\right\rangle \xi =\left\langle
-\xi ,u\right\rangle \left( -\xi \right) $ it defines a function from $S/%
\symbol{126}_{\pm }$ to $Hom\left( V,V\right) $.

\begin{definition}
We define the tangential projection $1$-form $\omega $. Let $\xi $ be a unit
strength vector field on $G$ with $\left[ \xi \right] _{\symbol{126}_{\pm
}}=\phi $. Then 
\begin{equation*}
\omega \left( g,u\right) =\left\langle \xi \left( g\right) ,u\right\rangle
\xi \left( g\right), \;\;\forall g\in G, \forall u
\end{equation*}
defines a vector 1-form. The 1-form depends on $\phi $, but is otherwise
independent of the choice of $\xi $.
\end{definition}

The 1-form $\omega $ is the projection of $u$ onto the line determined by $%
\phi \left( g\right) .$

\begin{theorem}
\begin{proposition}
\label{lem:oneformrep} The tangential projection $\omega $, defined $\mu $
a.e. on $G$, is a linear map from $V\rightarrow V$ with rank one. For almost
every $t$ one has 
\begin{equation*}
\gamma ^{\prime }\left( t\right) =\omega \left( \gamma \left( t\right)
,\gamma ^{\prime }\left( t\right) \right) 
\end{equation*}%
and as a result, using the fundamental theorem of calculus for Lipschitz
functions, 
\begin{eqnarray*}
\gamma \left( t\right)  &=&\int_{0<u<t}d\gamma _{u}+\gamma \left( 0\right) 
\\
&=&\int_{0<u<t}\omega \circ d\gamma _{u}+\gamma \left( 0\right) ,
\end{eqnarray*}%
for every $t\leq T$.
\end{proposition}
\end{theorem}

By approximating $\omega $ by other rank one 1-forms we will be able to
approximate $\gamma $ by (weakly) piecewise linear paths that also have
trivial signature. It will be easy to see that such paths are tree-like. The
set of tree-like paths is closed. This will complete the argument.

\subsection{Iterated integrals of iterated integrals}



We now prove that if $\gamma $ has a trivial signature $\left( 1,0,0,\ldots
\right) $, then it can always be approximated arbitrarily well by weakly
piecewise linear paths with shorter length and trivial signature. Our
approximations will all be line integrals of 1-forms against our basic path $%
\gamma $. Two key points we will need are that the integrals are continuous
against varying the 1-form, and that a line integral of a path with trivial
signature also has trivial signature. The Stone-Weierstrass theorem will
allow us to reduce this second problem to one concerning line integrals
against polynomial 1-forms, and in turn this will reduce to the study of
certain iterated integrals. The application of the Stone-Weierstrass theorem
requires a commutative algebra structure and this is provided by the
coordinate iterated integrals and the shuffle product. For completeness we
set this out below.

Suppose that we define 
\begin{equation*}
Z_{u}:=\underset{0<u_{1}<\ldots <u_{r}<u}{\int \cdots \int }d\gamma
_{u_{1}}\ldots d\gamma _{u_{r}}\in V^{\otimes r}
\end{equation*}
and 
\begin{equation*}
\tilde{Z}_{u}:=\underset{0<u_{1}<\ldots <u_{\tilde{r}}<u}{\int \cdots \int }%
d\gamma _{u_{1}}\ldots d\gamma _{u_{\tilde{r}}}\in V^{\otimes \tilde{r}},
\end{equation*}
then it is interesting as a general point, and necessary here, to consider
iterated integrals of $Z$ and $\tilde{Z}$%
\begin{equation*}
\underset{0<u_{1}<u_{2}<T}{\int \cdots \int }d\tilde{Z}_{u_{1}}dZ_{u_{2}}\in
V^{\otimes \tilde{r}}\otimes V^{\otimes r}.
\end{equation*}
It will be technically important to us to observe that such integrals can
also be expressed as linear combinations of iterated integrals of $\gamma $
so we do this with some care. Some of the results stated below follow from
the well known shuffle product and its relationship with multiplication of
coordinate iterated integrals. 

\begin{definition}
The truncated or $n$-signature $\mathbf{X}_{s,t}^{\left( n\right) }\mathbf{=}%
\left( 1,X_{s,t}^{1},X_{s,t}^{2},\ldots ,X_{s,t}^{n}\right) $ is the
projection of the signature $\mathbf{X}_{s,t}$ to the algebra $T^{\left(
n\right) }\left( V\right) :=\bigoplus_{r=0}^{n}V^{\otimes r}$ of tensors
with degree at most $n$.
\end{definition}

\begin{definition}
\label{defn:pre-coorditint}If $e$ is an element of the dual space $V^{\ast }$
to $V$, then $\gamma _{u}^{e}=\left\langle e,\gamma _{u}\right\rangle $ is a
scalar path and $d\gamma _{u}^{e}=\left\langle e,d\gamma _{u}\right\rangle $%
. If $\mathbf{e}=\left( e_{1},\ldots ,e_{r}\right) $ is a list of elements
of the dual space to $V$, then we define the coordinate iterated integral 
\begin{equation*}
X_{s,t}^{\mathbf{e}}:=\underset{s<u_{1}<\ldots <u_{r}<t}{\int \cdots \int }%
d\gamma _{u_{1}}^{e_{1}}\ldots d\gamma _{u_{r}}^{e_{r}}=\left\langle \mathbf{%
e},X_{s,t}^{r}\right\rangle .
\end{equation*}
\end{definition}

\begin{lemma}
The map $\mathbf{e}\rightarrow X_{s,t}^{\mathbf{e}}$ defined above extends
uniquely as a linear map from $T^{\left( n\right) }\left( V^{\ast }\right) $
to the space of real valued functions on paths of bounded variation.
\end{lemma}

\begin{proof}
Let $\mathbf{e}\in T^{\left( n\right) }\left( V^{\ast }\right) $. Since $%
T^{\left( n\right) }\left( V^{\ast }\right) $ is dual to $T^{\left( n\right)
}\left( V\right) $ the pairing $\mathbf{e}\rightarrow \left\langle \mathbf{%
e,X}_{s,t}^{\left( n\right) }\right\rangle $ defines a real number of each
path $\gamma .$If $\mathbf{e}=\left( e_{1},\ldots ,e_{r}\right) $ then this
coincides with $X_{s,t}^{\mathbf{e}}$; since such vectors span $T^{\left(
n\right) }\left( V^{\ast }\right) $ the result is immediate.
\end{proof}

We therefore extend definition \ref{defn:pre-coorditint}.

\begin{definition}
For any $n,$ $\mathbf{e}\in T^{\left( n\right) }\left( V^{\ast }\right) $ we
call $X_{s,t}^{\mathbf{e}}$ the $\mathbf{e}$-coordinate iterated integral of 
$\gamma $ over the interval $\left[ s,t\right] $.
\end{definition}



These functions on path space are important because they form an algebra
under pointwise multiplication and because they are like polynomials and so
it is easy to define a differentiation operator on this space. Given two
tensors $\mathbf{e}$, $\mathbf{f}$ there is a natural product $\mathbf{%
e\sqcup f}$, called the shuffle product, derived from the above. For basic
tensors 
\begin{eqnarray*}
\mathbf{e} &=&e_{1}\otimes \ldots \otimes e_{r}\in V^{\otimes r} \\
\mathbf{f} &=&f_{1}\otimes \ldots \otimes f_{s}\in V^{\otimes s}
\end{eqnarray*}%
and a shuffle $\left( \pi _{1},\pi _{2}\right) $ (a pair of increasing
injective functions from $(1,..,r),(1,..,s)$ to $(1,..,r+s)$ with disjoint
range) one can define a tensor of degree $r+s$: 
\begin{equation*}
\omega _{\left( \pi _{1},\pi _{2}\right) }=\omega _{1}\otimes \ldots \otimes
\omega _{r+s},
\end{equation*}%
where $\omega _{\pi _{1}\left( j\right) }=e_{j}$ for $j=1,\ldots r$ and $%
\omega _{\pi _{2}\left( j\right) }=f_{j}$ for $j=1,\ldots s$. Since the
ranges of $\pi _{1}$ and $\pi _{2}$ are disjoint a counting argument shows
that the union of the ranges is $1,\ldots ,r+s,$ and that $\omega _{k}$ is
well defined for all $k$ in $1,\ldots r+s$ and hence $\omega _{\left( \pi
_{1},\pi _{2}\right) }$ is defined. By summing over all shuffles 
\begin{equation*}
\mathbf{e\sqcup f}\mathbf{=}\sum_{\left( \pi _{1},\pi _{2}\right) }\omega
_{\left( \pi _{1},\pi _{2}\right) }
\end{equation*}%
one defines a multilinear map of $V^{\otimes r}\mathbf{\times }V^{\otimes s}%
\mathbf{\rightarrow }V^{\otimes \left( r+s\right) }$.

\begin{definition}
The unique extension of $\mathbf{\sqcup }$ to a map from $T\left( V\right) 
\mathbf{\times }T\left( V\right) \rightarrow T\left( V\right) $ is called
the shuffle product.
\end{definition}

The following is standard.

\begin{lemma}
The class of coordinate iterated integrals is closed under pointwise
multiplication. For each $\gamma $ the point-wise product of the $\mathbf{e}$%
-coordinate iterated integral and the $\mathbf{f}$-coordinate iterated
integral is the $\left( \mathbf{e\sqcup f}\right) $-coordinate iterated
integral: 
\begin{equation*}
X_{s,t}^{\mathbf{e}}X_{s,t}^{\mathbf{f}}=X_{s,t}^{\mathbf{e\sqcup f}}.
\end{equation*}%
%
%
%
%
%
%
%
%
%
%
%
%
\end{lemma}


\begin{corollary}
\label{cor:poly} Any polynomial in coordinate iterated integrals is a
coordinate iterated integral.
\end{corollary}

\begin{remark}
\bigskip It is at first sight surprising that any polynomial in the linear
functionals on T$\left( V\right) $ coincides with a unique linear functional
on $T\left( V\right) $ when restricted to signatures of paths and reflects
the fact that the signature of a path is far from being a generic element of
the tensor algebra.
\end{remark}

A slightly more demanding remark relates to iterated integrals of coordinate
iterated integrals.

\begin{proposition}
\label{prop:poly} The iterated integral 
\begin{equation}
\underset{s<u_{1}<\ldots <u_{r}<t}{\int \cdots \int }dX_{s,u_{1}}^{\mathbf{e}%
_{1}}\ldots dX_{s,u_{r}}^{\mathbf{e}_{r}}  \label{eq:iter_iter}
\end{equation}%
is itself a coordinate iterated integral.
\end{proposition}

\begin{proof}
A simple induction ensures that it suffices to consider the case 
\begin{equation*}
\underset{s<u_{1}<u_{r}<t}{\int \int }dX_{s,u_{1}}^{\mathbf{e}}dX_{s,u_{2}}^{%
\mathbf{f}},
\end{equation*}%
where 
\begin{eqnarray*}
\mathbf{e} &=&e_{1}\otimes \ldots \otimes e_{r}\in (V^{\ast })^{\otimes r} \\
\mathbf{f} &=&f_{1}\otimes \ldots \otimes f_{s}\in (V^{\ast })^{\otimes s}
\end{eqnarray*}%
and in this case 
\begin{equation*}
\underset{s<u_{1}<u_{2}<t}{\int \int }dX_{s,u_{1}}^{\mathbf{e}}dX_{s,u_{2}}^{%
\mathbf{f}}=\underset{%
\begin{array}{c}
s<v_{1}<\ldots <v_{r}<t \\ 
s<w_{1}<\ldots <w_{s}<t \\ 
v_{r}<w_{s}%
\end{array}%
}{\int \cdots \int }d\gamma _{v_{1}}^{e_{1}}\ldots d\gamma
_{v_{r}}^{e_{r}}d\gamma _{w_{1}}^{f_{1}}\ldots d\gamma _{w_{s}}^{f_{s}}.
\end{equation*}%
Expressing the integral as a sum of integrals over the regions where the
relative ordering of the $v_{i}$ and $w_{j}$ are preserved (i.e. all
shuffles for which the last card comes from the right hand pack) we have 
\begin{eqnarray*}
\underset{s<u_{1}<u_{2}<t}{\int \int }dX_{s,u_{1}}^{\mathbf{e}}dX_{s,u_{2}}^{%
\mathbf{f}} &=&X_{s,t}^{\left( \mathbf{e\sqcup \tilde{f}}\right) \otimes
f_{s}} \\
\mathbf{\tilde{f}} &\mathbf{=}&f_{1}\otimes \ldots \otimes f_{s-1}.
\end{eqnarray*}
\end{proof}

From this it is, of course, clear that

\begin{lemma}
\label{lem:itit} If a path has trivial signature, then all iterated
integrals of its iterated integrals are zero.
\end{lemma}

\subsection{Bounded, measurable, and integrable forms}

Recall that $\gamma $ is a path of finite length in $V$, and that it is
parameterized at unit speed. The occupation measure is $\mu $ and has total
mass equal to the length $T$ of the path $\gamma $. Let $\left( W,\left\Vert
{}\right\Vert \right) $ be a normed space with a countable base (usually $V$
itself). If $\omega $ is a $\mu $-integrable $1$-form with values in $W$
then we write 
\begin{equation*}
\left\Vert \omega \right\Vert _{L^{1}\left( V,\mathcal{B}\left( V\right)
\right) }=\int_{V}\left\Vert \omega \left( y\right) \right\Vert _{Hom\left(
V,W\right) }\mu \left( dy\right) =\int_{0}^{T}\left\Vert \omega \left(
\gamma _{t}\right) \right\Vert _{Hom\left( V,W\right) }dt.
\end{equation*}

\begin{proposition}
Let $\omega \in L^{1}\left( V,\mathcal{B}\left( V\right) ,\mu \right) $ be a 
$\mu $-integrable $1$-form with values in $W$. Then the indefinite line
integral $y_{t}:=\int_{0}^{t}\omega \left( d\gamma _{t}\right) $ is well
defined, linear in $\omega ,$ and a path in $W$ with $1$-variation at most $%
\left\Vert \omega \right\Vert _{L^{1}\left( V,\mathcal{B}\left( V\right)
,\mu \right) }$.
\end{proposition}

\begin{proof}
Since $\omega $ is a $1$-form defined $\mu $-almost surely, $\omega \left(
\gamma _{t}\right) \in Hom\left( V,W\right) $ (where $Hom\left( V,W\right) $
is equipped with the operator norm $\left\Vert \cdot \right\Vert $) is
defined $dt$ almost everywhere. Since $\omega $ is integrable, it is
measurable, and hence $\omega \left( \gamma _{t}\right) $ is measurable on $%
\left[ 0,T\right] $. Since $\gamma $ has finite variation and is
parameterized at unit speed, it is differentiable almost everywhere and its
derivative is measurable with unit length $dt$ almost surely. Hence $\omega
\left( \gamma _{t}\right) \left( \dot{\gamma}_{t}\right) $ is measurable and
dominated by $\left\Vert \omega \left( \gamma _{t}\right) \right\Vert $,
which is an integrable function, and hence $\omega \left( \gamma _{t}\right)
\left( \dot{\gamma}_{t}\right) $ is integrable. Thus the line integral can
be defined to be 
\begin{eqnarray*}
y_{t} &=&\int_{0}^{t}\omega \left( \gamma _{u}\right) \left( \dot{\gamma}%
_{u}\right) du \\
\left\Vert y_{t}-y_{s}\right\Vert &\leq &\int_{s}^{t}\left\Vert \omega
\left( \gamma _{u}\right) \right\Vert _{Hom\left( V,W\right) }\left\Vert 
\dot{\gamma}_{u}\right\Vert du \\
&=&\int_{s}^{t}\left\Vert \omega \left( \gamma _{u}\right) \right\Vert
_{Hom\left( V,W\right) }du
\end{eqnarray*}%
and so has $1$-variation bounded by $\left\Vert \omega \right\Vert
_{L^{1}\left( V,\mathcal{B}\left( V\right) \right) }$.
\end{proof}

\begin{proposition}
\label{prop:conv} Let $\omega _{n}\in L^{1}\left( V,\mathcal{B}\left(
V\right) ,\mu \right) $ be a uniformly bounded sequence of integrable $1$%
-forms with values in a vector space $W$. Suppose that they converge in $%
L^{1}\left( V,\mathcal{B}\left( V\right) ,\mu \right) $ to $\omega $, then
the signatures of the line integrals $\int \omega _{n}\left( d\gamma
_{t}\right) $ converge to the signature of $\int \omega \left( d\gamma
_{t}\right) $.
\end{proposition}

\begin{proof}
The $r$'th term in the iterated integral of the line integral $\int \omega
_{n}\left( d\gamma _{t}\right) $ can be expressed as 
\begin{equation*}
\underset{0<u_{1}<\ldots <u_{r}<T}{\int \cdots \int }\omega _{n}\left(
\gamma _{u_{1}}\right) \otimes \ldots \otimes \omega _{n}\left( \gamma
_{u_{r}}\right) \;\left( \dot{\gamma}_{u_{1}}\right) \ldots \left( \dot{%
\gamma}_{u_{r}}\right) \;du_{1}\ldots du_{r}
\end{equation*}
and since $\omega _{n}$ converge in $L^{1}\left( V,\mathcal{B}\left(
V\right) ,\mu \right) $, it follows that from the definition of $\mu $ that $%
\omega _{n}\left( \gamma _{u}\right) $ converge to $\omega \left( \gamma
_{u}\right) $ in $L^{1}\left( \left[ 0,T\right] ,\mathcal{B}\left( \mathbb{R}%
\right) ,du\right) $ almost everywhere. Thus $\omega _{n}\left( \gamma
_{u_{1}}\right) \otimes \ldots \otimes \omega _{n}\left( \gamma
_{u_{r}}\right) $ converges in $L^{1}\left( \left[ 0,T\right] ^{r},\mathcal{B%
}\left( \mathbb{R}\right) ,du_{1}\ldots du_{r}\right) $. Since $\left\| \dot{%
\gamma}_{u}\right\| =1$ for almost every $u$, Fubini's theorem implies that 
\begin{equation*}
\omega _{n}\left( \gamma _{u_{1}}\right) \otimes \ldots \otimes \omega
_{n}\left( \gamma _{u_{r}}\right) \;\left( \dot{\gamma}_{u_{1}}\right)
\ldots \left( \dot{\gamma}_{u_{r}}\right)
\end{equation*}
converges in $L^{1}\left( \left[ 0,T\right] ^{r},\mathcal{B}\left( \mathbb{R}%
\right) ,du_{1}\ldots du_{r}\right) $ to 
\begin{equation*}
\omega \left( \gamma _{u_{1}}\right) \otimes \ldots \otimes \omega \left(
\gamma _{u_{r}}\right) \;\left( \dot{\gamma}_{u_{1}}\right) \ldots \left( 
\dot{\gamma}_{u_{r}}\right).
\end{equation*}
Thus, integrating over $0<u_{1}<\ldots <u_{r}<T$, the proposition follows.
\end{proof}

\begin{corollary}
\label{cor:1formtriv} Let $\omega \in L^{1}\left( V,\mathcal{B}\left(
V\right) ,\mu \right) $. If $\gamma $ has trivial signature, then so does $%
\int \omega \left( d\gamma _{t}\right) $. That is to say, for each $r$, 
\begin{equation*}
\underset{0<u_{1},\ldots ,u_{r}<T}{\int \cdots \int }\omega \left( d\gamma
_{u_{1}}\right) \ldots \omega \left( d\gamma _{u_{r}}\right) =0\in
W^{\otimes r}.
\end{equation*}
\end{corollary}

\begin{proof}
It is a consequence of Proposition~\ref{prop:conv} that the set of $%
L^{1}\left( V,\mathcal{B}\left( V\right) ,\mu \right) $ forms producing line
integrals having trivial signature is closed. By Lusin's theorem, one may
approximate, in the $L^{1}\left( V,\mathcal{B}\left( V\right) ,\mu \right) $
norm, any integrable form by bounded continuous forms. If the initial form
is uniformly bounded then the approximations can be chosen to satisfy the
same uniform bound.

The support of $\mu $ is compact, so by the Stone Weierstrass theorem, we
can uniformly approximate these continuous forms by polynomial forms $\omega
=\sum_{i}p_{i}e_{i}$, where the $p_{i}$ are polynomials and $e_{i}$ are a
basis for $V^{\ast }$. Using the fact that 
\begin{equation*}
\frac{(\gamma _{T}^{e_{1}})^{r}}{r!}=\underset{0<u_{1},\ldots ,u_{r}<T}{\int
\cdots \int }d\gamma _{u_{1}}^{e_{1}}\dots d\gamma _{u_{r}}^{e_{1}},
\end{equation*}%
with Corollary~\ref{cor:poly} and Proposition~\ref{prop:poly}, we have that
the line integrals against these polynomial forms and their iterated
integrals can be expressed as linear combinations of coordinate iterated
integrals. If $\gamma $ has trivial signature, then by Lemma~\ref{lem:itit}
these integrals will all be zero. It follows from the $L^{1}\left( V,%
\mathcal{B}\left( V\right) ,\mu \right) $ continuity of the truncated
signature, that the signature of the path formed by taking the line integral
against any form $\omega $ in $L^{1}\left( V,\mathcal{B}\left( V\right) ,\mu
\right) $ will always be trivial.
\end{proof}

\subsection{Approximating rank one $1$-forms}

\begin{definition}
A vector valued 1-form $\omega $ is (at each point of $V$) a linear map
between vector spaces. We say the 1-form $\omega $ is of rank $k\in \mathbb{N%
}$ on the support of $\mu $ if $\dim \left( \omega \left( V\right) \right)
\leq k$ at $\mu $ almost every point in $V$.
\end{definition}

A linear multiple of a form has the same rank as the original form, but in
general the sum of two forms has any rank less than or equal to the sum of
the ranks of the individual components. However, we will now explain how one
can approximate any rank one $1$-form by piecewise constant rank one 1-forms 
$\omega$. Additionally we will choose the approximations so that, for some $%
\varepsilon >0,$ if $\omega \left( x\right) \neq \omega \left( y\right) $
and $\left| x-y\right| \leq \varepsilon $, then either $\omega \left(
x\right) $ or $\omega \left( y\right) $ is zero.

In other words $\omega $ is rank one and constant on patches which are
separated by thin barrier regions on which it is zero. The patches can be
chosen to be compact and so that the $\mu $-measure of the compliment is
arbitrarily small.

We will use the following easy consequence of Lusin's theorem for one forms
defined on a $\mu $-measurable set $K$. :

\begin{lemma}
\label{lem:lus} Let $\omega $ be a measurable 1-form $\omega $ in $%
L^{1}\left( V,\mathcal{B}\left( V\right) ,\mu \right) $. For each $%
\varepsilon >0$ there is a compact subset $L$ of $\gamma \left[ 0,T\right] $
so that $\omega $ restricted to $L$ is continuous, while $\int_{K\backslash
L}\left\Vert \omega \right\Vert _{Hom\left( V,W\right) }\mu \left( dx\right)
<\varepsilon $.
\end{lemma}

\begin{lemma}
If $\omega $ is a measurable 1-form in $L^{1}\left( V,\mathcal{B}\left(
V\right) ,\mu \right) $, then for each $\varepsilon >0$ there are finitely
many disjoint compact subsets $K_{i}$ of $K$ and a 1-form $\tilde{\omega}$,
that is zero off $\cup _{i}K_{i}$ and constant on each $K_{i}$, such that 
\begin{equation*}
\int_{K}\left\Vert \omega -\tilde{\omega}\right\Vert _{Hom\left( V,W\right)
}\mu \left( dx\right) \leq 4\varepsilon
\end{equation*}%
and with the property that $\tilde{\omega}$ is rank one if $\omega $ is.
\end{lemma}

\begin{proof}
Let $L$ be the compact subset introduced in Lemma~\ref{lem:lus}. Now $\omega
\left( L\right) $ is compact. Fix $\varepsilon >0$ and choose $l_{1},\ldots
,l_{n}$ so that 
\begin{equation*}
\omega \left( L\right) \subset \cup _{i=1}^{n}B\left( \omega \left(
l_{i}\right) ,\frac{\varepsilon }{\mu \left( L\right) }\right)
\end{equation*}%
and put 
\begin{equation*}
F_{j}=\omega ^{-1}\left( \cup _{i=1}^{j}B\left( \omega \left( l_{i}\right) ,%
\frac{\varepsilon }{\mu \left( L\right) }\right) \right) .
\end{equation*}%
Now choose a compact set $K_{j}\subset F_{j}\setminus F_{j-1}$ so that 
\begin{equation*}
\mu \left( \left( F_{j}\setminus F_{j-1}\right) \setminus K_{j}\right) \leq 
\frac{\varepsilon 2^{-j}}{\left\Vert \omega \right\Vert _{L^{\infty }\left(
L,\mathcal{B}\left( L\right) ,\mu \right) }}.
\end{equation*}%
Then the $K_{j}$ are disjoint 
and of diameter $\frac{2\varepsilon }{\mu \left( L\right) }$. Moreover 
\begin{eqnarray*}
L &=&F_{n} \\
\mu \left( L\backslash \cup _{i=1}^{n}K_{j}\right) &\leq &\frac{\varepsilon 
}{\left\Vert \omega \right\Vert _{L^{\infty }\left( L,\mathcal{B}\left(
L\right) ,\mu \right) }}
\end{eqnarray*}%
and 
\begin{equation*}
\int_{L\backslash \cup _{i=1}^{n}K_{j}}\left\Vert \omega \right\Vert
_{Hom\left( V,W\right) }\mu \left( dx\right) <\varepsilon .
\end{equation*}

For each non-empty $K_{j}$ choose $k_{j}\in K_{j}$. Define $\tilde{\omega}$
as follows. 
\begin{eqnarray*}
\tilde{\omega}\left( k\right) &=&\omega \left( k_{j}\right) ,\;\;\;k\in K_{j}
\\
\tilde{\omega}\left( k\right) &=&0,\;\;\;\;\;\;\;k\in K\backslash \cup
_{i=1}^{n}K_{j}.
\end{eqnarray*}%
Then 
\begin{eqnarray*}
\int_{L\backslash \cup _{i=1}^{n}K_{j}}\left\Vert \omega -\tilde{\omega}%
\right\Vert _{Hom\left( V,W\right) }\mu \left( dx\right) &<&\varepsilon \\
\int_{\cup _{i=1}^{n}K_{j}}\left\Vert \omega -\tilde{\omega}\right\Vert
_{Hom\left( V,W\right) }\mu \left( dx\right) &<&\frac{2\varepsilon }{\mu
\left( L\right) }\mu \left( L\right) ,
\end{eqnarray*}%
using Lemma~\ref{lem:lus} one has 
\begin{equation*}
\int_{K\backslash L}\left\Vert \omega -\tilde{\omega}\right\Vert _{Hom\left(
V,W\right) }\mu \left( dx\right) <\varepsilon
\end{equation*}%
and finally 
\begin{equation*}
\int_{K}\left\Vert \omega -\tilde{\omega}\right\Vert _{Hom\left( V,W\right)
}\mu \left( dx\right) \leq 4\varepsilon .
\end{equation*}

If $\omega $ had rank $1$ at almost every point of $K$, then it will have
rank $1$ everywhere on $L$ since $\omega $ is continuous. As either $\tilde{%
\omega}\left( k\right) =\omega \left( k_{j}\right) $ for some $k_{j}$ in $L$
or is zero, the form $\tilde{\omega}$ has rank one also.
\end{proof}

\begin{proposition}
\label{cor:1formapprox} Consider the set $\mathcal{P}$ of one forms on a set 
$K$ for which there exists finitely many disjoint compact subsets $K_{i}$ of 
$K$ so that the 1-form is zero off the $K_{i}$ and constant on each $K_{i}$.
The set of rank one 1-forms in  $\mathcal{P}$ is a dense subset in the $%
L^{1}\left( V,\mathcal{B}\left( V\right) ,\mu \right) $topology of the set
of rank one $1$-forms in $L^{1}\left( V,\mathcal{B}\left( V\right) ,\mu
\right) $. 
\end{proposition}

\section{Piecewise linear paths with no repeated edges}

\label{section6}We call a path $\gamma $ piecewise linear if it is
continuous, and if there is a finite partition 
\begin{equation*}
0=t_{0}<t_{1}<t_{2}<\ldots <t_{r}=T
\end{equation*}%
such that $\gamma $ is linear (or more generally, geodesic) on each segment $%
\left[ t_{i},t_{i+1}\right] .$

\begin{definition}
We say the path is non-degenerate if we can choose the partition so that $%
\left[ \gamma _{t_{i-1}},\gamma _{t_{i}}\right] $ and $\left[ \gamma
_{t_{i}},\gamma _{t_{i+1}}\right] $ are not collinear for any $0<i<r$ and if
the $\left[ \gamma _{t_{i-1}},\gamma _{t_{i}}\right] $ are non-zero for
every $0<i\leq r$.
\end{definition}

The positive length condition is automatic if the path is parameterised at
unit speed and $0<T$. If $\theta _{i}$ is the angle $\measuredangle \gamma
_{t_{i-1}}\gamma _{t_{i}}\gamma _{t_{i+1}}$, then $\gamma $ is
non-degenerate if we can find a partition so that for each $0<i<r$ one has 
\begin{equation*}
\left\vert \theta _{i}\right\vert \neq 0\func{mod}\pi .
\end{equation*}%
This partition is unique, and we refer to the $\left[ \gamma
_{t_{i-1}},\gamma _{t_{i}}\right] $ as the $i$-th linear segment in $\gamma .
$ We see, from the quantitative estimate in part 1 of Lemma \ref{spaced_out}%
, that if we choose $\theta =\frac{1}{2}\min \left\vert \theta
_{i}\right\vert $ and scale it so that the length of the minimal segment is
at least $K\left( \theta \right) =\log \left( \frac{2}{1-\cos \left\vert
\theta \right\vert }\right) $, then its development into hyperbolic space is
non-trivial and so its signature is not zero. That is to say, Lemma \ref%
{spaced_out} contains all the information you need to give a quantitative
form of Chen's uniqueness result in the context of piecewise linear paths:

\begin{theorem}
\label{cor:nondeg} If $\gamma $ is a non-degenerate piecewise linear path, $%
2\theta $ is the smallest angle between adjacent edges, and $D>0$ is the
length of the shortest edge then there is at least one $n$ for which 
\begin{equation*}
\left( \frac{2}{1-\cos \left\vert \theta \right\vert }\right) ^{\left( 1-%
\frac{1}{D}\right) }\leq n!\left\Vert \underset{0<u_{1}<\ldots <u_{n}<T}{%
\int \cdots \int }d\gamma _{u_{1}}\ldots d\gamma _{u_{\tilde{r}}}\right\Vert 
\end{equation*}%
and in particular $\gamma $ has non-trivial signature.
\end{theorem}

\begin{proof}
Choose $\alpha =K\left( \theta \right) /D$. Isometrically embed $V$ into $%
SO\left( I_{d}\right) $ and let $\Gamma _{\alpha }$ be the development of $%
\alpha \gamma $. Then $\Gamma _{\alpha ,t}o$ is a piecewise geodesic path in
hyperbolic space satisfying the hypotheses in Lemma \ref{spaced_out}. Thus
we can deduce that the distance $d\left( o,\Gamma _{\alpha }o\right) $ is at
least $K\left( \theta \right) >0.$ As in the discussion before
Theorem~\ref{thm:strongrecover} in Section 3.4 we have
\begin{eqnarray*}
e^{K\left( \theta \right) } &\leq &\left\Vert \Gamma _{\alpha }\right\Vert 
\\
&\leq &\sum_{n=0}^{\infty }\alpha ^{n}\left\Vert \underset{0<u_{1}<\ldots
<u_{n}<T}{\int \cdots \int }d\gamma _{u_{1}}\ldots d\gamma _{u_{\tilde{r}%
}}\right\Vert  \\
&=&\sum_{n=0}^{\infty }\frac{1}{n!}\left( \frac{K\left( \theta \right) }{D}%
\right) ^{n}n!\left\Vert \underset{0<u_{1}<\ldots <u_{n}<T}{\int \cdots \int 
}d\gamma _{u_{1}}\ldots d\gamma _{u_{\tilde{r}}}\right\Vert .
\end{eqnarray*}%
Now multiplying both side by $e^{-\frac{K\left( \theta \right) }{D}}$
we have
\begin{equation*}
e^{-\frac{K\left( \theta \right) }{D}}e^{K\left( \theta \right) }\leq e^{-%
\frac{K\left( \theta \right) }{D}}\sum_{n=0}^{\infty }\frac{1}{n!}\left( 
\frac{K\left( \theta \right) }{D}\right) ^{n}n!\left\Vert \underset{%
0<u_{1}<\ldots <u_{n}<T}{\int \cdots \int }d\gamma _{u_{1}}\ldots d\gamma
_{u_{\tilde{r}}}\right\Vert .
\end{equation*}%
Since any integrable function has at least one point where its value equals
or exceeds its average and since 
\begin{equation*}
1=e^{-\frac{K\left( \theta \right) }{D}}\sum_{n=0}^{\infty }\frac{1}{n!}%
\left( \frac{K\left( \theta \right) }{D}\right) ^{n}
\end{equation*}%
we can conclude an absolute lower bound on the $L^{1}$ norm of the
signature, against the Poisson measure and conclude that there is an $n$ for
which 
\begin{equation*}
e^{K\left( \theta \right) \left( 1-\frac{1}{D}\right) }\leq n!\left\Vert 
\underset{0<u_{1}<\ldots <u_{n}<T}{\int \cdots \int }d\gamma _{u_{1}}\ldots
d\gamma _{u_{\tilde{r}}}\right\Vert .
\end{equation*}%
Recalling the form of $K\left( \theta \right) $ we have the result.
\end{proof}

\begin{corollary}
\label{cor:chi}Any piecewise linear path $\gamma $ that has trivial
signature is tree-like with a height function $h$ having the same total
variation as $\gamma $.
\end{corollary}

\begin{proof}
We will proceed by induction on the number $r$ of edges in the minimal
partition 
\begin{equation*}
0=t_{0}<t_{1}<t_{2}<\ldots <t_{r}=T
\end{equation*}%
of $\gamma $. We assume that $\gamma $ is linear on each segment $\left[
t_{i},t_{i+1}\right] $ and that $\gamma $ is always parameterised at unit
speed.

We assume that $\gamma $ has trivial signature. Our goal is to find a
continuous real valued function $h$ with $h\geq 0,$ $h\left( 0\right)
=h\left( T\right) =0$, and so that for every $s$, $t\in \left[ 0,T\right] $
one has 
\begin{eqnarray*}
\left\vert h\left( s\right) -h\left( t\right) \right\vert  &\leq &\left\vert
t-s\right\vert  \\
\left\vert \gamma _{s}-\gamma _{t}\right\vert  &\leq &h\left( s\right)
+h\left( t\right) -2\inf_{u\in \left[ s,t\right] }h\left( u\right) .
\end{eqnarray*}

If $r=0$ the result is obvious; in this case $T=0$ and the function $h=0$
does the job.

Now suppose that the minimal partition into linear pieces has $r>0$ pieces.
By Corollary~\ref{cor:nondeg}, it must be a degenerate partition. In other
words one of the $\theta _{i}=\measuredangle \gamma _{t_{i-1}}\gamma
_{t_{i}}\gamma _{t_{i+1}}$ must have 
\begin{equation*}
\left| \theta _{i}\right| =0\func{mod}\pi .
\end{equation*}
If $\theta _{i}=\pi $ the point $t_{i}$ could be dropped from the partition
and the path would still be linear. As we have chosen the partition to be
minimal this case cannot occur and we conclude that $\theta _{i}=0$ and the
path retraces its trajectory for an interval of length 
\begin{equation*}
s=\min \left( \left| t_{i}-t_{i-1}\right| ,\left| t_{i+1}-t_{i}\right|
\right) >0.
\end{equation*}
Now $\gamma \left( t_{i}-u\right) =\gamma \left( t_{i}+u\right) $ for $u\in %
\left[ 0,s\right] $ and either $t_{i}-s=t_{i-1}$ or $t_{i}+s=t_{i+1}$.
Suppose that the former holds. Consider the path segments obtained by
restricting the path to the disjoint intervals 
\begin{eqnarray*}
\gamma _{-} &=&\gamma |_{\left[ 0,t_{i-1}\right] } \\
\gamma _{+} &=&\gamma |_{\left[ t_{i}+s,T\right] } \\
\tau &=&\gamma |_{\left[ t_{i}-s,t_{i}+s\right] },
\end{eqnarray*}
then $\gamma =\gamma _{-}\ast \tau \ast \gamma _{+}$ where $\ast $ denotes
concatenation.

Because the signature map $\gamma \rightarrow S\left( \gamma \right) $ is a
homomorphism one sees that the product of the signatures associated to the
segments is the signature of the concatenation of the paths and hence is
trivial, 
\begin{eqnarray*}
S\left( \gamma _{-}\right) \otimes S\left( \tau \right) \otimes S\left(
\gamma _{+}\right)  &=&S\left( \gamma \right)  \\
&=&1\oplus 0\oplus 0\oplus \ldots \in T\left( V\right) .
\end{eqnarray*}%
On the other hand the path $\tau $ is a linear trajectory followed by its
reverse and as reversal produces the inverse signature 
\begin{equation*}
S\left( \tau \right) =1\oplus 0\oplus 0\oplus \ldots \in T\left( V\right) .
\end{equation*}%
Thus 
\begin{equation*}
S\left( \gamma _{-}\right) \otimes S\left( \gamma _{+}\right) =1\oplus
0\oplus 0\oplus \ldots \in T\left( V\right) 
\end{equation*}%
and so the concatenation of $\gamma _{-}$ and $\gamma _{+}$ ($\gamma $ with $%
\tau $ excised) also has a trivial signature. As it is piecewise linear with
at least one less edge we may apply the induction hypothesis to conclude
that this reduced path is tree-like. Let $\tilde{h}$ be the height function
for the reduced path. Then define 
\begin{eqnarray*}
h\left( u\right)  &=&\tilde{h}\left( u\right) ,u\in \left[ 0,t_{i-1}\right] 
\\
h\left( u\right)  &=&\tilde{h}\left( u-2s\right) ,u\in \left[ t_{i}+s,T%
\right]  \\
h\left( u\right)  &=&s-\left\vert t_{i}-u\right\vert +\tilde{h}\left(
t_{i-1}\right) ,u\in \left[ t_{i}-s,t_{i}+s\right] .
\end{eqnarray*}%
It is easy to check that $h$ is a height function for $\gamma $ with the
required properties.
\end{proof}

The reader should note that the main result of the paper Theorem \ref%
{thm:main} linking the signature to tree-like equivalence relies on Chen's
result only through the above Corollary and hence only requires a version
for piecewise linear paths with no repeated edges. Our quantitative Theorem %
\ref{cor:nondeg} provides an independent proof of this result but, in this
context, is stronger than is necessary; Chen's non-quantitative result could
equally well have been used.

We end this section with two straightforward results which will establish
half of our main theorem.

\begin{lemma}
\label{lem:plapprox}\label{lemma:plapprox} If $\gamma $ is a Lipschitz
tree-like path with height function $h$, then one can find piecewise linear
Lipschitz tree-like paths converging in total variation to a
re-parameterisation of $\gamma $.
\end{lemma}

\begin{proof}
Without loss of generality we may re-parametrise time to be the arc length
of $h$. Since $h$ is of bounded variation, using the area formula (\ref%
{eq:areaformula}), we can find finitely many points $u_{n}$ within $\delta $
of one another and increasing in $\left[ 0,T\right] $ so that $h$ takes the
value $h\left( u_{n}\right) $ only finitely many times and only at the times 
$u_{n}$. Consider the path $\gamma _{n}$ that is linear on the intervals $%
\left( u_{n},u_{n+1}\right) $ and agrees with $\gamma $ at the times $u_{n}$%
. Define $h_{n}$ similarly. Then $h_{n}$ is a height function for $\gamma
_{n}$ and so $\gamma _{n}$ is a tree. The paths $\gamma _{n}$ converge to $%
\gamma $ uniformly, and in $p$-variation for all $p>1$. However, as we have
parameterised $h$ by arc length, it follows that the total variation of $%
\gamma $ is absolutely continuous with respect to arc length. As $\gamma
_{n} $ is a martingale with respect to the filtration determined by the
successive time partitions, applying the martingale convergence theorem, it
follows that $\gamma _{n}$ converges to $\gamma $ in $L_{1}$.
\end{proof}

\begin{corollary}
\label{cor:half} Any Lipschitz tree-like path has all iterated integrals
equal to zero.
\end{corollary}

\begin{proof}
For piecewise linear tree-like paths it is obvious by induction on the
number of segments that all the iterated integrals are 0. Since the process
of taking iterated integrals is continuous in $p$-variation norm for $p<2$,
and Lemma~\ref{lem:plapprox} proves that any Lipschitz tree-like path can be
approximated by piecewise linear tree-like paths in $1$-variation the result
follows.
\end{proof}

In the next section we introduce the concept of a weakly piecewise linear
path. After reading the definition, the reader should satisfy themselves
that the arguments of this section apply equally to weakly piecewise linear
paths.

\section{Weakly piecewise linear paths}

\label{section7}Paths that lie in lines are special.

\begin{definition}
A continuous path $\gamma _{t}$ is weakly linear (geodesic) on $\left[ 0,T%
\right] $ if there is a line $l$ (or geodesic $l$)\ so that $\gamma _{t}\in
l $ for all $t\in \left[ 0,T\right] $.
\end{definition}

Suppose that $\gamma $ is smooth enough that one can form its iterated
integrals.

\begin{lemma}
If $\gamma $ is weakly linear, then the $n$-signature of the path $\gamma
\left( t\right) _{t\in \left[ 0,T\right] }$ is 
\begin{equation*}
\sum_{n=0}^{\infty }\frac{\left( \gamma _{T}-\gamma _{0}\right) ^{\otimes n}%
}{n!}.
\end{equation*}
In particular the signature of a weakly linear path is trivial if and only
if the path has $\gamma _{T}=\gamma _{0}$ or, equivalently, that it is a
loop.
\end{lemma}

\begin{lemma}
\label{lem:hfntree} A weakly geodesic, and in particular a weakly linear,
path with $\gamma _{0}=\gamma _{T}$ is always tree-like.
\end{lemma}

\begin{proof}
By definition, $\gamma $ lies in a single geodesic. Define $h\left( t\right)
=d\left( \gamma _{0},\gamma _{t}\right) .$ Clearly 
\begin{eqnarray*}
h\left( 0\right) &=&h\left( T\right) =0 \\
h &\geq &0.
\end{eqnarray*}
If $h\left( u\right) =0$ at some point $u\in \left( s,t\right) $ then 
\begin{eqnarray*}
d\left( \gamma _{s},\gamma _{t}\right) &\leq &d\left( \gamma _{0},\gamma
_{s}\right) +d\left( \gamma _{0},\gamma _{t}\right) \\
&=&h\left( s\right) +h\left( t\right) -2\inf_{u\in \left[ s,t\right]
}h\left( u\right)
\end{eqnarray*}
while if $h\left( u\right) >0$ at all points $u\in \left( s,t\right) $ then $%
\gamma _{s}$ and $\gamma _{t}$ are both on the same side of $\gamma _{0}$ in
the geodesic. Assume that $d\left( \gamma _{0},\gamma _{s}\right) \geq
d\left( \gamma _{0},\gamma _{t}\right) $, then 
\begin{eqnarray*}
d\left( \gamma _{s},\gamma _{t}\right) &=&d\left( \gamma _{0},\gamma
_{s}\right) -d\left( \gamma _{0},\gamma _{t}\right) \\
&=&h\left( s\right) -h\left( t\right) \\
&\leq &h\left( s\right) +h\left( t\right) -2\inf_{u\in \left[ s,t\right]
}h\left( u\right).
\end{eqnarray*}
as required.
\end{proof}

There are two key operations, splicing and excising, which preserve the
triviality of the signature and (because we will prove it is the same thing)
the tree-like property. However, the fact that excision of tree-like pieces
preserves the tree-like property will be a consequence of our work.

\begin{definition}
If $\gamma \in V$ is a path taking $\left[ 0,T\right] $ to the vector space $%
V$, $t\in \left[ 0,T\right] $ and $\tau $ is a second path in $V$, then the
insertion of $\tau $ into $\gamma $ at the time point $t$ is the
concatenation of paths 
\begin{equation*}
\gamma |_{\left[ 0,t\right] }\ast \tau \ast \gamma |_{\left[ t,T\right] }.
\end{equation*}
\end{definition}

\begin{definition}
If $\gamma \in V$ is a path on $\left[ 0,T\right] ,$with values in a vector
space $V$, and $\left[ s,t\right] \subset \left[ 0,T\right] $, then $\gamma $
with the segment $\left[ s,t\right] $ excised is 
\begin{equation*}
\gamma |_{\left[ 0,s\right] }\ast \gamma |_{\left[ t,T\right] }.
\end{equation*}
\end{definition}

\begin{remark}
\textrm{Note that these definitions make sense for paths in manifolds as
well as in the linear case, but in this case concatenation requires the
first path to finish where the second starts. We will use these operations
for paths on manifolds, but it will always be a requirement for insertion
that $\tau $ is a loop based at $\gamma _{t},$ for excision we require that $%
\gamma |_{\left[ s,t\right] }$ is a loop.}
\end{remark}

We have the following two easy lemmas:

\begin{lemma}
\label{lem:hfnexcision} Suppose that $\gamma \in M$ is a tree-like path in a
manifold $M$, and that $\tau $ is a tree-like path in $M$ that starts at $%
\gamma _{t},$ then the insertion of $\tau $ into $\gamma $ at the point $t$
is also tree-like. Moreover, the insertion at the time point $t$ of any
height function for $\tau $ into any height function coding $\gamma $ is a
height function for $\gamma |_{\left[ 0,t\right] }\ast \tau \ast \gamma |_{%
\left[ t,T\right] }$.
\end{lemma}

\begin{proof}
Assume $\gamma \in M$ is a tree-like path on a domain $\left[ 0,T\right] $,
by definition there is a positive and continuous function $h$ so that for
every $s$, $\tilde{s}$ in the domain $\left[ 0,T\right] $ 
\begin{eqnarray*}
d\left( \gamma _{s},\gamma _{\tilde{s}}\right) &\leq &h\left( s\right)
+h\left( \tilde{s}\right) -2\inf_{u\in \left[ s,\tilde{s}\right] }h\left(
u\right), \\
h\left( 0\right) &=&h\left( T\right) =0.
\end{eqnarray*}
In a similar way, let the domain of $\tau $ be $\left[ 0,R\right] $ and let $%
g$ be the height function for $\tau $%
\begin{eqnarray*}
d\left( \tau _{s},\tau _{\tilde{s}}\right) &\leq &g\left( s\right) +g\left( 
\tilde{s}\right) -2\inf_{u\in \left[ s,\tilde{s}\right] }g\left( u\right), \\
g\left( 0\right) &=&g\left( R\right) =0.
\end{eqnarray*}
Now insert $g$ in $h$ at $t$ and $\tau $ in $\gamma $ at $t$. Let $\tilde{h}$%
, $\tilde{\gamma}$ be the resulting functions defined on $\left[ 0,T+R\right]
$. Then 
\begin{equation*}
\begin{array}{cc}
\tilde{\gamma}\left( s\right) =\gamma \left( s\right), & 0\leq s\leq t \\ 
\tilde{\gamma}\left( s\right) =\tau \left( s-t\right), & t\leq s\leq t+R \\ 
\tilde{\gamma}\left( s\right) =\gamma \left( s-R\right), & t+R\leq s\leq T+R,%
\end{array}%
\end{equation*}
and 
\begin{equation*}
\begin{array}{cc}
\tilde{h}\left( s\right) =h\left( s\right), & 0\leq s\leq t \\ 
\tilde{h}\left( s\right) =g\left( s-t\right), & t\leq s\leq t+R \\ 
\tilde{h}\left( s\right) =h\left( s-R\right), & t+R\leq s\leq T+R,%
\end{array}%
\end{equation*}
where the definition of these functions for $s\in \left[ t+R,T+R\right] $
uses the fact that $\tau $ and $g$ are both loops.

Now it is quite obvious that if $s,\tilde{s}\in \left[ 0,T+R\right]
\backslash \left[ t,t+R\right] $, then 
\begin{eqnarray*}
d\left( \tilde{\gamma}_{s},\tilde{\gamma}_{\tilde{s}}\right) &\leq &\tilde{h}%
\left( s\right) +\tilde{h}\left( \tilde{s}\right) -2\inf_{u\in \left[ s,%
\tilde{s}\right] \backslash \left[ t,t+R\right] }\tilde{h}\left( u\right) \\
&\leq &\tilde{h}\left( s\right) +\tilde{h}\left( \tilde{s}\right)
-2\inf_{u\in \left[ s,\tilde{s}\right] }\tilde{h}\left( u\right) \\
\tilde{h}\left( 0\right) &=&\tilde{h}\left( T+R\right) =0
\end{eqnarray*}
and that for $s,\tilde{s}\in \left[ t,t+R\right] $, 
\begin{eqnarray*}
d\left( \tilde{\gamma}_{s},\tilde{\gamma}_{\tilde{s}}\right) &=&d\left( \tau
_{s-t},\tau _{\tilde{s}-t}\right) \\
&\leq &g\left( s-t\right) +g\left( \tilde{s}-t\right) -2\inf_{u\in \left[
s-t,\tilde{s}-t\right] }g\left( u\right) \\
&=&\tilde{h}\left( s\right) +\tilde{h}\left( \tilde{s}\right) -2\inf_{u\in %
\left[ s,\tilde{s}\right] }\tilde{h}\left( u\right).
\end{eqnarray*}
To finish the proof we must consider the case where $0\leq s\leq t\leq 
\tilde{s}\leq t+R$ and the case where $0\leq t\leq s\leq t+R\leq \tilde{s}%
\leq T+R$. As both cases are essentially identical we only deal with the
first. In this case 
\begin{eqnarray*}
d\left( \tilde{\gamma}_{s},\tilde{\gamma}_{\tilde{s}}\right) &=&d\left(
\gamma _{s},\tau _{\tilde{s}-t}\right) \\
&\leq &d\left( \gamma _{s},\gamma _{t}\right) +d\left( \tau _{0},\tau _{%
\tilde{s}-t}\right) \\
&\leq &h\left( s\right) +h\left( t\right) -2\inf_{u\in \left[ s,t\right]
}h\left( u\right) +g\left( \tilde{s}-t\right) -g\left( 0\right) \\
&=&\tilde{h}\left( s\right) +\tilde{h}\left( \tilde{s}\right) -2\inf_{u\in %
\left[ s,t\right] }\tilde{h}\left( u\right) \\
&\leq &\tilde{h}\left( s\right) +\tilde{h}\left( \tilde{s}\right)
-2\inf_{u\in \left[ s,\tilde{s}\right] }\tilde{h}\left( u\right) .
\end{eqnarray*}
\end{proof}

\begin{remark}
\textrm{The argument above is straightforward and could have been left to
the reader. However, we draw attention to the converse result, which also
seems very reasonable: that a tree-like path with a tree-like piece excised
is still tree-like. This result seems very much more difficult to prove. The
point is that the height function one has initially, as a consequence of $%
\gamma $ being tree-like, may well not certify that $\tau $ is tree-like
even though there is a second height function defined on $\left[ s,t\right] $
that certifies that it is. A direct proof that there is a new height
function simultaneously attesting to the tree-like nature of $\gamma $ and $%
\tau $ seems difficult. Using the full power of the results in the paper, we
can do this for paths of bounded variation.} 
\end{remark}

\begin{lemma}
\label{lem:excision} Let be $\gamma $ a path defined on $\left[ 0,T\right] $
with values in $V$ and suppose that $\gamma |_{\left[ s,t\right] }$ has
trivial signature where $\left[ s,t\right] \subset \left[ 0,T\right] $. Then 
$\gamma $ has trivial signature if and only if $\gamma $ with the segment $%
\left[ s,t\right] $ excised has trivial signature.
\end{lemma}

\begin{proof}
This is also easy. Since the signature map is a homomorphism we see that 
\begin{eqnarray*}
\gamma &=&\gamma |_{\left[ 0,s\right] }\ast \gamma |_{\left[ s,t\right]
}\ast \gamma |_{\left[ t,T\right] } \\
S\left( \gamma \right) &=&S\left( \gamma |_{\left[ 0,s\right] }\right)
\otimes S\left( \gamma |_{\left[ s,t\right] }\right) \otimes S\left( \gamma
|_{\left[ t,T\right] }\right)
\end{eqnarray*}%
and by hypothesis $S\left( \gamma |_{\left[ s,t\right] }\right) $ is the
identity in the tensor algebra. Therefore 
\begin{eqnarray*}
S\left( \gamma \right) &=&S\left( \gamma |_{\left[ 0,s\right] }\right)
\otimes S\left( \gamma |_{\left[ t,T\right] }\right) \\
&=&S\left( \gamma |_{\left[ 0,s\right] }\ast \gamma |_{\left[ t,T\right]
}\right) .
\end{eqnarray*}
\end{proof}

\begin{definition}
A continuous path $\gamma ,$ defined on $\left[ 0,T\right] $ is weakly
piecewise linear (or more generally, weakly geodesic) if there are finitely
many times 
\begin{equation*}
0=t_{0}<t_{1}<t_{2}<\ldots <t_{r}=T
\end{equation*}
such that for each $0<i\leq r$ the path segment $\gamma _{\left[
t_{i-1},t_{i}\right] }$ is weakly linear (geodesic).\footnote{%
The geodesic will always be unique since the path has unit speed and $%
t_{i}<t_{i+1}$ so contains at least two distinct points.}
\end{definition}

Our goal in this section is to prove, through an induction, that a weakly
linear path with trivial signature is tree-like and construct the height
function. As before, every such path admits a unique partition so that

\begin{lemma}
If $\gamma $ is a weakly piecewise linear path, then there exists a unique
partition $0=t_{0}<t_{1}<t_{2}<\ldots <t_{r}=T$ so that the linear segments
associated to $\left[ \gamma _{t_{i-1}},\gamma _{t_{i}}\right] $ and $\left[
\gamma _{t_{i}},\gamma _{t_{i+1}}\right] $ are not collinear for any $0<i<r$.
\end{lemma}

We will henceforth only use this partition and refer to $r$ as the number of
segments in $\gamma $.

\begin{definition}
We say $\gamma $ is fully non-degenerate if, in addition, $\gamma
_{t_{i-1}}\neq \gamma _{t_{i}}$ for every $0<i\leq r$.
\end{definition}

\begin{lemma}
\label{lem:wlp2} If $\gamma $ is a weakly linear path with trivial signature
and at least one segment, then there exist $0<i\leq r$ so that $\gamma
_{t_{i-1}}=\gamma _{t_{i}}$ .
\end{lemma}

\begin{proof}
The arguments in the previous section on piecewise linear paths apply
equally to weakly piecewise linear and weakly piecewise geodesic paths. In
particular Corollary~\ref{spaced_out} only refers to the location of $\gamma 
$ at the times $t_{i}$ at which the path changes direction (by an angle
different from $\pi $).
\end{proof}

\begin{proposition}
\label{prop:main} Any weakly piecewise linear path $\gamma $ with trivial
signature is tree-like with a height function whose total variation is the
same as that of $\gamma $.
\end{proposition}

\begin{proof}
The argument is a simple induction using the lemmas above. If it has no
segments we are clearly finished with $h\equiv 0$. We now assume that any
weakly piecewise linear path $\gamma^{(r-1)}$, consisting of at most $r-1$
segments, with trivial signature is tree-like with a height function whose
total variation is the same as that of $\gamma^{(r-1)}$. Suppose that $%
\gamma^{(r)}$ is chosen so that it is a weakly piecewise linear path of $r$
segments with trivial signature but there was no height function coding it
as a tree-like path with total variation controlled by that of $\gamma^{(r)}$%
. Then, by Lemma~\ref{lem:wlp2}, in the standard partition there must be $%
0<i\leq r$ so that $\gamma^{(r)} _{t_{i-1}}=\gamma^{(r)}_{t_{i}}$, and by
assumption $t_{i-1}<t_{i}$. In other words, the segment $\gamma^{(r)} |_{%
\left[ t_{i-1},t_{i}\right] }$ is a weakly linear segment and a loop. It
therefore has trivial signature, is tree-like and the height function we
constructed for it in the proof of Lemma~\ref{lem:hfntree} was indeed
controlled by the variation of the loop.

Let $\hat{\gamma}$ be the result of excising the segment $\gamma |_{\left[
t_{i-1},t_{i}\right] }$ from $\gamma^{(r)} $. As $\gamma^{(r)} |_{\left[
t_{i-1},t_{i}\right] }$ has trivial signature, by Lemma~\ref{lem:excision}, $%
\hat{\gamma}$ also has trivial signature. On the other hand, $\hat{\gamma}$
is weakly piecewise linear with fewer edges than $\gamma $ (it is possible
that $\gamma $ restricted to $\left[ t_{i-2},t_{i-1}\right] $ and $\left[
t_{i},t_{i+1}\right] $ are collinear and so the number of edges drops by
more than one in the canonical partition - but it will always drop!). So by
induction, $\hat{\gamma}$ is tree-like and is controlled by some height
function $\hat{h}$ that has total variation controlled by the variation of $%
\hat{\gamma}$.

Now insert the tree-like path $\gamma ^{(r)}|_{\left[ t_{i-1},t_{i}\right] }$
into $\hat{\gamma}$. By Lemma~\ref{lem:hfnexcision} this will be tree-like
and the height function is simply the insertion of the height function for $%
\gamma ^{(r)}|_{\left[ t_{i-1},t_{i}\right] }$ into that for $\hat{\gamma}$
and by construction is indeed controlled by the variation of $\gamma ^{(r)}$
as required. Thus we have completed our induction.
\end{proof}

\section{Proof of the main theorem}

\label{section8}We can now combine the results of the last sections to
conclude the proof of our main theorem and its corollaries.

\begin{proof}[Proof of Theorem~\protect\ref{thm:main}]
Corollary~\ref{cor:half} establishes that tree-like paths have trivial
signature.

Thus we only need to establish that if the path of bounded variation has
trivial signature, then it is tree-like. By Lemma~\ref{lem:oneformrep} we
can write the path as an integral against a rank one 1-form. By Corollary~%
\ref{cor:1formapprox} we can approximate any rank one 1-form by a sequence
of rank one 1-forms with the property that each 1-form is piecewise constant
on finitely many disjoint compact sets and 0 elsewhere. By integrating $%
\gamma $ against the sequence of 1-forms we can construct a sequence of
weakly piecewise linear paths approximating $\gamma $ in bounded variation.
By Corollary~\ref{cor:1formtriv}, these approximations have trivial
signature. By Proposition~\ref{prop:main} this means that these weakly
piecewise linear paths must be tree-like. Hence we have a sequence of
tree-like paths which approximate $\gamma $. By re-parametrizing the paths
at unit speed and using Lemma~\ref{lem:cmpctness} $\gamma $ must be
tree-like, completing the proof.
\end{proof}

\begin{proof}[Proof of Corollary~\protect\ref{cor:equiv}]
Recall that we defined $X\sim Y$, by the relation that $X$ then $Y$ run
backwards is tree-like. The transitivity is the part that is not obvious.
However, we can now say $X\sim Y$ if and only if the signature of $\mathbf{X}%
\mathbf{Y}^{-1}$ is trivial. As multiplication in the tensor algebra is
associative, it is now simple to check the conditions for an equivalence
relation. Denoting the signature of $X$ by $\mathbf{X}$ etc. one sees that 
\newline
1. The path run backward has signature $\mathbf{Y}\mathbf{X}^{-1}=-\mathbf{X}%
\mathbf{Y}^{-1}= \mathbf{0}$.\newline
2. $\mathbf{X}\mathbf{X}^{-1}=\mathbf{0}$ by definition.\newline
3. If $X\sim Y$ and $Y\sim Z$, then $\mathbf{XY}^{-1}=\mathbf{0}$ and $%
\mathbf{YZ}^{-1}=\mathbf{0}$. Thus 
\begin{equation*}
\mathbf{0} = \left( \mathbf{XY}^{-1}\right) \left( \mathbf{YZ}^{-1}\right) =%
\mathbf{X}\left( \mathbf{Y}^{-1}\mathbf{Y}\right) \mathbf{Z}^{-1} = \mathbf{X%
} \left( \mathbf{0}\right) \mathbf{Z}^{-1} =\mathbf{XZ}^{-1}
\end{equation*}
and hence $X\sim Z$ as required.

It is straightforward to see that the equivalence classes form a group.
\end{proof}

\begin{proof}[Proof of Corollary~\protect\ref{cor:minimisers}]
In order to deduce the existence and uniqueness of minimisers for the length
within each equivalence class we observe that;\newline
1. We can re-parameterise the paths to have unit speed and thereafter to be
constant. Then by the compactness of the equivalence classes of paths with
the same signature and bounded length, any sequence of paths will have a
subsequential uniform limit with the same signature. As length is
lower-semicontinuous in the uniform topology, the limit of a sequence of
paths with length decreasing to the minimum will have length less than or
equal to the minimum. We have seen, through a subsubsequence where the
height functions also converge, that it will also be in the same equivalence
class as far as the signature is concerned, so it is a minimiser.\newline
2. Within the class of paths with given signature and finite length there
will always be \emph{at least one} minimal element. Let $X$ and $Y$ be two
minimisers parameterised at unit speed, and let $h$ be a height function for 
$XY^{-1}$. Let the time interval on which $h$ is defined be $\left[ 0,T%
\right] $ and let $\tau $ denote the time at which the switch from $X$ to $Y$
occurs. The function $h$ is monotone on $\left[ 0,\tau \right] $ and on $%
\left[ \tau ,T\right] $ for otherwise there would be an interval $\left[ s,t%
\right] \subset \left[ 0,\tau \right] $ with $h\left( s\right) =h\left(
t\right) $. Then the function $u\rightarrow h\left( u\right) -h\left(
s\right) $ is a height function confirming that the restriction of $X$ to $%
\left[ s,t\right] $ is tree-like. Now we know from the associativity of the
product in the tensor algebra that the signature is not changed by excision
of a tree-like piece. Therefore, $X$ with the interval $\left[ s,t\right] $
excised is in the same equivalence class as $X$ but has strictly shorter
length. Thus $X$ could not have been a minimiser - as it is, we deduce the
function $h$ is strictly monotone. A similar argument works on $\left[ \tau
,T\right] .$

Let $\sigma :\left[ 0,\tau \right] \rightarrow \left[ \tau ,T\right] $ be
the unique function with 
\begin{equation*}
h\left( t\right) =h\left( \sigma \left( t\right) \right).
\end{equation*}
Then $\sigma $ is continuous decreasing and $\sigma \left( 0\right) =T$ and $%
\sigma \left( \tau \right) =\tau $. Moreover, $X_{u}=Y_{T-\sigma \left(
u\right) }$ and so we see that $($up to reparameterisations), the two paths
are the same.

Hence we have a unique minimal element!
\end{proof}

\appendix
\section{}

\label{appendix1}

\subsection{Trees and paths - background information}

We have shown in this paper that trees have an important role as the
negligible sets of control theory, quite analogous to the null sets of
Lebesgue integration. The trees we need to consider are \emph{analytic}
objects in flavour, and not the finite combinatorial objects of
undergraduate courses. 
In this appendix we collect together a few related ways of looking at them,
and prove a basic characterisation generalising the concept of height
function.

We first recall that\\
(1). Graphs $\left( E,V\right) $ that are acyclic and connected are
generally called \emph{trees}. If such a tree is non-empty and has a
distinguished vertex $\mathbf{v}$ it is called a \emph{rooted tree}.\\
(2). A rooted tree induces and is characterised by a \emph{partial order on 
}$V$\emph{\ with least element }$\mathbf{v}$. The partial order is defined
as follows 
\begin{equation*}
a\preceq b\quad\text{if the circuit free path from the root } \mathbf{v}%
\rightarrow b\text{ goes through }a.
\end{equation*}
This order has the property that for each fixed $b$ \emph{the set }$\left\{
a\preceq b\right\} $\emph{\ is totally ordered} by $\preceq$.

Conversely any partial order on a finite set $V$ with a least element $v$
and the property that for each $b$ the set $\left\{ a\preceq b\right\} $ is
totally ordered defines a unique rooted tree on $V$. One of the simplest
ways to construct a tree is to consider a (finite) collection $\Omega$ of
paths in a graph with all paths starting at a fixed vertex, and with the
partial order that $\omega\preceq\omega^{\prime}$ iff $\omega$ is an initial
segment of $\omega^{\prime}.$\\
(3). Alternatively, let $\left( E,V\right) $ be a graph extended into a
continuum by assigning a length to each edge. Let $d\left( a,b\right) $ be
the infimum of the lengths of paths\footnote{%
the sum of the lengths of the edges} between the two vertices $a,\;b$ in the
graph. Then $g$ is a geodesic metric on $V$. Trees are exactly \emph{the
graphs that give rise to 0-hyperbolic metrics} in the sense of Gromov (see
for example \cite{Kapovich95}).\\
(4). There are many ways to enumerate the edges and nodes of a finite
rooted tree. One way is to think of a family tree recording the descendants
of a single individual (the root). Start with the root. At the root, if all
children have been visited stop, at any other node, if all the children have
been visited, move up to the parent. If there are children who have not been
visited, then visit the oldest unvisited child. At each time $n$ the
enumeration either moves up an edge or down an edge - each edge is visited
exactly twice. Let $h\left( n\right) $ denote the distance from the top of
the family tree after $n$ steps in this enumeration with the convention that 
$h\left( 0\right) =0$, then $h$ is similar to the path of a random walk,
moving up or down one unit at each step, except that it is positive and
returns to zero exactly as many times as there are edges coming from the
root. Hence $h\left( 2\left| E\right| \right) =0$.

\emph{The function }$h$\emph{\ completely describes the rooted tree.} The
function $h$ directly yields the nearest neighbour metric on the tree. If $h$
is a function such that $h\left( 0\right) =0$, it moves up or down one unit
at each step, is positive and $h\left( 2\left| E\right| \right) =0$, then $d$
defined by 
\begin{equation*}
d\left( m,n\right) =h\left( m\right) +h\left( n\right) -2\inf _{u\in\left[
m,n\right] }h\left( u\right),
\end{equation*}
is a pseudo-metric on $\left[ 0,2\left| V\right| \right] $. If we identify
points in $\left[ 0,2\left| V\right| \right] $ that are zero distance apart
and join by edges the equivalence classes of points that are distance one
apart, then one recovers an equivalent rooted tree.

Put less pedantically, let the enumeration be $a$ at step $n$ and $b$ at
step $m$ and define 
\begin{equation*}
d\left( a,b\right) =h\left( m\right) +h\left( n\right) -2\inf_{u\in \left[
m,n\right] }h\left( u\right) ,
\end{equation*}
then it is simple to check that $d$ is well defined and is a metric on
vertices making the set of vertices a tree.

Thus excursions of simple (random) walks are a convenient (and well studied)
way to describe abstract graphical trees. This particular choice for \emph{%
coding a tree with a positive function on the interval} can be extended to
describe continuous trees. This approach was used by
Le Gall \cite{snake} in his development of the Brownian snake associated to
the measure valued Dawson-Watanabe process. 

\subsection{$\mathbf{R}$-trees are coded by continuous functions}

One of the early examples of a continuous tree is the evolution of a
continuous time stochastic process, where, as is customary in probability,
one identifies the evolution of two trajectories until the first time they
separate. (This idea dates back at least to Kolmogorov and his introduction
of filtrations). Another popular and equivalent approach to continuous trees
is through $\mathbf{R}$-trees (\cite{rtree} p425 and the references there).

Interestingly, analysts and probabilists have generally rejected the
abstract tree as too wild an object, and usually add extra structure,
essentially a second topology or Borel structure on the tree that comes from
thinking of the tree as a family of paths in a space which also has some
topology. This approach is critical to the arguments used here, as we prove
our tree-like paths are tree-like by approximating them with simpler
tree-like paths. (They would never converge in the `hyperbolic' metric). In
contrast, group theorists and low dimensional topologists have made a great
deal of progress by studying specific symmetry groups of these trees and do
not seem to find their hugeness too problematic.

Our goal in this subsection of the appendix is to prove the simple
representation: that the general $\mathbf{R}$-tree arises from identifying
the contours of a continuous function on a locally connected and connected
space. The height functions we considered on $\left[ 0,T\right] $ are a
special case.

\begin{definition}
An $\mathbf{R}$-tree is a uniquely arcwise connected metric space, in which
the arc between two points is isometric to an interval.
\end{definition}

Such a space is locally connected, for let $B_{x}$ be the set of points a
distance at most $1/n$ from $x$. If $z\in B_{x}$, then the arc connecting $x$
with $z$ is isometrically embedded, and hence is contained in $B_{x}$. Hence 
$B_{x}$ is the union of connected sets with non-empty common intersection
(they contain $x$) and is connected. The sets $B_{x}$ form a basis for the
topology induced by the metric. Observe that if two arcs meet at two points,
then the uniqueness assertion ensures that they coincide on the interval in
between.

Fix some point $v$ as the `root' and let $x$ and $y$ be two points in the $%
\mathbf{R}$-tree. The arcs from $x$ and $y$ to $v$ have a maximal interval
in common starting at $v$ and terminating at some $v_{1}$, after that time
they never meet again. One arc between them is the join of the arcs from $x$
to $v_{1}$ to $y$ (and hence it is the arc and a geodesic between them).
Hence 
\begin{equation*}
d\left( x,y\right) =d\left( x,v\right) +d\left( y,v\right) -2d\left(
v,v_{1}\right) .
\end{equation*}

\begin{example}
Consider the space $\Omega $ of continuous paths $X_{t}\in E$ where each
path is defined on an interval $\left[ 0,\xi \left( \omega \right) \right) $
and has a left limit at $\left[ 0,\xi \left( \omega \right) \right) $.
Suppose that if $X\in \Omega $ is defined on $\left[ 0,\xi \right) $, then $%
X|_{\left[ 0,s\right) }\in \Omega $ for every $s$ less than $\xi $. Define 
\begin{equation*}
d\left( \omega ,\omega ^{\prime }\right) =\xi \left( \omega \right) +\xi
\left( \omega ^{\prime }\right) -2\sup \left\{ t<\min \left( \xi \left(
\omega \right) ,\xi \left( \omega ^{\prime }\right) \right) |\;\omega \left(
s\right) =\omega ^{\prime }\left( s\right) \;\forall s\leq t\right\} .
\end{equation*}
Then $\left( \Omega ,d\right) $ is an $R$-tree.
\end{example}

We now give a way of constructing $\mathbf{R}$-trees. 
The basic idea for this is quite easy, but the core of the argument lies in
the detail so we proceed carefully in stages.

Let $I$ be a connected and locally connected topological space, and $%
h:I\rightarrow\mathbb{R}$ be a positive continuous function that attains its
lower bound at a point $v\in I$.

\begin{definition}
For each $x\in I$ and $\lambda\leq h\left( x\right) $ define $C_{x,\lambda}$
to be the maximal connected subset of $\left\{ y\;|\;h\left( y\right)
\geq\lambda\right\} $ containing $x$.
\end{definition}

\begin{lemma}
\label{lem:cxyexist} The sets $C_{x,\lambda}$ exist, and are closed.
Moreover, if $C_{x,\lambda }\cap C_{x^{\prime},\lambda^{\prime}}\neq\phi$
and $\lambda\leq\lambda ^{\prime}$, then 
\begin{equation*}
C_{x^{\prime},\lambda^{\prime}}\subset C_{x,\lambda}.
\end{equation*}
\end{lemma}

\begin{proof}
An arbitrary union of connected sets with non-empty intersection is
connected, taking the union of all connected subsets of $\left\{
y\;|\;h\left( y\right) \geq\lambda\right\} $ containing $x$ constructs the
unique maximal connected subset. Since $h$ is continuous the closure $%
D_{x,\lambda}$ of $C_{x,\lambda}$ is also a subset of $\left\{ y\;|\;h\left(
y\right) \geq\lambda\right\}$. The closure of a connected set is always
connected hence $D_{x,\lambda}$ is also connected. It follows from the fact
that $C_{x,\lambda}$ is maximal that $C_{x,\lambda}=D_{x,\lambda}$ and so is
a closed set.

If $C_{x,\lambda}\cap C_{x^{\prime},\lambda^{\prime}}\neq\phi$ and $%
\lambda\leq\lambda^{\prime}$, then 
\begin{equation*}
x\in C_{x,\lambda}\cup C_{x^{\prime},\lambda^{\prime}}\subset\left\{
y\;|\;h\left( y\right) \geq\lambda\right\},
\end{equation*}
and since $C_{x,\lambda}\cap C_{x^{\prime},\lambda^{\prime}}\neq\phi$, the
set $C_{x,\lambda}\cup C_{x^{\prime},\lambda^{\prime}}$ is connected. Hence
maximality ensures $C_{x,\lambda}=C_{x,\lambda}\cup
C_{x^{\prime},\lambda^{\prime}}$ and hence $C_{x^{\prime},\lambda^{\prime}}%
\subset C_{x,\lambda}$.
\end{proof}

\begin{corollary}
Either $C_{x,\lambda}$ equals $C_{x^{\prime},\lambda}$ or it is disjoint
from it.
\end{corollary}

\begin{proof}
If they are not disjoint, then the previous Lemma can be applied twice to
prove that $C_{x^{\prime},\lambda}\subset C_{x,\lambda}$ and $%
C_{x,\lambda}\subset C_{x^{\prime},\lambda}.$
\end{proof}

\begin{corollary}
If $C_{x,\lambda}=C_{x^{\prime},\lambda}$, then $C_{x,\lambda^{\prime%
\prime}}=C_{x^{\prime},\lambda^{\prime\prime}}$ for all $\lambda^{\prime%
\prime }<\lambda$.
\end{corollary}

\begin{proof}
The set $C_{x,\lambda},C_{x^{\prime},\lambda}$ are nonempty and have
nontrivial intersection. $C_{x,\lambda}\subset C_{x,\lambda^{\prime\prime}}$
and $C_{x^{\prime},\lambda}\subset C_{x^{\prime},\lambda^{\prime\prime}}$
hence $C_{x,\lambda^{\prime\prime}}$ and $C_{x^{\prime},\lambda^{\prime%
\prime }}$ have nontrivial intersection. Hence they are equal.
\end{proof}

\begin{corollary}
$y\in C_{x,\lambda}$ if and only if $C_{y,h\left( y\right) }\subset
C_{x,\lambda}$.
\end{corollary}

\begin{proof}
Suppose that $y\in C_{x,\lambda}$, then $C_{y,h\left( y\right) }$ and $%
C_{x,\lambda}$ are not disjoint. It follows from the definition of $%
C_{x,\lambda}$ and $y\in C_{x,\lambda}$ that $h\left( y\right) \geq\lambda$.
By Lemma~\ref{lem:cxyexist} $C_{y,h\left( y\right) }\subset C_{x,\lambda}$.
Suppose that $C_{y,h\left( y\right) }\subset C_{x,\lambda}$, since $y\in
C_{y,h\left( y\right) }$ it is obvious that $y\in C_{x,\lambda}$.
\end{proof}


\begin{definition}
The set $C_{x}:=C_{x,h\left( x\right) }$ is commonly referred to as the 
\emph{contour} of $h$ through $x$.
\end{definition}

The map $x\rightarrow C_{x}$ induces a partial order on $I$ with $x\preceq y$
if $C_{x}\supseteq C_{y}$. If $h$ attains its lower bound at $x$, then $%
C_{x}=I$ since $\left\{ y\;|\;h\left( y\right) \geq h\left( x\right)
\right\} =I$ and $I$ is connected by hypothesis. Hence the root $v\preceq y$
for all $y\in I$.

\begin{lemma}
Suppose that $\lambda\in\left[ h\left( v\right) ,h\left( x\right) \right]$,
then there is a $y$ in $C_{x,\lambda}$ such that $h\left( y\right) =\lambda$
and, in particular, there is always a contour ($C_{x,\lambda }$) at height $%
\lambda$ through $y$ that contains $x$.
\end{lemma}

\begin{proof}
By the definition of $C_{x,\lambda}$ it is the maximal connected subset of $%
h\geq\lambda$ containing $x$; assume the hypothesis that there is no $y$ in $%
C_{x,\lambda}$ with $h\left( y\right) =\lambda$ so that it is contained in $%
h>\lambda,$ hence C$_{x,\lambda}$ is a maximal connected subset of $%
h>\lambda $. Now $h>\lambda$ is open and locally connected, hence its
maximal connected subsets of $h>\lambda$ are open and $C_{x,\lambda}$ is
open. However it is also closed, which contradicts the connectedness of the $%
I$. Thus we have established the existence of the point $y$.
\end{proof}

The contour is obviously unique, although $y$ is in general not. If we
consider the equivalence classes $x\symbol{126}y$ if $x\preceq y$ and $%
y\preceq x$, then we see that the equivalence classes $\left[ y\right] _{%
\symbol{126}}$ of $y\preceq x$ are totally ordered and in one to one
correspondence with points in the interval $\left[ h\left( v\right) ,h\left(
x\right) \right]$.

\begin{lemma}
\label{lem:loccon} If $z\in C_{y,\lambda}$ and $h\left( z\right) >\lambda$,
then $z$ is in the interior of $C_{y,\lambda}$. If $C_{x^{\prime},\lambda^{%
\prime}}\subset C_{x,\lambda}$ with $\lambda^{\prime}>\lambda$, then $%
C_{x,\lambda}$ is a neighbourhood of $C_{x^{\prime},\lambda^{\prime}}$.
\end{lemma}

\begin{proof}
$I$ is locally connected, and $h$ is continuous, hence there is a connected
neighbourhood $U$ of $z$ such that $h\left( z\right) \geq\lambda$. By
maximality $U\subset C_{z,\lambda}$. Since $C_{z,\lambda}\cap
C_{y,\lambda}\neq\phi$ we have $C_{z,\lambda}=C_{y,\lambda} $ and thus $%
U\subset C_{y,\lambda}$. Hence $C_{y,\lambda}$ is a neighbourhood of $z$.
The last part follows trivially once by noting that for all $z\in
C_{x^{\prime },\lambda^{\prime}}$ we have $h\left( z\right)
\geq\lambda^{\prime}>\lambda$ and hence $C_{y,\lambda}$ is a neighbourhood
of $z$.
\end{proof}

We now define a pseudo-metric on $I$. Lemma~\ref{lem:loccon} (the only place
we will use local connectedness) is critical to showing that the map from $I$
to the resulting quotient space is continuous.

\begin{definition}
If $y$ and $z$ are points in $I$, define $\lambda\left( y,z\right) \leq
\min\left( h\left( y\right) ,h\left( z\right) \right) $ such that $%
C_{y,\lambda}=C_{z,\lambda}$%
\begin{equation*}
\lambda\left( y,z\right) =\sup\left\{ \lambda\;|\;C_{y,\lambda
}=C_{z,\lambda},\;\lambda\leq h\left( y\right) ,\;\lambda\leq h\left(
z\right) \right\} .
\end{equation*}
\end{definition}

The set 
\begin{equation*}
\left\{ \lambda\;|\;C_{y,\lambda}=C_{z,\lambda},\;\lambda\leq h\left(
y\right) ,\;\lambda\leq h\left( z\right) \right\}
\end{equation*}
is a non-empty interval $\left[ h\left( v\right) ,\lambda\left( y,z\right) %
\right] $ or $\left[ h\left( v\right) ,\lambda\left( y,z\right) \right) $
where $\lambda\left( y,z\right) $ satisfies 
\begin{equation*}
h\left( v\right) \leq\lambda\left( y,z\right) \leq\min\left( h\left(
y\right) ,h\left( z\right) \right) .
\end{equation*}
Clearly $\lambda\left( x,x\right) =h\left( x\right) .$

\begin{lemma}
The function $\lambda $ is lower semi-continuous 
\begin{equation*}
\liminf_{z\rightarrow z_{0}}\lambda \left( y,z\right) \geq \lambda \left(
y,z_{0}\right) .
\end{equation*}
\end{lemma}

\begin{proof}
Fix $y,\;z_{0}$ and choose some $\lambda^{\prime}<\lambda\left(
y,z_{0}\right) $. By the definition of $\lambda(y,z_0)$ we have that $%
C_{y,\lambda^{\prime}}=C_{z_{0},\lambda^{\prime}}$. Since $h\left(
z_{0}\right) \geq\lambda^{\prime}$ there is a neighbourhood $U$ of $z_{0}$
so that $U\subset C_{z_{0},\lambda^{\prime}}$. For any $z\in U$ one has $%
z\in C_{z,\lambda^{\prime}}\cap C_{z_{0},\lambda^{\prime}}$. Hence $%
C_{z_{0},\lambda^{\prime}}=C_{z,\lambda^{\prime}}$ and $C_{y,\lambda^{%
\prime}}=C_{z,\lambda^{\prime}}$. Thus $\lambda\left( y,z\right)
\geq\lambda^{\prime}$ for $z\in U$ and hence 
\begin{equation*}
\liminf_{z\rightarrow z_{0}}\lambda\left( y,z\right) \geq\lambda^{\prime}.
\end{equation*}
Since $\lambda^{\prime}<\lambda\left( y,z_{0}\right) $ was arbitrary 
\begin{equation*}
\liminf_{z\rightarrow z_{0}}\lambda\left( y,z\right) \geq\lambda(y,z_0)
\end{equation*}
and the result is proved.
\end{proof}

\begin{lemma}
The following inequality holds 
\begin{equation*}
\min\left\{ \lambda\left( x,z\right) ,\lambda\left( y,z\right) \right\}
\leq\lambda\left( x,y\right) .
\end{equation*}
\end{lemma}

\begin{proof}
If $\min\left\{ \lambda\left( x,z\right) ,\lambda\left( y,z\right) \right\}
=h\left( v\right)$, then there is nothing to prove. Recall that 
\begin{equation*}
\left\{ \lambda\;|\;C_{y,\lambda}=C_{z,\lambda},\;\lambda\leq h\left(
y\right) ,\;\lambda\leq h\left( x\right) \right\}
\end{equation*}
is connected and contains $h\left( v\right)$. Suppose $h\left( v\right)
\leq\lambda<\min\left\{ \lambda\left( x,z\right) ,\lambda\left( y,z\right)
\right\}$, then it follows that the identity $C_{x,\lambda }=C_{z,\lambda}$
holds for $\lambda$. Similarly $C_{y,\lambda}=C_{z,\lambda}$. As a result $%
C_{x,\lambda}=C_{y,\lambda}$ and $\lambda\left( x,y\right) \geq\lambda$.
\end{proof}

\begin{definition}
Define $d$ on $I\times I$ by 
\begin{equation*}
d\left( x,y\right) =h\left( x\right) +h\left( y\right) -2\lambda\left(
x,y\right).
\end{equation*}
\end{definition}

\begin{lemma}
The function $d$ is a pseudo-metric on $I$. If $\left( \tilde{I},d\right) $
is the resulting quotient metric space, then the projection $I\rightarrow 
\tilde{I}$ from the topological space $I$ to the metric space is continuous.
\end{lemma}

\begin{proof}
Clearly $d$ is positive, symmetric and we have remarked that for all $x$, $%
\lambda\left( x,x\right) =h\left( x\right) $ hence it is zero on the
diagonal. To see the triangle inequality, assume 
\begin{equation*}
\lambda\left( x,z\right) =\min\left\{ \lambda\left( x,z\right)
,\lambda\left( y,z\right) \right\}
\end{equation*}
and then observe 
\begin{align*}
d\left( x,y\right) & =h\left( x\right) +h\left( y\right) -2\lambda\left(
x,y\right) \\
& \leq h\left( x\right) +h\left( y\right) -2\lambda\left( x,z\right) \\
& =h\left( x\right) +h\left( z\right) -2\lambda\left( x,z\right) +h\left(
y\right) -h\left( z\right) \\
& \leq d\left( x,z\right) +\left| h\left( y\right) -h\left( z\right) \right|
\end{align*}
but $\lambda\left( y,z\right) \leq\min\left( h\left( y\right) ,h\left(
z\right) \right) $ and hence 
\begin{align*}
\left| h\left( y\right) -h\left( z\right) \right| & =h\left( y\right)
+h\left( z\right) -2\min\left( h\left( y\right) ,h\left( z\right) \right) \\
& \leq h\left( y\right) +h\left( z\right) -2\lambda\left( y,z\right) \\
& =d\left( y,z\right)
\end{align*}
hence 
\begin{equation*}
d\left( x,y\right) \leq d\left( x,z\right) +d\left( y,z\right)
\end{equation*}
as required.
\end{proof}

We can now introduce the equivalence relation $x\symbol{126}y$ if $d\left(
x,y\right) =0$ and the quotient space $I/\symbol{126}$. We write $I/\symbol{%
126}=\tilde{I}$ and $i:I\rightarrow\tilde{I}$ for the canonical projection.
The function $d$ projects onto $\tilde{I}\times\tilde{I}$ and is a metric
there.

It is tempting to think that $x\symbol{126}y$ if and only if $C_{x}=C_{y}$
and this is true if $I$ is compact Hausdorff. However the definitions imply
a slightly different criteria: $x\symbol{126}y$ iff 
\begin{equation*}
h\left( x\right) =h\left( y\right) =\lambda\text{ and }C_{x,\lambda
^{\prime\prime}}=C_{y,\lambda^{\prime\prime}}\text{ for all }\lambda
^{\prime\prime}<\lambda.
\end{equation*}
The stronger statement $x\symbol{126}y$ if and only if $C_{x}=C_{y}$ is not
true for all continuous functions $h$ on $\mathbb{R}^{2}$ as it is easy to
find a decreasing family of closed connected sets there whose limit is a
closed set that is not connected.

Consider again the new metric space $\tilde{I}$ that has as its points the
equivalence classes of points indistinguishable under $d$. We now prove that
the projection $i$ taking $I$ to $\tilde{I}$ is continuous. Fix $y\in I$ and 
$\varepsilon >0.$ Since $\lambda \left( y,.\right) $ is lower
semi-continuous and $h$ is (upper semi)continuous there is a neighbourhood $%
U $ of $y$ so that for $z\in U$ one has $\lambda \left( y,z\right) >\lambda
\left( y,y\right) -\varepsilon /4$ and $h\left( z\right) <h\left( y\right)
+\varepsilon /2$. Thus $d\left( y,z\right) <\varepsilon $ for $z\in U$.
Hence $\tilde{d}\left( i\left( y\right) ,i\left( z\right) \right)
<\varepsilon $ if $z\in U$. The function $i$ is continuous and as continuous
images of compact sets are compact we have the following.

\begin{corollary}
If $I$ is compact, then $\tilde{I}$ is a compact metric space.
\end{corollary}

To complete this section we will show $\tilde{I}$ is a uniquely arcwise
connected metric space, in which the arc between two points is isometric to
an interval and give a characterisation of compact trees.

\begin{proposition}
If $I$ is a connected and locally connected topological space, and $%
h:I\rightarrow \mathbb{R}$ is a positive continuous function that attains
its lower bound, then its \textquotedblleft contour tree\textquotedblright\
the metric space $\left( \tilde{I},\tilde{d}\right) $ is an $\mathbf{R}$%
-tree. Every $\mathbf{R}$ -tree can be constructed in this way.
\end{proposition}

\begin{proof}
It is enough to prove that the metric space $\tilde{I}$ we have constructed
is really an $\mathbf{R}$-tree and that every $\mathbf{R}$-tree can be
constructed in this way. Let $\tilde{x}$ any point in $\tilde{I}$ and $x\in I
$ satisfy $i\left( x\right) =\tilde{x}$. Then $h\left( x\right) $ does not
depend on the choice of $x$. Fix $h\left( v\right) <\lambda <h\left(
x\right) $. We have seen that there is a $y$ such that $h\left( y\right) $=$%
\lambda $ and $y\prec x$ moreover any two choices have the same contour
through them and hence the same $\tilde{y}\left( \lambda \right) $. In this
way we see that there is a map from $\left[ h\left( v\right) ,h\left(
x\right) \right] $ into $\tilde{I}$ that is injective. Moreover, it is
immediate from the definition of $d$ that it is an isometry and that $\tilde{%
I}$ is uniquely arc connected.

Suppose that $\Omega $ is an $\mathbf{R}$-tree, then we may fix a base
point, and for each point in the tree consider the distance from $V$ it is
clear that this continuous function is just appropriate to ensure that the
contour tree is the original tree.
\end{proof}

\begin{remark}
\textrm{1. In the case where $I$ is compact, obviously $\tilde{I}$ is both
complete and totally bounded as it is compact. }

\textrm{2. An $R$-tree is a metric space; it is therefore possible to
complete it. Indeed the completion consists of those paths, all of whose
initial segments are in the tree\footnote{\textrm{\textrm{We fix a root and
identify the tree with the geodesic arc from the root to the point in the
tree.}}}; we have not identified a simple sufficient condition on the
continuous function and topological space $\Omega $ to ensure this. An $R$%
-tree is totally bounded if it is bounded and for each $\varepsilon >0$
there is an $N$ so that for each $t$ the paths that extend a distance $t$
from the root have at most $N$ ancestral paths between them at time $%
t-\varepsilon $. In this way we see that the $R$-tree that comes out of
studying the historical process for the Fleming-Viot or the Dawson Watanabe
measure-valued processes is, with probability one, a compact $R$-tree for
each finite time.}
\end{remark}

\begin{lemma}
\bigskip Given a compact $R$-tree, there is always a height function on a
closed interval that yields the same tree as its quotient.
\end{lemma}

\begin{proof}
As the tree is compact, path connected and locally path connected, there is
always as based loop mapping $\left[ 0,1\right] $ onto the tree. Let $h$
denote the distance from the root. Its pullback onto the interval $\left[ 0,1%
\right] $ is a height function and the natural quotient is the original
tree. In this way we see that there is always a version of Le Gall's snake 
\cite{snake} traversing a compact tree.
\end{proof}



\end{document}